\documentclass{article}
\usepackage{epsfig}
\usepackage{eepic}
\usepackage{pict2e}
\usepackage{latexsym}
\usepackage{amsmath}
\usepackage{amssymb}
\usepackage[english]{babel}
\usepackage{hyperref}
\usepackage{xcolor}

\usepackage{geometry}
\newtheorem{thm}{Theorem}[subsection]

\newtheorem{lem}[thm]{Lemma}
\newtheorem{prop}[thm]{Proposition}
\newtheorem{defn}[thm]{Definition}
\newtheorem{rem}[thm]{Remark}
\newtheorem{expl}[thm]{Example}

\newcommand{\s}{\mathcal{S}}

\newcommand{\M}{\mathcal{M}}

\newcommand{\Real}{\mathbb{R}}

\DeclareMathOperator*{\bigboxplus}{\text{\Huge $_{\boxplus}$}}

\DeclareMathOperator*{\bigsmileplus}{\text{\Large ${\stackrel{_{+}}{\smile}}$}}
\DeclareMathOperator*{\bigsmileminus}{\text{\Large ${\stackrel{-}{\smile}}$}}

\DeclareMathOperator*{\Card}{\mathrm{Card}}
\DeclareMathOperator*{\Ls}{\mathrm{Ls}}
\DeclareMathOperator*{\Linf}{\mathrm{Li}}
\DeclareMathOperator*{\Lim}{\mathrm{Lim}}
\DeclareMathOperator*{\sgn}{\mathrm{sgn}}

\DeclareMathOperator*{\cl}{\mathrm{cl}}
\DeclareMathOperator*{\inter}{\mathrm{int}}

\begin{document}

\title{On the Multiary Algebraic Formulation of an Idempotent Symmetric Limit Convex Structure}
 \date{April 14, 2026}

\author{\thanks {LAMPS-Laboratory of Mathematics and Physics, University of Perpignan, 52 avenue Villeneuve,
66000 Perpignan, France.}  Walter Briec  }

\maketitle
\begin{abstract}

In \cite{bh}, $\mathbb{B}$-convexity was defined as an appropriate {Painlevé-Kuratowski} limit of linear convexities. More recently, an alternative algebraic formulation over the entire Euclidean vector space was proposed in \cite{b15} and \cite{b17}. The issue with the definition presented in \cite{bh} is that it was not developed from an algebraic perspective, but rather as the upper limit of a sequence of generalized convex polytopes whose form was not explicitly given. In this paper, we build on recent work and  provide an algebraic formulation for these limiting polytopes. Consequently, we deduce a  multiary algebraic form of $\mathbb{B}$-convexity that involves an idempotent, non-associative algebraic structure, extending the formalism proposed in \cite{b15} to an arbitrary number of points. Among other things, we demonstrate that these limiting polytopes do not satisfy the idempotent symmetrical convex structure defined in \cite{b15}. In the context of this formalism, we derive a general separation result in $\mathbb{R}^n$ by approximating convex sets by polytopes.  We conclude by clarifying some points regarding the external representation of polytopes proposed in \cite{b19} and analyze the structure of the polytopes that arise in this context.
\end{abstract}

{\bf AMS:} 06D50, 32F17\\


\section{Introduction}\label{SECMXVSP}

$\mathbb{B}$-convexity was introduced in ~\cite{bh} as a suitable  Painlevé-Kuratowski limit of linear convexities. These generalized convexities are derived from a $\varphi_p$-transformation defined as $\varphi_p(\lambda) = \lambda^{2p+1}$, following the scheme proposed in \cite{ben}. When restricted to the Euclidean orthant $\mathbb{R}_+^n$, $\mathbb{B}$-convexity\footnote{ The symbol $\mathbb{B}$, used to define a related notion of convexity, may be somewhat confusing given the concept of convexity introduced by Anatole Beck in Banach spaces \cite{b62}. Some clarifications regarding the choice of this notation are provided in \cite{bh, b24}. More importantly, note that the operation $\boxplus$, pronounced “$\mathbb{B}$oxplus,” visually resembles the Chebyshev unit ball (the $\ell_\infty$ ball) when depicted in two dimensions using Cartesian coordinates. This ball plays a central role in the proposed framework and arises as the Hausdorff{-Pompeiu} limit of the sequence of $\ell_p$ balls as $p \to \infty$.}  can be expressed in terms of a suitable lattice structure, reducing to Max-Times convexity. Max-Times convexity is homeomorphic to Max-Plus convexity, which arises in tropical mathematics (see references \cite{s88}, \cite{MP}, \cite{ms}, \cite{p94}). These convexities involve an idempotent semi-ring algebraic structure. In this context, $\mathbb{B}$-convex functions were analyzed in \cite{adil}, and Hahn-Banach-type separation properties were established in \cite{bh3}. Additionally, fixed point results were obtained in \cite{bh2} and \cite{by22}, and alternative idempotent convex structures were proposed in \cite{adilYe} and \cite{adilYeTi}.

It is well known, however, that symmetry, associativity, and idempotence are generally incompatible, except for trivial groups that reduce to the neutral element. In \cite{b15}, a special class of idempotent magmas was considered in which associativity was relaxed to preserve both symmetry and idempotence. This was achieved by introducing a non-associative operation denoted by $\mathbb{B}$ox-plus $\boxplus$, which benefits from an $n$-ary extension, for any $n$. It has been shown that this notion of convexity can be equivalently characterized using the  {Painlevé-Kuratowski} limit of the generalized convex hull of two points, as defined in \cite{bh}.

In this paper, we compare this definition to the one proposed in \cite{bh}, known as $\mathbb{B}$-convexity. A significant result of this study is the demonstration that an algebraic formulation of this convexity concept can be provided based on the non-associative and idempotent convex structure defined in \cite{b15}, whose algebraic properties have been more recently analyzed in detail in \cite{b20}. This implies a multiary extension of the binary notion of convexity proposed in \cite{b15}, formulated within a non-associative framework. It enables the direct treatment of sets of 
$m$ points, which cannot be derived from the binary case due to the non-associative nature of the underlying algebraic structure. 

We achieve this by employing a suitable notion of a determinant. Additionally, we show that $\mathbb{B}$-polytopes, defined as the limits of the generalized convex hulls of a finite number of points, are not convex in the sense of the convexity defined in \cite{b15} (i.e., for two points). 

With an algebraic definition of $\mathbb{B}$-convexity in the general case, it appears that this also imposes a stronger definition that enjoys many interesting properties. Paralleling some known properties    established   in the context of usual convexity, we deduce a limit separation theorem from some earlier results established in \cite{b17} and \cite{b20}.  {To achieve this objective, we rely on the fact that the separation of two $\mathbb B$-polytopes follows directly from the standard separation theorem for polytopes, upon taking a suitable limit.} Proceeding in this direction, we derive a separation result  by approximating convex sets by polytopes.

The paper unfolds as follows. Section 2 presents the framework of $\mathbb{B}$-convexity induced by the $\mathbb{B}$ox-plus operator $\boxplus$ and some related algebraic properties. Section 3 first establishes a preliminary inclusion property for $\mathbb{B}$-polytopes using a suitable notion of determinant and some earlier properties noted in \cite{b19}.  {Then, t}he converse inclusion is shown. We  deduce a  multiary algebraic form of $\mathbb{B}$-convexity that involves an idempotent, non-associative algebraic structure, extending the formalism proposed in \cite{b15} to an arbitrary number of points. Moreover, we provide examples demonstrating that that these limiting polytopes do not satisfy the idempotent symmetrical convex structure defined in \cite{b15}.  In Section 4, it is shown that compact $\mathbb{B}$-convex sets can be approximated by $\mathbb{B}$-polytopes. Along this line, we derive a separation result in Section 4. We conclude this section by clarifying the external representation of polytopes by a finite number of half-spaces, as considered in \cite{b19}.

\section{Idempotent and Non-Associative Convex Structure}

\subsection{Generalized Means and Convexity}
The concept of generalized means was introduced by Hardy, Littlewood, and Polya \cite{hlp}, and later extended into the field of optimization theory by Avriel \cite{avr1, avr2} and Ben-Tal \cite{ben}. In particular, Avriel introduced $r$-convex functions, which generalize the notion of log-convexity. This concept was further developed from an algebraic perspective by Ben-Tal \cite{ben}, who demonstrated that by transposing an algebraic structure via an isomorphism, a specific form of convexity can be derived within a vector space.\\
In ~\cite{bh}, $\mathbb{B}$-convexity is introduced as a limit of
linear convexities.
More precisely, for  {any} $p\in \mathbb{N}$,  a bijection $\varphi_p :\mathbb{R} \longrightarrow
\mathbb{R}$ defined by:
\begin{equation}\varphi_p:\lambda \longrightarrow \lambda^{2p+1}\end{equation}
 and $\phi_p (x_1,... ,x_n)=(\varphi_p(x_1),... ,\varphi_p(x_n))$ was considered; this is basically the approach of
Ben-Tal  \cite{ben} and Avriel   \cite{avr1}.
A field structure on
$\mathbb{R}$, for which $\varphi_p$ becomes a field isomorphism, can then be defined. Given this change of notation via
$\varphi_p$ and $\phi_p$ we can define  {an} $\Real$-vector space structure on $\Real^n$ by:
$\lambda\stackrel{\varphi_p}{\cdot} x=\phi_p^{-1}(\varphi_p(\lambda)\cdot\phi_p(x))=\lambda \cdot x$ for all $\lambda \in \Real$ and
$x  \stackrel{\varphi_p}{+} y =\phi_p^{-1}(\phi_p(x )+\phi_p(y))$, for all $x ,y\in \Real^n$. We call these two operations the 
 {indexed external product}  and the indexed sum (indexed by $\varphi_p$).

For all positive natural numbers $m$, let us denote $[m]=\{1,...,m\}$. The $\varphi_p$-sum, denoted $\stackrel{\varphi_p}{\sum}$,    of
$(x_1,...,x_m)\in \mathbb{R}^m$ is defined by
 \begin{equation}\stackrel{\varphi_p}{\sum_{i\in [m]}}x_i=
\phi_p^{-1}\Big(\sum_{i\in [m]}\phi_p (x_i)\Big)\end{equation} and
the
$\varphi_p$-convex hull of  a finite set $A=\{x^{(1)},...,x^{(m)}\}\subset \mathbb{R}^n$ is
defined by:
\begin{equation}Co^{\varphi_p}(A)=\Big\{\stackrel{\varphi_p}{\sum_{j\in [m]}}t_i\stackrel{\varphi_p}{\cdot}x^{(j)} :
\stackrel{\varphi_p}{\sum_{j\in [m]}}t_j=1, t_j\geq
0\Big\}\end{equation}
which can be rewritten:
\begin{equation}Co^{\varphi_p}(A)=\Big\{\phi_p^{-1}\Big(\sum_{j\in [m]}t_j^{2p+1}{\cdot }\phi_p\big(x^{(j)}\big)\Big) :
\big(\sum_{j \in [m]}t_j^{2p+1}\big)^{\frac{1}{2p+1}}=1 , j\in [m] , t_j\geq
0\Big\}.\end{equation}
It should be noted that a closely related notion of convexity has recently been analyzed in \cite{iy23}.

This   definition can be made   more explicit showing that:
\begin{equation}Co^{\varphi_p}(A)=\phi_p^{-1}\Big(Co\big(\phi_p(A)\big)\Big)\end{equation}
For simplicity, throughout the paper we denote for all $x,y\in \Real^n$:

\begin{equation}x\stackrel{p}{+}y=x\stackrel{\varphi_p}{+}y.\end{equation}
Moreover, for all $L\subset \mathbb{R}$, we simplify the notations
denoting $Co^{p}(L)=Co^{\varphi_p}(L)$. The Kuratowski-Painlev\'e\footnote{   { We have adopted the standard terminology “Painlevé–Kuratowski limits” throughout the manuscript. We also note that the terminology ``Peano--Kuratowski'' appears in the historical literature; see, for instance, the comments in \cite[Chapter~4]{rw98} and in \cite[Volume~I, Chapter~1, Comments]{t23}, where these terminological issues are discussed. }  } lower limit of the sequence of sets $\{A^{(p)}\}_{p\in \mathbb N}$,
denoted $\Linf_{p\to\infty} A^{(p)}$, is the set of points $x$ for which
there exists a sequence $\{x^{(p)}\}_{p\in \mathbb N}$ of points such that $x^{(p)}\in A^{(p)}$
for all $p$ and $x = \lim_{p\to\infty}x^{(p)}$; the Kuratowski-Painlev\'e upper limit of the sequence of sets $\{A^{(p)}\}_{p\in \mathbb N}$,
denoted $\Ls_{p\to\infty }A^{(p)}$, is the set of points $x$ for which
there exists a subsequence $\{x^{(p_k)}\}_{k\in \mathbb N}$ of points such that $x^{(p_k)}\in A^{(p_k)}$
for all $k$ and $x = \lim_{k\to\infty}x^{(p_k)}$ ; a sequence
$\{A^{(p_k)}\}_{k\in \mathbb{N}}$ of subsets of $\Real^n$ is said to
converge, in the Kuratowski-Painlev\'e sense, to a set $A$ if
$\Ls_{p\to\infty}A^{(p)} = A = \Linf_{p\to\infty}A^{(p)}$, in which case we
write $A = \Lim_{p\to\infty}A^{(p)}$.
We now define $\mathbb{B}$-convexity and its characterization in terms of
Kuratowski-Painlev\'e limit of $\varphi_p$-convex sets. For a non-empty finite subset $A=\{x^{(1)},...,x^{(m)} \}\subset \Real_{}^{n}$ we let
$Co^\infty(A)=\Ls_{p\longrightarrow \infty}Co^p(A)$. When $A\subset \Real_+^n$, it was shown in ~\cite{bh} that  we have\footnote{ {In the following, we denote $\Real_+=[0,+\infty[$ and $\Real_{++}=]0,+\infty[$.}}:
$$Co^\infty(A)=
\Big\{\bigvee_{j\in [m]}t_j x^{(j)}: t_j \in [0,1], \max_{j\in [m]} t_j =1\Big\}=\Lim_{p\longrightarrow \infty}Co^p(A),$$
where $\vee $ stands for the lattice operation related to the canonical partial order of $\Real^n.$ Moreover it is also proved
that the sequence $\{Co^p(A)\}_{p\in \mathbb N}$ is Hausdorff-{Pompeiu} convergent {(see  \cite {rw98}and \cite{t23} for definitions and terminologies)}. Such a formulation is useful to apprehend Max-Times convexity that is defined on the idempotent semi-module $(\Real_+^n, \vee, \cdot)$. This convexity is homeomorphic to Max-Plus convexity via a logarithmic transformation (see \cite{bh2} ).\\

$(i)$ A subset $C$  of $\Real^n$ is {\bf $\mathbb B$-convex} if for all finite subsets $A$ of $C$, $Co^\infty(A)\subset C$. \\

 $(ii)$ A subset $C$  of $\Real_+^n$ {\bf is Max-Times convex} ($\mathbb B$-convex over $\Real_+^n$) if for all  $x,y\in C$ and all $t\in [0,1]$   $x\vee ty\in C$. \\

$\mathbb{B}$-convex sets possess a number of remarkable properties. If $L \subset \Real_+^n$ is 
$\mathbb{B}$-convex (and therefore Max-Times convex), then 
$L = \lim_{p \to \infty} Co^p(L) = \bigcap_{p=0}^{\infty} Co^p(L). $ For any set $S \subset \Real_+^n$, we denote 
$Co^{\infty}(S) = \Lim_{p \to \infty} Co^p(S). $ Properties of this operation are developed in~\cite{bh}. In particular, for all $S \subset \Real_+^n$, the set $Co^{\infty}(S)$ is $\mathbb{B}$-convex.

\subsection{ An Algebraic Structure Extended to the whole Euclidean Vector Space}\label{Rec}

In \cite{b15} it was shown that for all $\lambda,\mu\in \Real$ we have:
 {\begin{equation}\label{base}\lim_{p\longrightarrow
+\infty}\lambda \stackrel{p}{+}\mu=
\left\{\begin{matrix}x\ &\hbox{ if } &|\lambda|&>&|\mu|\\
\frac{1}{2}(\lambda +\mu )&\hbox{ if }&|\lambda|&=&|\mu|\\
\mu& \hbox{ if }& |\lambda|&<&|\mu|.\end{matrix}\right.\end{equation}} Along
this line one can introduce the binary operation $\boxplus$
defined for all $\lambda,\mu\in \Real$ by:
 {\begin{equation}
\lambda \boxplus \mu=\lim_{p\longrightarrow +\infty}\lambda \stackrel{p}{+}\mu.
\end{equation}}

In the following it is established that, though the operation
$\boxplus$ does not satisfy associativity,  it can be extended by
constructing a non-associative algebraic structure which returns to
a given $n$-tuple a real value. For all $x\in \mathbb R^n$ and any
subset $I$ of $[n]$, let us  consider the map $\xi_I[x]:\Real
\longrightarrow \mathbb Z$ defined for all $\alpha \in \Real$ by
\begin{equation}\label{defxi}\xi_I[x](\alpha)= \Card
\{i\in I: x_i=\alpha\}- \Card \{i\in I: x_i=-\alpha\}.\end{equation}
This map measures the symmetry of the occurrences of a given value
$\alpha $ in the components of a vector $x$.

For all $x\in \mathbb R^n$, let $\mathcal J_I(x)$ be a subset of $I$
defined by
\begin{equation}\mathcal J_I(x)=\Big\{j\in I: \xi_I[x](x_j)\not=0\Big\}=I\backslash \big (\psi_I[x]^{-1}(0)\big),\end{equation}
where $\psi[x](j)= \xi_I[x](x_j)$ for all $j\in I$.
$\mathcal J_{I} (x)$ is called {\bf the residual  index set } of
$x$. It is obtained by dropping from $I$ all the $i$'s such that
$\Card \{j\in I: x_j=x_i\}= \Card \{j\in I: x_j=-x_i\}$.

For all positive natural numbers $n$ and
for any subset $I$ of $[n]$, let $\digamma_I: \Real^n
\longrightarrow \Real $ be the map defined for all $x\in \Real^n$ by

\begin{equation}\digamma_{ I}(x)=\left\{\begin{matrix}\max_{i\in \mathcal
J_I (x)}x_{i} &\text{if}\quad  \mathcal
J_I (x) \not=\emptyset    &   \text{and}\quad \xi_I[x]\big (\max_{i\in
\mathcal J_I(x)}|x_i|\big)>0 \\
\min_{i\in \mathcal J_I(x)}x_i &\hbox{ if }\quad \mathcal
J_I (x) \not=\emptyset  &\text{and}\quad \xi_I[x]\big(\max_{i\in \mathcal J_I(x)}|x_i|\big)<0 \\
 0 &\text{ if }\quad \mathcal
J_I (x)   =\emptyset, &.
\end{matrix}\right.\end{equation}
where $\xi_I[x]$ is the map defined in \eqref{defxi} and $\mathcal
J_I(x)$ is the residual index set of $x$. The operation that takes
an $n$-tuple $(x_1,....,x_n)$ of $\Real^n$ and returns a single real
element $\digamma_I(x_1,...,x_n)$ is called a $n$-ary extension of
the binary operation $\boxplus.$ For all natural numbers $n\geq 1$, if $I$ is
a nonempty subset of $[n]$, then, for all $n$-tuples $x=(x_1,...,x_n)$, one can define the operation:

\begin{equation}\label{defrecOp}\bigboxplus_{i\in
I}x_i= \lim_{p\longrightarrow
\infty}\stackrel{\varphi_p}{\sum_{i\in I}}x_i=\digamma_{I}(x).
\end{equation}
As an example, we have  
\begin{equation}
\digamma_{[12]} (8,5,-8,3,-5,2,1,-3,3,-5,5,1) =
\not 8 \boxplus \not 5 \boxplus \not\!\!\!{-8} \boxplus \not 3 \boxplus \not\!\!\!{-5} \boxplus 2 \boxplus 1 \boxplus \not\!\!\!{-3} \boxplus -3 \boxplus \not\!\!\!{-5} \boxplus \not 5 \boxplus 1 = -3.
\end{equation}

This operation encompasses as a special case the binary
operation defined in equation \eqref{base}. Clearly, for all
$(x_1,x_2)\in \Real^2$, $\bigboxplus_{i\in \{1,2\}}x_i=x_1\boxplus
x_2.$

There are some basic properties that can be inherited from the
above algebraic structure.  If all the elements of the family
$\{x_i\}_{i\in I}$ are mutually non symmetric, then:
$$\bigboxplus_{i\in
I}x_i=\arg\max_{|\cdot|,\lambda }\big\{|\lambda|: \lambda \in
\{x_i\}_{i\in I}\big\}.$$

The algebraic structure $(\Real,\boxplus,\cdot)$  can be extended
to $\Real^n$. Let us consider $m$ vectors $x^{(1)},...,x^{(m)}\in
\Real^n$, and define
\begin{align}
\bigboxplus_{j\in [m]}x^{(j)}&=\Big(\bigboxplus_{j\in [m]}
x_{1}^{(j)},...,\bigboxplus_{j\in [m]} x_{n}^{(j)}\Big).
\end{align}

Along this line an extended (relaxed for two points) definition of $\mathbb B$-convexity over $\Real^n$ was suggested in \cite{b15}. \\

$(iii)$ A subset of $\Real^n$ $C$
is {\bf idempotent symmetric convex} if for all $x,y\in C$ and all $t\in [0,1]$, $x\boxplus ty\in C$. \\

It has been equivalently established in  \cite{b15} that a set
$C\subset \Real^n$ is equivalently idempotent symmetric convex set if
for all $x,y\in C$, $Co^{\infty}(\{x,y\})\subset C$. In the remainder we will denote $Co^{\infty}(\{x,y\})=Co^\infty(x,y)$.  Moreover, for all $x,y\in \Real^n$,
$Co^\infty(x,y)$ is the smallest idempotent symmetric convex set that contains $x$ and $y$. In general, for any subset $ E $ of \( \mathbb{R}^n \), we denote by $ \mathrm{IS}[E] $ the smallest idempotent  symmetric convex set containing $ E $.

\subsection{Limit Systems and Determinant}\label{cramer}

\noindent In this section we present some results which were extensively studied in \cite{b20} from an algebraic viewpoint (see also \cite{b10} and \cite{butheged84} for a related concept of permanent). Let us denote $\M_{n}(\Real)$ the set of all the matrices
defined over the group $\Real^n$. Let us denote over $\M_{n}(\Real)$ the
map $\Phi_p:\M_{n}(\Real)\longrightarrow \M_{n}(\Real)$ defined
for all matrices $A= \left(a_{i,j}\right)_{\substack
{i=1...n\\j=1..n}} \in \M_{n}(\Real)$ by:
\begin{align}
\Phi_p(A)&=\left(\varphi_p(a_{i,j})\right)_{\substack
{i=1...n\\j=1...n}}= \left({a_{i,j}}^{2p+1}\right)_{\substack
{i=1...n\\j=1...n}}.\end{align} Its reciprocal is the map
$\Phi_p^{-1}:\M_{n}(\Real)\longrightarrow \M_{n}(\Real)$ defined
by:\begin{align}
\Phi_p^{-1}(A)&=\left(\varphi_p^{-1}(a_{i,j})\right)_{\substack
{i=1...n\\j=1...n}}=\left({a_{i,j}}^{\frac{1}{2p+1}}\right)_{\substack
{i=1...n\\j=1...n}}.
\end{align}
$\Phi_p$ is a natural extension of the map $\phi_p$ from $\Real^n$
to $\mathcal M_n(\Real)$. $\Phi_p(A)$ is the $2p+1$ Hadamard power
of matrix $A$. In addition, we introduce the matrix
product:

\begin{equation}
A \stackrel{p}{.}x=\sum_{j\in [n]}^{\varphi_p} x_{j}{.}a^{(j)},
\end{equation}
where $a^{(j)}$ is the $j$-th column of $A$. It is straightforward to
show that this formulation is equivalent to the following:

\begin{equation}
A \stackrel{p}{.}x=\phi_p^{-1}\big(\Phi_p(A).\phi_p(x)\big).
\end{equation}

If $f : \mathbb{R}^n \to \mathbb{R}^n$ is a linear endomorphism, then the map $f^{(p)} : \mathbb{R}^n \to \mathbb{R}^n$, defined for all $x \in \mathbb{R}^n$ by $f^{(p)}(x) = \phi_p^{-1} \circ f \circ   \phi_p(x) $,
is called the $\varphi_p$-linear transformation associated with $f$. By construction, this map is $\varphi_p$-linear. It is easy to see that the map $x\mapsto A\stackrel{p}{\cdot}x$ is
$\varphi_p$-linear. Conversely,  if $g$ is a $\varphi_p$-linear
map, then it can be represented by a matrix $A$ such that $g(x)=
A\stackrel{p}{\cdot}x$ for all $x\in \Real^n$.
Notice that the identity matrix $I$ is invariant with respect of
$\Phi_p$.

A $\varphi_p$-linear endomorphism is invertible if and only if its matrix representation
$\Phi_p(A)$ is invertible.  Therefore, it is useful to introduce a
suitable notion of determinant.  For all $n\times n$ matrices $A$,
let $|A|$ denotes its determinant.  The $\varphi_p$-determinant of a matrix $A\in \mathcal M_n(\Real)$ is defined by:
\begin{equation}|A|_p=\varphi_p^{-1}|\Phi_p(A)|.\end{equation}
Let $\mathfrak S_n$ denotes the set of all the permutations defined on $[n]$. Using the Leibnitz formula yields
\begin{equation}|A|_p=\Big(\sum_{\sigma\in \mathfrak S_n}\sgn(\sigma) \prod_{i\in [n]}a_{i,\sigma(i)}^{2p+1}\Big)^{\frac{1}{2p+1}}.\end{equation}
 This implies that if $f$ is a linear endomorphism having a matrix representation $A$, then
 $|f^{(p)}|_p=|A|_p$.  Moreover, we have:
\begin{equation}\lim_{p\longrightarrow \infty}|A|_p: =|A|_\infty=\bigboxplus_{\sigma\in \mathfrak S_n} \big(\sgn(\sigma)
 \prod_{i\in [n]}a_{i,\sigma(i)}\big).\end{equation}
 Let $f:\Real^n\longrightarrow \Real^n$ be a linear
endomorphism having a matrix representation $A$. For all $p\in
\mathbb N$, let $f^{(p)}$ be its $\varphi_p$-linear
transformation.   If $|A|_\infty\not=0$, then there is some
$p_0\in \mathbb N$ such that, for all $p\geq p_0$,  $f^{(p)}$ is
$\varphi_p$-invertible and for all $b\in \Real^n$, there
exists a solution $x^{(p)}$ to the system
$A\stackrel{p}{\cdot}x=b$  expressed as:
\begin{equation}\label{firstform}x^{(p)}_j=\frac{|A^{(j)}|_p}{|A|_p}=\frac{\left|\Phi_p(A^{(j)})\right|^{\frac{1}{2p+1}}}{\left |\Phi_p(A)\right |\frac{1}{2p+1}},\end{equation}
where $A^{(j)}$ is obtained from $A$ by dropping column $j$ and
replacing it with $b$.

Moreover, we have :
  \begin{equation} \lim_{p\longrightarrow \infty}x^{(p)}=x^\star\quad 
  \text{with}\quad  \text{for all}\quad  j\in[n] \quad x^\star_{j}=\frac{|A^{(j)}|_\infty}{|A|_\infty}.\end{equation}
  It follows that if for each $p$, $A\stackrel{p}{\cdot }x^{(p)}=b$, then:
  \begin{equation}\label{cramp}
 \frac{1}{ |A |_p}\stackrel{\varphi_p}{\sum_{j\in [n]}} \Big(\stackrel{\varphi_p}{\sum_{\sigma \in \mathfrak S_n}}\sgn(\sigma)\prod_{i\in [n]}a_{j,i,\sigma(i)}  \Big)a^{(j)}=b 
  \end{equation}
  where the $a_{j,i,k}$'s stand for the components of the matrix $A^{(j)}$. Taking the limit yields:
 \begin{equation}\label{cramplim}
  \bigboxplus_ {\substack{\sigma \in \mathfrak S_n\\ j\in [n]}}\lambda_{\sigma, j} a^{(j)}=b.
  \end{equation}
where $\lambda_{\sigma, j}= \sgn(\sigma)\frac{\prod_{i\in [n]}a_{j,i,\sigma(i)} }{ |A |_\infty}$ for each $j$. As in some properties of classical algebra, \( b \) can be expressed as a \(\mathbb{B}\)ox-plus combination of the \( a^{(j)} \)'s. However, contrary to the linear case, several terms (\( n! \)) involving each \( a^{(j)} \) are considered. These terms cannot be grouped because the operation \(\boxplus\) is not associative. This idea is systematically used in the remainder and also applies to the algebraic convex hull for two points further given in Proposition \ref{fond}, which was established in \cite{b15}.

Notice that since the integer numbers \( p \) are not fixed a priori, it is not possible to use the associativity property due to the lack of an addition table for each of them. Therefore, even to solve a \( \varphi_p \)-linear system, it is necessary to keep the history of all intermediate operations, and this is exactly the same thing for the cancellative sum $\boxplus$ considered to solve the limit system.

\section{Algebraic formulation}\label{fund}

In the following, we show that the limit polytope $Co^\infty(A)$ can be expressed algebraically using the $n$-ary operation $\bigboxplus$. The formula established below is closely related to the standard formula for polytopes, with the formal substitution $\sum \mapsto \bigboxplus$. However, it differs in an essential way: the nonnegative ``sum'' cannot be regrouped due to the non-associativity of the operation $\boxplus$. The key idea is   briefly summarized in the following. Suppose that  $A=\{x^{(1)},...,x^{(n+1)}\}$ is a finite subset of $\Real^n$.  Let $\Lambda=\{\lambda_{i,j}\}_{i,j\in [n+1]}$ with $\lambda_{i,j}=x_i^{(j)}$ for $i\in [n]$ and $j\in [n+1]$ and $\lambda_{n+1,j}=1$, for all $j$. Suppose moreover that there is some $p_0$ such that for all $p\geq p_0$, $x\in Co^p(A)$ and $|\Lambda|_p\not=0$ for some $p$. For each $p$, there is a solution to the system:
\begin{equation}
\Lambda\stackrel{p}{\cdot}t=\begin{pmatrix}x\\1\end{pmatrix}. 
\end{equation}
From the Cramer formula we obtain $t_j^{(p)}=\frac{|\Lambda^{(j)}|_p}{|\Lambda|_p}$ for all $j$, where $\Lambda^{(j)}$ is obtained from $\Lambda$ replacing the $j$-th column with $\begin{pmatrix}x\\1\end{pmatrix}$. It follows that from \eqref{cramp}:
\begin{equation}
x=\sum_{j\in [n+1]}t_i^{(p)}x^{(j)}=\frac{1}{|\Lambda|_p}\stackrel{\varphi_p}{\sum_{j\in [n+1]}}|\Lambda^{(j)}|_p x^{(j)}=\frac{1}{|\Lambda|_p}\stackrel{\varphi_p}{\sum_{j\in [n+1]}}\stackrel{\varphi_p}{\sum_{\sigma\in \mathfrak S_{n+1}}}\prod_{i \in [n+1]}\mathrm{sgn}(\sigma)\lambda_{i,\sigma(i)} ^{(j)}x^{(j)}. 
\end{equation}

Now, note that all the components $\mathrm{sgn}(\sigma)\lambda_{i,\sigma(i)} ^{(j)}$ are independent of $p$. Therefore, setting $\alpha_{\sigma,j}=\frac{1}{|\Lambda|_\infty}\prod_{i \in [n+1]}\mathrm{sgn}(\sigma)\lambda_{i,\sigma(i)} ^{(j)}$ for all $\sigma$, from \eqref{cramplim} we obtain taking the limit on the right side:
\begin{equation}
x=\sum_{j\in [n+1]}t_i^{(p)}x^{(j)}=\frac{1}{|\Lambda|_\infty}{\bigboxplus_{\substack{j\in [n+1]\\\sigma\in \mathfrak S_{n+1}}}}\prod_{i\in [n+1]}\mathrm{sgn}(\sigma)\lambda_{i,\sigma(i)} ^{(j)}x^{(j)}= {\bigboxplus_{\substack{j\in [n+1]\\\sigma\in \mathfrak S_{n+1}}}}\alpha_{\sigma,j}x^{(j)}. 
\end{equation}
In addition, note that since $t_j^{(p)}\geq 0$ for all $p$, we have $\bigboxplus_{\sigma\in \mathfrak S_{n+1}}\alpha_{\sigma,j}\geq 0$. Moreover, since  $\stackrel{\varphi_p}{\sum\limits_{j\in [n+1]}}t_j^{(p)}=1$, this implies that $\bigboxplus_{\substack{j\in [n+1]\\\sigma\in \mathfrak S_{n+1}}}\alpha_{\sigma,j} =\lim_{p\longrightarrow \infty}\stackrel{\varphi_p}{\sum\limits_{j\in [n+1]}}t_j^{(p)}=1.$
One can recognize the form of the standard convex hull, except that, for each $j$, the $\alpha_{\sigma,j}$ cannot be grouped due to the lack of associativity. However, in this particular situation, we have assumed that, for each $p$, $x \in Co^p(A)$; this is not true in general. Moreover, note that selecting, for each $p$, a point $x^{(p)}$ such that $x^{(p)} \to x$ is not sufficient, as there is no guarantee that the generalized sum will converge to the $\mathbb{B}$ox-Plus sum —its components are no longer independent of $p$.

However, this strategy can be useful for analyzing the limiting extreme points of the intersection between a $\mathbb{B}$-polytope and an $n$-dimensional orthant.
In the next section, we briefly recall some previously published results in the case of two points.

\subsection{Intermediate Points and Copositivity: the Case of two Points}

For simplicity, we shall denote $Co^{(p)}(x,y)$ instead of $Co^{(p)}(\{x,y\})$ in the following. The concept of intermediate points was considered in \cite{b15} to show that the sequence of sets $\{Co^{(p)}(x,y)\}_{p \in \mathbb{N}}$ has a Painlevé-Kuratowski limit, denoted by $Co^\infty(x,y)$, which can be expressed through an explicit algebraic formula in terms of the idempotent and non-associative algebraic structure defined in the paper. Two vectors $ u, v \in \mathbb{R}^n $ are said to be copositive if all their components have the same sign, i.e., $  u \odot v \geq 0 $, where $ \odot $ denotes the Hadamard product, defined as $ u \odot v = (u_1v_1, \ldots, u_nv_n) $.

To achieve this, a notion of intermediate points was introduced in \cite{b15}. Below, we briefly summarize the key results. To show that $Co^\infty(x,y)$ is the {Painlevé-Kuratowski} limit of the sequence $\{Co^{(p)}(x,y)\}_{p \in \mathbb{N}}$, it is useful to define the concept of $i$-intermediate points of order $p$, where $i\in [n]$ is such that $x_iy_i<0$.

  For all $(x,y)\in
\Real^n\times \Real^n$, let  $ \mathcal I(x,y)$ be the subset
defined by $ \mathcal I(x,y)=\{i\in [n]: x_iy_i<0\}$ and let
$n(x,y)$ be its cardinal.
For any natural number $p$, let us consider the map $\gamma^{(p)}:
\Real^n\times \Real^n\times [0,+\infty]\longrightarrow \Real$
defined by:
\begin{equation}\label{gammap}\gamma^{(p)}(x,y,t)=\Big(\frac{1}{1\stackrel{p}{+}t}\Big) x\stackrel{p}{+}
\Big(\frac{t}{1\stackrel{p}{+}t}\Big) y ,\quad \text{for all } t\geq
0
\end{equation} and by  $\gamma^{(p)}(x,y,+\infty)=y$.
For all $i\in \mathcal I(x,y)$, we say that
$\gamma^{(p)}(x,y,t_{i}^\star)$ is {\bf a $i$-intermediate point
of order $p$ }between $x$ and $y$ if
$\gamma_i^{(p)}(x,y,t_{i}^{\star})=0$.
For all $p\in \mathbb N$, there is a unique 
$i$-intermediate point of order $p$
\begin{equation}\label{gammap}\gamma^{(p)}(x,y,t_{i}^\star)=\Big(\frac{|y_i|}
{|x_i|\stackrel{p}{+}|y_i|}\Big)x\stackrel{p}{+}
\Big(\frac{|x_i|}{|x_i|\stackrel{p}{+}|y_{i}|}\Big)y,\end{equation} with
$t_i^\star=-\frac{x_i}{y_i}=|\frac{x_i}{y_i}|$.  For all $(x,y)\in
\Real^n\times \Real^n$, let us consider the map $\gamma:
\Real^n\times \Real^n\times [0,+\infty]\longrightarrow \Real^n$
defined by:
\begin{equation}\label{gamma}\gamma(x,y,t)=(\max\{1,t\})^{-1}\big(x
\boxplus {t}y\big),\quad \text{for all } t\geq 0
\end{equation} and by $\gamma(x,y,+\infty)=y$.

 Remark that $\gamma(x,y,0)=x$. For all $i\in \mathcal I(x,y)$ and all $t_{i}^\star \in \Real_{++}$,
a point $\gamma_{}\in \Real^n$ is called an {\bf $i$-intermediate
point} between $x$ and $y$ if there is some $t_{i}^\star\in
]0,+\infty[$ such that $\big(\gamma (x,y,t_{i}^\star)\big)_i=0.$
 For all $i\in \mathcal I(x,y)$, the map $t\mapsto
\gamma_i(x,y,t)$ has a unique  zero
$t_{i}^{\star}=-\frac{x_i}{y_i}=|\frac{x_i}{y_i}|>0$ and there is a
unique  $i$-intermediate point \begin{equation}\label{gammaty}\gamma
(x,y,t_{i}^\star)=\Big(\frac{|y_i|}{\max\{|x_i|,|y_i|\}}\,x\Big)
\boxplus \Big(\frac{|x_i|}{\max\{|x_i|,|y_i|\}}\,y\Big). \end{equation} Moreover, $\lim_{t\longrightarrow 0}\gamma^{}(x,y,t)=x$ and
$\lim_{t\longrightarrow +\infty}\gamma^{}(x,y,t)=y$. 

It was shown in \cite{b15} that for all $i\in \mathcal I(x,y)$
\begin{equation}\label{limgammap}
\lim_{p\longrightarrow \infty}\gamma^{(p)}(x,y,t_{i}^\star)=\gamma (x,y,t_{i}^\star):=\Big(\frac{|y_i|}{\max\{|x_i|,|y_i|\}}\,x\Big)
\boxplus \Big(\frac{|x_i|}{\max\{|x_i|,|y_i|\}}\,y\Big).
\end{equation}

\begin{center}
\begin{minipage}[h][5cm][t]{7cm}
\setlength{\unitlength}{0.99mm}
{\scriptsize \begin{picture}(90,20)(0,10)
\linethickness{0.127mm}
{\linethickness{0.15mm} \put(-30,0){\vector(1,0){60}}
\put(0,-20){\vector(0,1){45}}}
\put(20,20){\circle*{1}}
\put(20,22){$x^{(1)}$}
\put(20,19.8){\line(-1,-1){10}}
\put(10,10){\line(-1,0){20}}
\put(-10,-5){\circle*{1}}
\put(-16,-5){$y^{(1)}$}
\put(-10,10){\line(0,-1){15}}

\put(-10,20){\circle*{1}}
\put(-8,22){$x^{(2)}$}
\put(-20,10){\circle*{1}}
\put(-18,10){$y^{(2)}$}
\put(-10,20){\line(-1,0){10}}
\put(-20,20){\line(0,-1){10}}

\put(-4,4){\line(1,-1){10}}
\put(-4,4){\circle*{1}}
\put(6,-6){\circle*{1}}
\put(-7,6){$x^{(6)}$}
\put(8,-7){$y^{(6)}$}

\put(20,10){\circle*{1}}
\put(20,12){$x^{(5)}$}
\put(20,-5){\circle*{1}}
\put(20,-8){$y^{(5)}$}
\put(20,10){\line(0,-1){15}}

\put(20,-20){\circle*{1}}
\put(21,-23){$x^{(4)}$}
\put(-10,-10){\circle*{1}}
\put(-15,-10){$y^{(4)}$}
\put(20,-20){\line(-1,1){10}}
\put(10,-10){\line(-1,0){20}}

\put(-20,-20){\circle*{1}}
\put(-25,-20){$x^{(3)}$}
\put(-10,-15){\circle*{1}}
\put(-8,-16){$y^{(3)}$}
\put(-20,-20){\line(1,1){5}}
\put(-15,-15){\line(1,0){5}}

\put(0,0){\circle*{1}}
\put(1,1){$0 $}
\put(-0.5,27){$x_2$}
\put(32,-0.5){$x_1$}
\put(-30,-26){{\bf Figure \ref{fund}.1: Examples of
sets $Co^\infty(x,y)$.}}

\end{picture}}
{\scriptsize \begin{picture}(70,0)(-60,5.0)
\linethickness{0.127mm}
{\linethickness{0.15mm} \put(-20,0){\vector(1,0){50}}
\put(0,-21){\vector(0,1){45}}}
\put(30,30){\circle*{1}}
\put(30,32){$x $}
\put(30,30){\line(-1,-1){15}}
\put(15,15){\line(0,-1){35}}
\put(15,-20){\line(-1,0){30}}
\put(-15,-20){\circle*{1}}
\put(-16,-23){$y$}
\put(-15,-20){\circle*{1}}
\put(-16,-23){$y$}

\qbezier(30,30)(15,15)(14.1,0)
\qbezier(14.1,0)(15,-20)(-15,-20)
\put(14.0,0){\circle*{1}}
\put(15,0){\circle*{1}}
\put(10.1,2){$\gamma_{2}^{(p)}$}
\put(17,2){$\gamma_{2}$}
\put(0,-18.1){\circle*{1}}
\put(0,-20){\circle*{1}}
\put(-5,-17){$\gamma_{1}^{(p)}$}
\put(2,-22){$\gamma_{1}$}

\put(0,0){\circle*{1}}
\put(1,1){$0 $}
\put(-0.5,27){$x_2$}
\put(32,-0.5){$x_1$}
\put(-20,-26){{\bf
Figure \ref{fund}.2: Limit of $Co^p(x,y)$.}}

\end{picture}}
\end{minipage}
\end{center} 

\bigskip
\bigskip
Moreover, for all $x,y\in \mathbb R^n$ such that $\mathcal
I(x,y)\not=\emptyset$, there is a sequence of indexes  $\{i_m\}_{m\in [n(x,y)]}\subset \mathcal
I(x,y)$   such that setting  $t_{i_0}^\star=0$, $t_{i_{n(x,y)+1}}^\star=+\infty$ and
$t_{i_m}^\star=-\frac{x_{i_{m}}}{y_{i_{m}}}$ for all $m\in
[n(x,y)]$:
\begin{equation}\label{conseq} \gamma^{(p)}
\big(x,y,t^\star_{i_{m}}\big)\odot \gamma^{(p)}
\big(x,y,t^\star_{i_{m+1}}\big)\geq 0\quad \text{ and }\quad \gamma
\big(x,y,t^\star_{i_{m}}\big)\odot \gamma
\big(x,y,t^\star_{i_{m+1}}\big)\geq 0.
\end{equation} The sequences of intermediate points defined in equation \eqref{conseq} (whether at order 
$p$ or at the limit) are copositive, meaning that the consecutive terms of these sequences have the same sign. This property is important because it allows for determining the algebraic form of the limit of the convex envelope for two points. Using the fact that:
\begin{equation}\label{gammapdecomp}
Co^p(x,y)=\bigcup_{m=0}^{n(x,y)}Co^p\Big(\gamma^{(p)}
\big(x,y,t_{i_m}^\star\big),\gamma^{(p)}
\big(x,y,t_{i_{m+1}}^\star\big)\Big),
\end{equation}
the following result was established in \cite{b15}:
\begin{equation}
 Co^{\infty}(x,y)=\Lim_{p\longrightarrow +\infty}
Co^p(x,y)=\bigcup_{m=0}^{n(x,y)}
Co^\infty\Big(\gamma_{}(x,y,t_{i_m}^\star),\gamma_{}(x,y,t_{i_{m+1}}^\star)\Big). \label{fond}
\end{equation}

In general given $a_1, ..., a_m\in \Real$, if $(I_1,I_2)$ is a partition of $[m] $ there is no guarantee that 
$\Big(\bigboxplus\limits_{i\in  I_1 }a_i\Big)\boxplus \Big(\bigboxplus\limits_{i\in I_2}a_i\Big)$ is identical to $\bigboxplus\limits_{i\in  I_1\cup I_2 }a_i$. However, this is true if $\bigboxplus\limits_{i\in I_1}a_i$ and $\bigboxplus\limits_{i\in I_2}a_i$ are not symmetric or are equal to $0$. Since for all $m$, $\gamma_{}(x,y,t_{i_m}^\star)$ and $\gamma_{}(x,y,t_{i_{m+1}}^\star)$ are copositive, we have for all $r,s\in [0,1]$  with $\max\{r,s\}=1$ the relation:
{ \begin{align}
r \gamma_{}(x,y,t_{i_m}^\star)\boxplus s \gamma_{}(x,y,t_{i_{m+1}}^\star)&=\big[(rw_m)\,x
\boxplus (r w_m)\,y\big]\boxplus \big[(sw_{m+1})\,x
\boxplus (sw_{m+1})\,y\big]\notag \\&=  (rw_m)\,x
\boxplus (r w_m)\,y\boxplus (sw_{m+1})\,x
\boxplus (sw_{m+1})\,y,
\end{align}}
where  $w_k=\frac{|y_{i_k}|}{\max\{|x_{i_k}|,|y_{i_k}|\}}$ for $k\in \{m,m+1\}$. If $\max\{t,r,s,w\}=1$ for the real numbers $t,r,s,w\geq 0$, then:
\begin{equation}
\lim_{p\longrightarrow \infty} \frac{1}{t\stackrel{p}{+}r\stackrel{p}{+}s\stackrel{p}{+}w}\Big( {tx\stackrel{p}{+}r x\stackrel{p}{+}sy\stackrel{p}{+}wy}\Big)=t x\boxplus r x\boxplus s y\boxplus w y.
\end{equation}
 As an immediate consequence of Proposition \ref{fond}, one can then give the following
formulation of a $\mathbb B$-polytope in the case of two points\footnote{We use the convention that for all $(\alpha,\beta,\gamma,\delta)\in \mathbb R^4$, $ \alpha\boxplus \beta\boxplus \gamma\boxplus \delta=\digamma_{[4]}(\alpha,\beta,\gamma,\delta)$.}:
\begin{prop}\cite{b15}\label{fond}
For all $x,y\in \Real^n$,
$$Co^{\infty}(x,y)=\Big \{t x\boxplus r x\boxplus s y\boxplus w y:
\max\{t,r,s,w\}=1,t,r,s,w\geq 0\Big\}.$$

\end{prop}
 
The copositivity of the sequence of intermediate points plays an important role. Let $n$ be a positive natural number and let $I$ be a nonempty subset of $[n]$. Let $\mathfrak{P}(I) = \{I_j : j \in [m]\}$ be a partition of $I$ into $m$ nonempty subsets $I_j$.
 In \cite{b15} it was shown that if
  for all $(j,k)\in [m]\times
[m]$
\begin{equation}\label{recompos}\Big(\bigboxplus_{i\in I_j} x_{i}\Big)+\Big(\bigboxplus_{i\in I_k} x_{i}\Big)\not=0\;\text{ then }
\;\bigboxplus_{j\in [m]}\Big(\bigboxplus_{i\in I_j}
x_{i}\Big)=\bigboxplus_{i\in I}x_{i}.\end{equation}

\begin{lem} Suppose that there exists $m$ subsets of $\mathbb N$, $I_1, I_2...I_m$, with $x^{(i)}\in\Real^n$ for all $i\in I_j$. Suppose moreover that for all $j,k$, $\bigboxplus_{i\in I_j}x^{(i)}$ and $\bigboxplus_{i\in I_k}x^{(i)}$ are copositive i.e. have   the same sign which equivalently means  $\Big(\bigboxplus_{i\in I_j}x^{(i)}\Big)\odot \Big(\bigboxplus_{i\in I_k}x^{(i)}\Big)\geq 0 $. Then:
$$\bigboxplus_{j\in [m]}\Big(\bigboxplus_{i\in I_j}
x^{(i)}\Big)=\bigboxplus_{\substack{  i\in I_j\\j\in [m]}}x^{(i)}. $$

\end{lem}
 
\noindent{\bf Proof:}
Fix some $h\in [n]$. If for all $ j,k$ $\Big(\bigboxplus_{i\in I_j}x_h^{(i)}\Big)\odot \Big(\bigboxplus_{i\in I_k}x_h^{(i)}\Big)>0$ then $\Big(\bigboxplus_{i\in I_j}x_h^{(i)}\Big)+ \Big(\bigboxplus_{i\in I_k}x_h^{(i)}\Big)\not=0$
and from \cite{b15} and equation \eqref{recompos}
$$\bigboxplus_{j\in [m]}\Big(\bigboxplus_{i\in I_j}
x_h^{(i)}\Big)=\bigboxplus_{\substack{  i\in I_j\\j\in [m]}}x_h^{(i)}. $$ 
Suppose that $\Big(\bigboxplus_{i\in I_{j_0}}x_h^{(i)}\Big)\odot \Big(\bigboxplus_{i\in I_{k_0}}x_h^{(i)}\Big)=0$ for some ${j_0}$ and ${k_0}$, then this implies that either $\bigboxplus_{i\in I_{j_0}}x_h^{(i)}=0$ or $\bigboxplus_{i\in I_{k_0}}x_h^{(i)}=0$. Let $J_0=\{j\in [m]:\bigboxplus_{i\in I_{j}}x_h^{(i)}=0\}$. We have: 
 
$$\bigboxplus_{j\in [m]  }\Big(\bigboxplus_{i\in I_j}
x_h^{(i)}\Big)=\bigboxplus_{j\in [m]\backslash  J_0  }\Big(\bigboxplus_{i\in I_j}
x_h^{(i)}\Big)=\bigboxplus_{\substack{  i\in I_j\\j\in [m]\backslash  J_0 }}x_h^{(i)} =\bigboxplus_{\substack{  i\in I_j\\j\in [m]}}x_h^{(i)}.  $$

Since the result can be applied to each $h$, we obtain the result. $\Box$\\

\subsection{ Intersection with $n$-dimensional Orthant, Extreme  Points in Limits and Inclusion Property }   \label{Polalg}
 
In this subsection, we show that any $\mathbb{B}$-polytope in $\mathbb{R}^n$ is a subset of another set that can be described using the idempotent and non-associative algebraic structure considered in the paper. To do this, we provide, for all natural numbers $p$, an explicit algebraic form for the extreme points of the $\varphi_p$-polytope $Co^p(A) \cap K$, where $K$ is an $n$-dimensional orthant of $\mathbb{R}^n$. In this context, it is shown that these extreme points have an explicit limit, which allows the inclusion of $Co^\infty(A) \cap K$ in this set. To achieve this, the notion of limit determinant defined in \ref{cramer} and in \cite{b20} is very useful.

For all natural numbers $m>1$, let  $\{e_1,...,e_m\}$ denote the canonical basis of $\mathbb{R}^m$. For all nonempty subsets $J$ of $[m]$, let us denote $\Delta_J^{(p)}$ the $\varphi_p$-simplex defined by:
\begin{equation}
\Delta_J^{(p)}=\Big\{\sum_{j\in J}t_je_j:\stackrel{\varphi_p}{\sum\limits_{j\in J}}t_j=1, t_j\geq 0, j\in J \Big\}.
\end{equation}

Let $A=\{x^{(1)},...,x^{(m)}\}$ be a finite subset of $\mathbb{R}^n$. For all subsets $J$ of ${[m]}$, let $\Gamma_{J}: \mathbb{R}^{n\times m}\times \mathbb{R}_+^m\longrightarrow \mathbb{R}^n$ be the map defined by:
\begin{equation}
\Gamma_J^{(p)}\big({x^{(1)},...,x^{(m)},t}\big)=\stackrel{\varphi_p}{\sum_{j\in J}}t_jx^{(j)}.
\end{equation}

Let $I$ be a subset of $[n]$ such that $\Card\, I  = \Card\, J  - 1$.  Let \( t_{_I,_J}^{(p)} \in \Delta_J^{(p)} \). We define the point
\begin{equation}
\zeta_{_I,_J}^{(p)} := \Gamma_J^{(p)}\big(x^{(1)}, \dots, x^{(m)}, t_{_I,_J}^{(p)} \big)
\end{equation}
as an \textbf{intermediate point of order \( p \) for \( A \), indexed by the pair \( (I, J) \)}, if \( t_{_I,_J}^{(p)} \) is the unique element of \( \Delta_J^{(p)} \) satisfying
\begin{equation}
\left( \Gamma_J^{(p)}\big(x^{(1)}, \dots, x^{(m)}, t_{_I,_J}^{(p)}\big) \right)_i = 0 \quad \text{for all } i \in I.
\end{equation}
In this case, \( t_{_I,_J}^{(p)} \) is referred to as an \textbf{intermediate value of order \( p \) for \( A \), indexed by \( (I, J) \)}. By convention, each initial point \( x^{(j)} \) is considered to be a \( (\{j\}, \emptyset) \)-intermediate point.

\begin{center}
{\scriptsize

\unitlength 1mm % = 2.845pt
\linethickness{0.4pt}
\ifx\plotpoint\undefined\newsavebox{\plotpoint}\fi % GNUPLOT compatibility
\begin{picture}(127.866,56.149)(0,0)
\put(95.159,5.076){\vector(0,1){45.289}}
%\vector(81.8,16.42)(118.39,35.269)
\put(118.39,35.269){\vector(2,1){.07}}\multiput(81.8,16.42)(.0654561717,.0337191413){559}{\line(1,0){.0654561717}}
%\end
\put(104.688,52.583){\circle*{.771}}
\put(84.334,36.378){\circle*{.771}}
\put(85.126,49.939){\circle*{.771}}
\put(95.089,28.067){\circle*{.771}}
\put(111.896,50.28){\circle*{.771}}
\put(101.796,16.203){\circle*{.771}}
\put(83.621,17.102){\circle*{.771}}
\put(92.888,7.464){\circle*{.771}}
%\emline(84.968,50.195)(95.897,55.653)
\multiput(84.968,50.195)(.067462963,.033691358){162}{\line(1,0){.067462963}}
%\end
%\dashline{1}(115.064,19.917)(113.955,19.064)
\multiput(114.994,19.847)(-.0426538,-.0328077){13}{\line(-1,0){.0426538}}
%\end
%\dashline{1}(81.72,30.493)(92.808,20.087)
\multiput(81.65,30.423)(.0343282,-.0322167){19}{\line(1,0){.0343282}}
\multiput(82.954,29.198)(.0343282,-.0322167){19}{\line(1,0){.0343282}}
\multiput(84.259,27.974)(.0343282,-.0322167){19}{\line(1,0){.0343282}}
\multiput(85.563,26.75)(.0343282,-.0322167){19}{\line(1,0){.0343282}}
\multiput(86.868,25.526)(.0343282,-.0322167){19}{\line(1,0){.0343282}}
\multiput(88.172,24.302)(.0343282,-.0322167){19}{\line(1,0){.0343282}}
\multiput(89.477,23.077)(.0343282,-.0322167){19}{\line(1,0){.0343282}}
\multiput(90.781,21.853)(.0343282,-.0322167){19}{\line(1,0){.0343282}}
\multiput(92.085,20.629)(.0343282,-.0322167){19}{\line(1,0){.0343282}}
%\end
%\dashline{1}(81.641,30.322)(94.947,36.719)
\multiput(81.571,30.252)(.069302,.033318){12}{\line(1,0){.069302}}
\multiput(83.234,31.051)(.069302,.033318){12}{\line(1,0){.069302}}
\multiput(84.897,31.851)(.069302,.033318){12}{\line(1,0){.069302}}
\multiput(86.56,32.651)(.069302,.033318){12}{\line(1,0){.069302}}
\multiput(88.224,33.45)(.069302,.033318){12}{\line(1,0){.069302}}
\multiput(89.887,34.25)(.069302,.033318){12}{\line(1,0){.069302}}
\multiput(91.55,35.049)(.069302,.033318){12}{\line(1,0){.069302}}
\multiput(93.213,35.849)(.069302,.033318){12}{\line(1,0){.069302}}
%\end
%\dashline{1}(92.65,20.002)(92.808,7.806)
\put(92.58,19.932){\line(0,-1){.9382}}
\put(92.604,18.055){\line(0,-1){.9382}}
\put(92.628,16.179){\line(0,-1){.9382}}
\put(92.653,14.303){\line(0,-1){.9382}}
\put(92.677,12.426){\line(0,-1){.9382}}
\put(92.701,10.55){\line(0,-1){.9382}}
\put(92.726,8.674){\line(0,-1){.9382}}
%\end
%\emline(85.047,50.11)(89.64,39.448)
\multiput(85.047,50.11)(.033525547,-.077824818){137}{\line(0,-1){.077824818}}
%\end
%\emline(89.64,39.448)(79.265,33.563)
\multiput(89.64,39.448)(-.059285714,-.033628571){175}{\line(-1,0){.059285714}}
%\end
%\emline(79.503,33.734)(81.404,30.664)
\multiput(79.503,33.734)(.03335088,-.05385965){57}{\line(0,-1){.05385965}}
%\end
\put(81.784,30.332){\line(0,-1){11.381}}
%\emline(92.71,7.44)(81.784,19.017)
\multiput(92.71,7.44)(-.0337222222,.0357314815){324}{\line(0,1){.0357314815}}
%\end
%\emline(96.129,55.579)(127.866,44.591)
\multiput(96.129,55.579)(.0973527607,-.0337055215){326}{\line(1,0){.0973527607}}
%\end
%\emline(127.866,44.591)(106.758,26.866)
\multiput(127.866,44.591)(-.0401292776,-.0336977186){526}{\line(-1,0){.0401292776}}
%\end
\put(106.758,26.866){\line(0,-1){11.446}}
%\emline(106.758,15.42)(92.785,7.44)
\multiput(106.758,15.42)(-.058957806,-.033670886){237}{\line(-1,0){.058957806}}
%\end
%\dashline{1}(94.42,36.481)(106.609,26.931)
\multiput(94.35,36.411)(.0421765,-.033045){17}{\line(1,0){.0421765}}
\multiput(95.784,35.287)(.0421765,-.033045){17}{\line(1,0){.0421765}}
\multiput(97.218,34.164)(.0421765,-.033045){17}{\line(1,0){.0421765}}
\multiput(98.652,33.04)(.0421765,-.033045){17}{\line(1,0){.0421765}}
\multiput(100.086,31.917)(.0421765,-.033045){17}{\line(1,0){.0421765}}
\multiput(101.52,30.793)(.0421765,-.033045){17}{\line(1,0){.0421765}}
\multiput(102.954,29.67)(.0421765,-.033045){17}{\line(1,0){.0421765}}
\multiput(104.388,28.546)(.0421765,-.033045){17}{\line(1,0){.0421765}}
\multiput(105.822,27.422)(.0421765,-.033045){17}{\line(1,0){.0421765}}
%\end
%\dashline{1}(106.609,26.931)(92.71,20.195)
\multiput(106.539,26.861)(-.0668221,-.0323846){13}{\line(-1,0){.0668221}}
\multiput(104.801,26.019)(-.0668221,-.0323846){13}{\line(-1,0){.0668221}}
\multiput(103.064,25.177)(-.0668221,-.0323846){13}{\line(-1,0){.0668221}}
\multiput(101.327,24.335)(-.0668221,-.0323846){13}{\line(-1,0){.0668221}}
\multiput(99.589,23.493)(-.0668221,-.0323846){13}{\line(-1,0){.0668221}}
\multiput(97.852,22.651)(-.0668221,-.0323846){13}{\line(-1,0){.0668221}}
\multiput(96.114,21.809)(-.0668221,-.0323846){13}{\line(-1,0){.0668221}}
\multiput(94.377,20.967)(-.0668221,-.0323846){13}{\line(-1,0){.0668221}}
%\end
%\vector(108.319,9.533)(79.703,39.489)
\put(79.703,39.489){\vector(-1,1){.07}}\multiput(108.319,9.533)(-.03370553592,.03528386337){849}{\line(0,1){.03528386337}}
%\end
{\linethickness{0.15mm}  
\put(36,9.527){\vector(0,1){40}}
}
\put(66,49.527){\line(-1,-1){15}}
\put(56,39.827){\line(-1,0){40}}
\put(16,39.827){\line(0,-1){5}}
\put(16,34.827){\line(4,-3){10}}
\put(26,27.527){\line(0,-1){15.5}}
\put(51,34.827){\line(0,-1){22.5}}
\put(51,12.227){\line(-1,0){25.2}}
\put(36,19.527){\line(0,-1){10}}
\put(33.5,16.527){$0$}
\put(36,19.527){\circle*{1}}
\put(36,19.527){\circle*{1}}
\put(51,24.527){\circle*{1}}
\put(54,24.527){$x^{(1)}$}
\put(16,34.827){\circle*{1}}
\put(10,33.827){$x^{(2)}$}
\put(26,12.227){\circle*{1}}
\put(20,11.527){$x^{(3)}$}
\put(66,49.527){\circle*{1}}
\put(69,49.527){$x^{(4)}$}
\put(35.5,51.527){$x_2$}
\put(68,19.027){$x_1$}
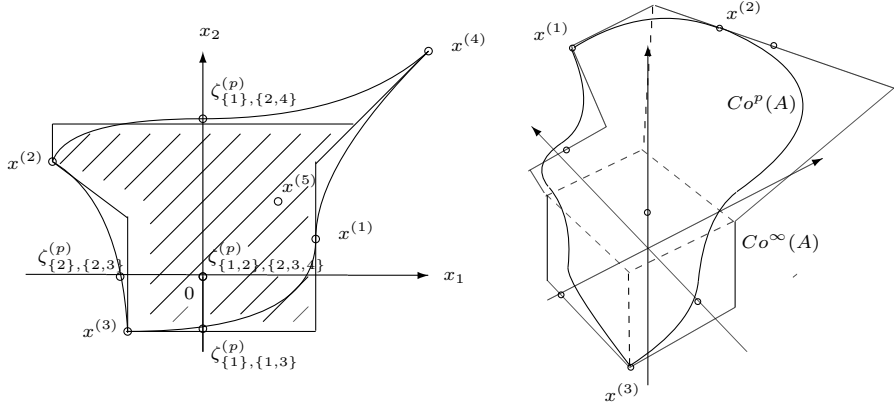
\put(0,0){{\bf Figure \ref{Polalg}.1:} 2-dimensional $\varphi_p$-polytope and intermediate points.}
\put(88.501,0){{\bf Figure \ref{Polalg}.2:} 3-dimensional simplex in limit.}
\put(21,38.527){\line(-1,-1){4}}
\put(25,38.527){\line(-1,-1){5}}
\put(29,38.527){\line(-1,-1){6.2}}
\put(33,38.527){\line(-1,-1){7.4}}
\put(37,38.527){\line(-1,-1){8.9}}
\put(41,38.527){\line(-1,-1){12}}
\put(45,38.527){\line(-1,-1){16}}
\put(49,38.527){\line(-1,-1){20}}
\put(53,38.527){\line(-1,-1){24}}
\put(50,31.527){\line(-1,-1){8}}
\put(50,27.527){\line(-1,-1){13}}
\put(46,19.527){\line(-1,-1){5}}
\put(50,19.527){\line(-1,-1){6}}
%\emline(50,15.527)(48,13.527)
\multiput(50,15.527)(-.03333333,-.03333333){60}{\line(0,-1){.03333333}}
%\end
%\emline(34,15.527)(32,13.527)
\multiput(34,15.527)(-.03333333,-.03333333){60}{\line(0,-1){.03333333}}
%\end
{
\qbezier(66,49.527)(56,40.527)(36,40.527)
\qbezier(36,40.527)(18,40.527)(16,34.827)
\qbezier(16,34.827)(25,29.527)(26,12.227)
\qbezier(26,12.227)(51,12.227)(51,24.527)
\qbezier(51,24.527)(51,34.527)(66,49.527)
}
\put(25,19.527){\circle*{1}}
\put(14,21.527){$\zeta^{(p)}_{\{2\},\{2,3\}}$}
\put(36,12.587){\circle*{1}}
\put(37,8.527){$\zeta^{(p)}_{\{1\},\{1,3\}}$}
\put(36,40.527){\circle*{1}}
\put(37,43.527){$\zeta^{(p)}_{ \{1\},\{2,4\}}$}
\put(36.6,21.527){$\zeta^{(p)}_{\{1,2\},\{2,3,4\}}$}
\put(46,29.527){\circle*{1}}
\put(46.5,30.727){$x^{(5)}$}
%\vector(12.439,19.824)(65.953,19.676)
\put(65.953,19.676){\vector(1,0){.07}}\multiput(12.439,19.824)(10.7028,-.0296){5}{\line(1,0){10.7028}}
%\end
\put(82.379,52.156){\makebox(0,0)[cc]{$x^{(1)}$}}
\put(91.818,4.142){\makebox(0,0)[cc]{$x^{(3)}$}}
\put(108.021,54.98){\makebox(0,0)[cc]{$x^{(2)}$}}
%\dashline{1}(95.853,55.507)(94.592,36.482)
\put(95.783,55.437){\line(0,-1){.9513}}
\put(95.657,53.534){\line(0,-1){.9512}}
\put(95.531,51.632){\line(0,-1){.9513}}
\put(95.404,49.729){\line(0,-1){.9513}}
\put(95.278,47.827){\line(0,-1){.9513}}
\put(95.152,45.924){\line(0,-1){.9513}}
\put(95.026,44.022){\line(0,-1){.9513}}
\put(94.9,42.119){\line(0,-1){.9512}}
\put(94.774,40.217){\line(0,-1){.9513}}
\put(94.648,38.314){\line(0,-1){.9513}}
%\end
\qbezier(104.604,52.569)(95.147,56.149)(85.336,50.006)
\qbezier(104.781,52.658)(125.597,43.951)(106.991,30.826)
\qbezier(106.902,30.737)(102.704,27.467)(101.864,18.893)
\qbezier(101.864,18.893)(101.29,13.192)(92.76,7.668)
\qbezier(84.982,49.829)(88.827,40.725)(83.656,37.278)
\qbezier(83.656,37.278)(79.59,34.494)(81.888,31.533)
\qbezier(82.065,31.356)(85.247,27.467)(84.54,20.75)
\qbezier(84.54,20.75)(84.452,18.496)(92.849,7.58)
\put(110.322,42.437){\makebox(0,0)[cc]{$Co^p(A)$}}
\put(112.839,23.937){\makebox(0,0)[cc]{$Co^\infty(A)$}}
\end{picture}

}
\end{center}

It is immediate to see that this notion generalizes the concept of $i$-intermediate points defined in the earlier subsection and in \cite{b15}. For all subsets $H \subseteq \mathbb{N} \backslash \{0\}$, let $\nu_H$ denote the unique increasing bijection defined from $H$ to $[\Card \, H]$, i.e., for all $l,m \in H$, if $m > l$, then $\nu_H(m) > \nu_H(l)$. For simplicity, let us denote $n_J = \Card \, J$. By construction, $\nu_J(J)=[n_J]$ and $\nu_J^{-1}([n_J])=J$.

Suppose that $I\subset [n]$ and $J\subset [m]$ such that $\Card \, I=\Card\, J\,-1$. This implies that $n_J=\Card \, I\,+1$. For any finite subset $A=\{x^{(1)},...,x^{(m)}\}$ of $\mathbb{R}^n$, let:
\[
\Lambda_{_I,_J} =\big(\lambda_{i,j}\big)_{\substack{i,j\in [n_J]}}
\]
be the square matrix of $\mathcal{M}_{ n_J}(\mathbb{R})$ defined by:

\begin{equation}
\lambda_{i,j} =\begin{cases}x_{\nu_I^{-1}(i)}^{(\nu_J^{-1}(j))} & \text{if } (i,j)\in [n_J-1]\times [n_J], \\
1 & \text{if } i=n_J, j \in [n_J]. \end{cases}
\end{equation}

In the remainder, we say that \textbf{$\Lambda_{_I,_J}$ is an intermediate matrix for $A$ indexed on $(I,J)$}. By construction, it follows that $t_{_I,_J}^{(p)}\in \Delta_J^{(p)}$ is an intermediate value of order $p$ for $A$ if and only if there exists some $s^\star\in \Delta _{[n_J]}^{(p)}$ satisfying the $\varphi_p$-matrix equation:
\begin{equation}
\Lambda_{_I,_J}\stackrel{p}{\cdot}s=e_{ n_J},
\end{equation}
such that $s^\star=\sum\limits_{j\in J}t_{\nu_J(j)}e_{\nu_J(j)}$.

In the following, for all $k\in [n_J]$, let $\Lambda_{_I,_J}^{(k)}=\Big(\lambda^{(k)}_{i,j}\Big)_{\substack{i,j\in [n_J]}}$ denote the matrix obtained by replacing the $k$-th column of $\Lambda_{_I,_J}$ with $e_{ n_J }$. Along this line, one can retrieve as a special case the intermediate points in the case where $m=2$.

\begin{rem} \textnormal{Let $x^{(1)},x^{(2)}\in \Real^n$. Suppose that there is some $i\in [n]$ such that $x_i^{(1)}x^{(2)}_i<0$.
In such a case, we have
$$\Lambda_{\{i\},\{1,2\}}=\begin{pmatrix}x_{i}^{(1)}&x_{i}^{(2)}\\1&1\end{pmatrix}.$$
Since $\left
|\begin{matrix}x_{i}^{(1)}&x_{i}^{(2)}\\1&1\end{matrix}\right|_p
=\Big(\big({x_{i}^{(1)}}\big)^{2p+1}-\big({x_{i}^{(2)}}\big)^{2p+1}\Big)^{\frac{1}{2p+1}}\not=0$,
the Cramer solutions of the system
$$\begin{pmatrix}x_{i}^{(1)}&x_{i}^{(2)}\\1&1\end{pmatrix}
\stackrel{p}{\cdot}\begin{pmatrix}t_1\\t_2\end{pmatrix}=\begin{pmatrix}0\\1\end{pmatrix}$$
are
 $$t_{\{i\},\{1,2\},1}^{(p)}=\frac{\left| \begin{matrix}0&x_{i}^{(2)}\\1&1\end{matrix}\right|_p}
 {\left|\begin{matrix}x_{i}^{(1)}&x_{i}^{(2)}\\1&1\end{matrix}\right|_p}
=\frac{|x_{i}^{(2)}|}{(|x_{i}^{(1)}|^{2p+1}+|x_{i}^{(2)}|^{2p+1})^\frac{1}{2p+1}}$$
and $$ t_{\{i\},\{1,2\},2}^{(p)}=\frac{\left|
\begin{matrix}x_{i}^{(1)}&0\\1&1\end{matrix}\right|_p}{\left|
\begin{matrix}x_{i}^{(1)}&x_{i}^{(2)}\\1&1\end{matrix}\right|_p}=
\frac{|x_{i}^{(1)}|}{(|x_{i}^{(1)}|^{2p+1}+|x_{i}^{(2)}|^{2p+1})^\frac{1}{2p+1}}.$$
We retrieve the solutions established in \cite{b15}. In particular
$$\zeta_{\{i\},\{1,2\}}^{(p)}=\Gamma_{ \{1,2\}}\Bigg(x^{(1)}, x^{(2)},\Big(\frac{|x_{i}^{(2)}|}{(|x_{i}^{(1)}|^{2p+1}+|x_{i}^{(2)}|^{2p+1})^\frac{1}{2p+1}}, \frac{|x_{i}^{(1)}|}{(|x_{i}^{(1)}|^{2p+1}+|x_{i}^{(2)}|^{2p+1})^\frac{1}{2p+1}}\Big)\Bigg).$$}
\end{rem}

In the following, given a $\varphi_p$-convex subset $C$ of
$\Real^n$, we say that $z$ is a $\varphi_p$-extreme point if there
do not exist $u,v\in C$ and $(s,t)\in \Delta_2^{(p)}$, such that
$z=su\stackrel{p}{+}tv$. Given  {an} $n$-dimensional orthant $K$, the
subset $Co^p(A)\cap K$ has a $\varphi_p$-internal representation.
This means that there exists a finite subset $\mathcal E_K^{(p)}(A)$ of
$\varphi_p$-extreme points such that \begin{equation}Co^p(A)\cap
K=Co^p\Big(\mathcal E_K^{(p)}(A)\Big ).\end{equation} In the following, we say that {\bf  an  intermediate point is properly indexed on   $(I,J)$   } if  $\Lambda_{_I,_J}$ is $\varphi_p$-invertible. In the next lemma, we show that any $\varphi_p$-extreme point is a an intermediate point of order $p$ properly indexed on some pair $(I,J)$. In the following, $\mathcal Z^{(p)}(A)$ denotes the set of all the intermediate points of order $p$ for $A$ that are properly indexed.

\begin{lem}Let \( p \) be a natural number and suppose that \( Co^p(A) \cap K \) is nonempty. If \( z \) is a \( \varphi_p \)-extreme point of \( Co^p(A) \cap K \), then there exists a pair \( (I, J) \) such that \( z \) is an   intermediate point of order \( p \) properly indexed on the pair \( (I, J) \). 

\end{lem}
{\bf Proof:} Suppose that there exists some $j$ such that $z = x^{(j)}$. Then, by definition, $z$ is an  intermediate point of order $p$ properly indexed on the pair $({\emptyset, \{j\}})$, and the result follows immediately.

Suppose now that this is not the case. This implies that $z$ is not a \( \varphi_p \)-extreme point of $Co^p(A)$. Since $z$ is a \( \varphi_p \)-extreme point of $Co^p(A) \cap K$,  it belongs to $Co^p(A)$, and there exists a subset $J \subset [m]$ and a set of points $\{x^{(j)}\}_{j \in J}$ having full rank $n_J$, along with coefficients $(t_j)_{j \in J}$ such that $t_j > 0$ for all $j \in J$, satisfying
$$
\stackrel{\varphi_p}{\sum_{j \in J}} t_j = 1, \quad \text{and} \quad z = \stackrel{\varphi_p}{\sum_{j \in J}} t_j x^{(j)}.
$$
If \( z \) belongs to the interior of \( K \), then, since it is not a \( \varphi_p \)-extreme point of \( Co^p(A) \), there exist points \( z', z'' \in \mathrm{int} \, K \cap Co^p(A) \) such that \( z \in Co^p(z',z'') \). This contradicts the assumption that \( z \) is a \( \varphi_p \)-extreme point of \( Co^p(A) \cap K \). Therefore, \( z \in \partial K \), and there exists a subset \( I \) of \( [n] \) such that \( z_i = 0 \) for all \( i \in I \).
  Suppose now that $\mathrm{Card}\,I \geq n_J - 1$. Since the rank of $\{x^{(j)}\}_{j \in J}$ is $n_J$, the associated system $
\Lambda_{I, J} \stackrel{p}{\cdot} s = e_{n_J}$ has a unique solution. The problem is then resolved by selecting $I' \subset I$ with $\mathrm{Card}\,I' = n_J - 1$.

Now, suppose that $\mathrm{Card}\, I < n_J - 1$, with $z_i \neq 0$ for $i \notin I$. We will derive a contradiction. Let $H$ be a $\varphi_p$-hyperplane of dimension $n_J - 1$ passing through $x^{(j)}$ for all $j \in J$. Since $t_j > 0$, there exists a ball $B_\infty(z, r]$ for some $r > 0$ such that
\[
B(z, r] \cap H \subset \mathrm{rint} \Big(Co^p(A) \cap H \Big),
\]
where $\mathrm{rint}$ denotes the relative interior of $H$.

Define $F = \{ u \in \mathbb{R}^n : u_i = 0, \forall i \in I \}$. The dimension of $F$ is at least $1$. Since $z \in H \cap F$, there exists a $\varphi_p$-line $\Delta_z$ passing through $z$, and there exist two points $z', z''$ such that $\|z'\|_\infty = \|z''\|_\infty = r$ in $B_\infty(z, r] \cap \Delta_z$, with $z$ belonging to the $\varphi_p$-segment $]z', z''[$.

For all $i \notin I$, we have $\epsilon_i z_i > 0$, where $\epsilon_i$ is either $-1$ or $+1$. Therefore, there exist $u', u'' \in [z', z'']$ such that $\epsilon_i u'_i > 0$, $\epsilon_i u''_i > 0$, and $z \in ]u', u''[$. However, since $]u', u''[ \subset Co^p(A) \cap K$, this contradicts the assumption that $z$ is an extreme point.  
 $\Box$\\

It was shown in \cite{b19}, that the intermediate points have an algebraic representation in term of determinant.  Let $A=\{x^{(1)},...,x^{(m)}\}$ be a finite subset of $\Real^n$. If there is some
 $t_{_I,_J}^{(p)}\in \Delta_J^{(p)}$ such that $\zeta_{_I,_J}^{(p)}=\Gamma_J^{(p)}\big ({x^{(1)},...,x^{(m)},t_{_I,_J}^{(p)}}\big)=\stackrel{\varphi_p}{\sum\limits_{j\in J}}t_{_I,_J,j}^{(p)}\,x^{(j)}$ is an intermediate point of order $p$ for $A$  properly indexed on     $(I,J)$, then
\begin{equation}\label{transitlim}t_{_I,_J}^{(p)}=\stackrel{\varphi_p}{\sum_{j\in J}}\frac{\big|{\Lambda_{_I,_J}^{(\nu_J(j))} }\big |_p}
{\big |{\Lambda_{_I,_J}}\big |_p}\,e_{j}\, \quad \text{ and
}\quad \,\zeta_{_I,_J}^{(p)}=
\stackrel{\varphi_p}{\sum_{j\in J}}\frac{\big
|{\Lambda_{_I,_J}^{(\nu_J(j))}}\big |_p}{\big |\Lambda_{_I,_J}\big
|_p}x^{(j)}.\end{equation}
Any initial point of the form $x^{(j)}$ will be assumed to be properly indexed on $(\emptyset, \{j\})$ and denoted $\zeta_{\emptyset, \{j\}}$.

In the following, we say that $\zeta_{_I,_J}\in\Real^n$ is an {\bf
 intermediate point for $A$ in limit indexed on     $(I,J)$}, if there exists a
natural number $p_{_I,_J}$ and a sequence of  
 intermediate points $\big\{\zeta_{_I,_J}^{(p)}\big\}_{p\geq p_{_I,_J}}$ properly indexed on     $(I,J)$,
with $t_{_I,_J}^{(p)}\in \Delta_J^{(p)}$ for all $p\geq p_{_I,_J}$
such that:
\begin{equation}\zeta_{_I,_J}=\lim_{p\longrightarrow \infty}\zeta_{_I,_J}^{(p)}=\lim_{p\longrightarrow \infty}\Gamma_{J}^{(p)}\big ({x^{(1)},...,x^{(m)},t_{_I,_J}^{(p)}}\big).\end{equation}
Let us denote $\mathcal Z^\infty(A)$ the set of all the
 intermediate points in limit.

Let us denote as $\mathfrak S_{n_J}$ the set of all the permutations involved in the computation of the determinants $|{\Lambda_{_I,_J}^{(\nu_J(j))}}\big |_p$ for all $j\in J$. The next result is then derived from equation \eqref{transitlim}. In particular note that:

\begin{equation}\label{Declim}\zeta_{_I,_J}^{(p)}=\Gamma_J^{(p)}({x^{(1)},...,x^{(m)},t_{_I,_J}^{(p)}})=\frac{1}{\big |\Lambda_{_I,_J}\big
|_p} \stackrel{\varphi_p}{\sum_{j\in J}}\stackrel{\varphi_p}{\sum_{\sigma\in \mathfrak S_{n_J} }}  \mathrm{sign}(\sigma)\prod_{i\in [n_J]}\lambda_{i, \sigma(i)}^{(\nu_J(j))} x^{(j)}.\end{equation}

 One can see  {from} equation \eqref{Declim} that each term of the $\varphi_p$-generalized sum is independent of $p$, the limit of $\big |\Lambda_{_I,_J}\big
|_p$ being computed independently. An important property established in \cite{b20} is that for all $x\in \Real^n$, we have: 
\begin{equation}
\bigboxplus_{i\in [n]}x_i=0\quad \Longrightarrow \quad \stackrel{\varphi_p}{\sum_{i\in [n]}}x_i=0\quad \text{for all }  p\in \mathbb N.
\end{equation}
Equivalently, if there is some $p\in \mathbb N$ such that $\stackrel{\varphi_p}{\sum_{i\in [n]}}x_i\not=0$ then this implies that  $\bigboxplus_{i\in [n]}x_i\not=0$. 
\begin{lem}\label{algform}Let $A=\{x
^{(1)},...,x^{(m)}\}$ be a finite subset of $\Real^n$. Let
$(I,J)\subset [n]\times [m]$, and suppose that there exists some
natural number $p$ and some $t_{_I,_J}^{(p)}\in \Delta_J^{(p)}$ such
that $\zeta_{_I,_J}^{(p)}$
is a    intermediate point of order $p$ properly indexed on     $(I,J)$.   Then:\\
$(a)$ We have $|\Lambda_{_I,_J}|_\infty\not=0$;\\
$(b)$ There exists a natural number $p_{_I,_J}$ such that  $\Lambda_{_I,_J}$ is $\varphi_p$
invertible for all $p\geq p_{_I,_J}$;\\
$(c)$ There exists   an intermediate point in limit $\zeta_{_I,_J}$ indexed on  $(I,J)$
such that:
\begin{align*} \zeta_{_I,_J} =\lim_{p\longrightarrow \infty}\stackrel{\varphi_p}
{\sum_{j\in J}}\frac{\big | \Lambda_{_I,_J}^{(\nu_J(j))}\big
|_p}{\big |\Lambda_{_I,_J}\big |_p}x^{(j)}
 =\bigboxplus_{\substack{j\in
J\\\sigma\in \mathfrak S_{n_J}} }\alpha_{\sigma,
j}x^{(j)},\end{align*} where  $$\alpha_{\sigma,j}=\frac{1}{\big |\Lambda_{_I,_J}\big |_\infty}\sgn(\sigma)
 \prod_{i\in [n_J]}\lambda^{(\nu_J(j))}_{i,  \sigma(i) },
 \quad {\bigboxplus_{\substack{j\in
J\\\sigma\in \mathfrak S_{n_J}} }}\alpha_{\sigma, j}=1$$ and
$\bigboxplus_{\sigma\in \mathfrak S_{n_J} }\alpha_{\sigma, j}\geq 0$, for all $j$.
\end{lem}
\begin{expl} {\small Suppose that $n=3$ and $m=3$. Let $A=\{x^{(1)}, x^{(2)}, x^{(3)}  \}$
with $x^{(1)}=(-2,2,3)$, $x^{(2)}=( 3,-1,-4)$, 
and $x^{(3)}=(-1,-3,-1) $.   We have
{\footnotesize
$$\Gamma^{(p)}_{\{1,2,3\}}(x^{(1)}, x^{(2)}, x^{(3)};
t_1,t_2,t_3)=\stackrel{\varphi_p}{\sum_{j\in \{1,2,4\}}}t_ix_i=t_1
\begin{pmatrix}-2\\ 2\\3\end{pmatrix}
\stackrel{p}{+}t_2\begin{pmatrix} 3\\-1\\-4\end{pmatrix}
\stackrel{p}{+}t_3\begin{pmatrix}-1\\-3\\-1\end{pmatrix}.
$$ }
Suppose that $I=\{1,2\}$. Then: {\footnotesize $$
\Lambda_{\{1,2\},\{1,2,3\}}  =\begin{pmatrix}-2& 3& -1\\
2&-1&-3\\1&1&1 \end{pmatrix};
 \quad \Lambda_{\{1,2\},\{1,2,3\}}^{(1)}  =\begin{pmatrix}0& 3&-1\\0&-1&-3\\
1&1&1\end{pmatrix};$$
$$\Lambda_{\{1,2\},\{1,2,3\}}^{(2)}
=\begin{pmatrix}-2& 0&-1\\2&0&-3\\ 1&1&1\end{pmatrix};
 \quad \Lambda_{\{1,2\},\{1,2,3\}}^{(3)}  =\begin{pmatrix}-2& 3&0\\ 2&-1&0\\ 1&1&1\end{pmatrix} .$$}
The Cramer solution of the $\varphi_p$-linear system
$$\Lambda_{\{1,2\},\{1,2,3\}} \stackrel{p}{\cdot} \begin{pmatrix}s_1\\s_2\\s_3\end{pmatrix}=\begin{pmatrix}0 \\ 0\\ 1\end{pmatrix}$$
yields {\footnotesize
$$t_{\{1,2\},\{1,2,3\},1}^{(p)}=s_1^\star=\frac{\left|\Lambda_{\{1,2\},\{1,2,3\}}^{(1)}
\right|_p}{\left|\Lambda_{\{1,2\},\{1,2,3\}} \right|_p}=\frac{
\big(-9^{2p+1}- 1\big )^{\frac{1}{2p+1}}}{\big( -9^{2p+1}
-2^{2p+1}-1\big)^{\frac{1}{2p+1}}}.$$
$$t_{\{1,2\},\{1,2,3\},2}^{(p)}=s_2^\star=\frac{\left|\Lambda_{\{1,2\},\{1,2,3\}}^{(2)} \right|_p}{\left|\Lambda_{\{1,2\},\{1,2,3\}} \right|_p}=\frac{ \big(-2^{2p+1}-  6^{2p+1}\big)^{\frac{1}{2p+1}}}{\big( -9^{2p+1} -2^{2p+1}-1\big)^{\frac{1}{2p+1}}}.$$
$$t_{{\{1,2\},\{1,2,3\}},3}^{(p)}=s_3^\star=\frac{\left|\Lambda_{\{1,2\},\{1,2,3\}}^{(3)} \right|_p}{\left|\Lambda_{\{1,2\},\{1,2,3\}} \right|_p}=\frac{ \big(2^{2p+1}-  6^{2p+1}\big )^{\frac{1}{2p+1}}}{\big( -9^{2p+1}
-2^{2p+1}-1\big)^{\frac{1}{2p+1}}}.$$} One can check that for all
$p\in \mathbb N$, $t^{(p)}_{\{1,2\},\{1,2,4\}}\geq 0$. Moreover
$\lim_{p\longrightarrow \infty }t_{\{1,2\},\{1,2,3\},1}^{(p)}=1$,
$\lim_{p\longrightarrow \infty }t_{\{1,2\},\{1,2,3\},2}^{(p)}=
\lim_{p\longrightarrow \infty
}t_{\{1,2\},\{1,2,3\},3}^{(p)}=\frac{2}{3}$. The corresponding $\varphi_p$ intermediate point indexed on $(I,J)$ is then:
{\scriptsize \begin{align*} &\Gamma^{(p)}_{\{1,2,4\}}(x^{(1)}, x^{(2)}, x^{(3)};
t_1^{(p)},t_2^{(p)},t_3^{(p)})\\&=\frac{
\big(-9^{2p+1}- 1\big )^{\frac{1}{2p+1}}}{\big( -9^{2p+1}
-2^{2p+1}-1\big)^{\frac{1}{2p+1}}}x^{(1)}\stackrel{p}{+}\frac{ \big(-2^{2p+1}-  6^{2p+1}\big)^{\frac{1}{2p+1}}}{\big( -9^{2p+1} -2^{2p+1}-1\big)^{\frac{1}{2p+1}}}x^{(2)}\stackrel{p}{+}\frac{ \big(2^{2p+1}-  6^{2p+1}\big )^{\frac{1}{2p+1}}}{\big( -9^{2p+1}
-2^{2p+1}-1\big)^{\frac{1}{2p+1}}}x^{3}\end{align*}}
Taking the limit yields:
{\scriptsize \begin{align*}  \zeta_{\{1,2\},\{1,2,3\}} &=\frac{
 (-9)}{  (-9)\boxplus (-2)}x^{(1)}\boxplus \frac{
(- 1)}{  (-9)\boxplus (-2)}x^{(1)}\boxplus\frac{  (-2) }{  (-9)\boxplus (-2)} x^{(2)}\boxplus \frac{    (-  6)}{  (-9)\boxplus (-2)} x^{(2)}\boxplus \frac{  2 }{ (-9)\boxplus (-2)}x^{(3)}\boxplus \frac{    (-  6)}{  (-9)\boxplus (-2)} x^{(3)}\\& =\frac{
   9  }{  9\boxplus 2 }x^{(1)}\boxplus \frac{
     1}{  9\boxplus 2 }x^{(1)}\boxplus\frac{   2  }{  9\boxplus 2 } x^{(2)}\boxplus\frac{     6 }{  9\boxplus 2 } x^{(2)}\boxplus \frac{  6 }{ 9 \boxplus 2}x^{(3)}\boxplus \frac{   -  2 }{ 9 \boxplus 2}x^{(3)}\\&=\frac{
    1}{  9\boxplus 2 }\Big[ 9 x^{(1)} \boxplus    x^{(1)}\boxplus   2 x^{(2)} \boxplus   6   x^{(2)}\boxplus   6 x^{(3)} \boxplus (-  2) x^{(3)}\Big]=x^{(1)}\boxplus \frac{1}{9}x^{(1)}\boxplus \frac{2}{9} x^{(2)}\boxplus \frac{2}{3}x^{(2)}\boxplus \frac{2}{3} x^{(3)}\boxplus (-\frac{2}{9})x^{(3)}.\end{align*}}
 The complete calculation yields:
    { \scriptsize \begin{align*}  \zeta_{\{1,2\},\{1,2,3\}}&=\frac{
   9  }{  9\boxplus 2 }\begin{pmatrix}-2\\2\\3   
   \end{pmatrix}\boxplus \frac{
     1}{  9\boxplus 2 }\begin{pmatrix}-2\\2\\3   
   \end{pmatrix}\boxplus\frac{   2  }{  9\boxplus 2 } \begin{pmatrix}3\\-1\\-4   
   \end{pmatrix}\boxplus\frac{     6 }{  9\boxplus 2 } \begin{pmatrix}3\\-1\\-4   
   \end{pmatrix}\boxplus \frac{  6 }{ 9 \boxplus 2}\begin{pmatrix}-1\\-3\\-1   
   \end{pmatrix}\boxplus \frac{   -  2 }{ 9 \boxplus 2}\begin{pmatrix}-1\\-3\\-1   
   \end{pmatrix}\\
   &=\frac{1}{9}\begin{pmatrix} (-18)\boxplus(-2) \boxplus 6\boxplus 18  \boxplus(- 6) \boxplus 2\\
    18 \boxplus 2 \boxplus (-2)\boxplus (-6)  \boxplus(- 18) \boxplus 6\\
     27 \boxplus 3 \boxplus (-8)\boxplus (-24)  \boxplus(- 6) \boxplus 2\\
   \end{pmatrix}=\begin{pmatrix}0\\0\\{3}   
   \end{pmatrix}.
   \end{align*}} 
   Clearly, we have found $t_{1,1}=1, t_{1,2}=\frac{2}{9},t_{2,1}=\frac{2}{9},t_{2,2}=\frac{2}{3},t_{3,1}=\frac{2}{3}, t_{3,2}=-\frac{2}{9}$,
   with $t_{1,1}\boxplus t_{1,2}=1\geq 0$,  $t_{2,1}\boxplus t_{2,2}=\frac{2}{3}\geq 0$ and $t_{3,1}\boxplus t_{3,2}=\frac{2}{3}\geq 0$; moreover $t_{1,1}\boxplus t_{1,2}\boxplus t_{2,1}\boxplus t_{2,2}\boxplus t_{3,1}\boxplus t_{3,2}=1.$ }
\end{expl}

In the following, for all finite subsets $A$ of $\Real^n$, let us
denote $\mathcal K(A)$ the collection of all the $n$-dimensional
orthant $K$ such that $Co^\infty(A)\cap K\not=\emptyset$.  Let $A=\{x^{(1)},...,x^{(m)}\}$ be a finite subset of $K$. For any natural number $p$ let $A^{(p)}=\{x ^{(1,p)},...,x ^{(m,p)}\}$ be a
finite collection of $m$ vectors in $K$. It was shown in \cite{b15} that if for $j=1,...,m$, $\lim_{p\longrightarrow \infty}x^{(j,p)}=x^{(j)}$, then \begin{equation}
\Lim_{p\to\infty}Co^{p}(A^{(p)}) =\Big\{\bigboxplus_{j\in [m]}t_jx^{(j)}:
\max_j t_j=1, t_j\geq 0\Big\}=\mathbb B[A].\end{equation}
The next result is  {quite} immediate.  
 
 \begin{prop}Let $A=\{x
^{(1)},...,x^{(m)}\}$ be a finite subset of $\Real^n$. Suppose that $Co^\infty(A)\cap K$ is a nonemty set. Suppose that there is an infinite set of natural numbers $p$ such that  $Co^{p}(A)\cap K\not=\emptyset$.\\
$(a)$ There exists a collection $\mathcal C$ of pairs $(I,J)$ independent of $p$, and an increasing sequence of natural numbers $\{p_k\}_{k\in \mathbb N}$ such that for all $k$
 $$Co^{p_k}(A)\cap K=Co^{p_k}\Big(\big\{\zeta_{_I,_J}^{(p_k)}:(I,J)\in \mathcal C\big \}\Big)$$
 where for each $(I,J)$, $\zeta_{_I,_J}^{(p_k)}$ is an intermediate point of order $p_k$ properly indexed on  $(I,J)$.\\
\noindent $(b)$ There exists a collection $\mathcal C$ of pairs $(I,J)$ such that $$\Ls_{p\longrightarrow \infty}\big[Co^p(A)\cap K\big]\subset \Big\{\bigboxplus_{(I,J)\in \mathcal C} t_{ _I,_J }\zeta_{_I,_J}: \max_{ (I,J) \in \mathcal C}t_{ _I,_J }=1, t_{ _I,_J} \geq 0\Big\}, $$
with $\zeta_{_I,_J}\in K$ for all $(I,J)\in \mathcal C$.\\
 
\end{prop}
{\bf Proof:}  $(a)$ For all $p$ such that $Co^{p}(A)\cap K\not=\emptyset$, let us denote $\mathcal C^{(p)}$ the collection of pairs $(I,J)$ such that:
$$  Co^{p}(A)\cap K=\big\{\zeta_{_I,_J}^{(p)}: (I,J)\in \mathcal C^{(p)}\big\}, $$
where for each $(I,J)$, $\zeta_{_I,_J}^{(p)}$ is an intermediate point of order $p$ properly indexed. Since for all $p$, $ \mathcal C^{(p)}$ is bounded by a finite number of pairs $(I,J)$, there is a collection $\mathcal C$   and a sequence $ \{ p_{k} \}_{k\in \mathbb N}$ such that for all $k$
$$\mathcal C^{(p_{k })}=\mathcal C. $$
 Therefore, for all $k$
 $$Co^{p_k}(A)\cap K=Co^{p_k}\Big(\big\{\zeta_{_I,_J}^{(p_k)}:(I,J)\in \mathcal C\big \}\Big).$$
$(b)$ Suppose that $x\in \Ls_{p\longrightarrow \infty }[Co^p(A)\cap K]$ there is a sequence $\{p_{q}\}_{q\in \mathbb N}$ and a sequence 
$\{x^{(p_ q )}\}_{q\in \mathbb N}$ with $x^{(p_{q})}\in Co^{p_q}(A)\cap K$ for all $q$ and such that $x\in\Ls_{q\longrightarrow \infty} x^{(p_{q})}$.   From $(a)$   there exists an increasing subsequence $\{p_{q_r}\}_{r\in \mathbb N}$ and a collection $\mathcal C$ of pairs $(I,J)$ such that for all $r$
$$\mathcal C^{(p_{q_r })}=\mathcal C. $$
It follows that for all $r$
$$x^{p_{q_r }}\in Co^{p_{q_r}}\Big(\big\{\zeta_{_I,_J}^{(p_{q_r})}:(I,J)\in \mathcal C\big \}\Big).$$
Now note that $|\Lambda_{_I,_J}|_p\not=0$ for some $p$ implies $|\Lambda_{_I,_J}|_\infty \not=0$. 
It follows that for all $(I,J)\in \mathcal C$
$$\lim_{r\longrightarrow \infty}\zeta_{_I,_J}^{(p_{q_r})}=\zeta_{_I,_J}. $$
Since for all $(I,J)$ the vectors $\zeta_{_I,_J}^{(p_{q_r})}$ are copositive, it follows that:
$$\Lim_{r\longrightarrow \infty}Co^{p_{q_r}}\Big(\big\{\zeta_{_I,_J}^{(p_{q_r})}:(I,J)\in \mathcal C\big \}\Big)=Co^{\infty}\Big(\big\{\zeta_{_I,_J}:(I,J)\in \mathcal C\big \}\Big).$$
Therefore $$x\in \Big\{\bigboxplus_{(I,J)\in \mathcal C} t_{_{_I,_J}}\zeta_{_I,_J}: \max_{_{(I,J)}\in \mathcal C}t_{_{_I,_J}}=1, t_{_{_I,_J}}\geq 0\Big\}. $$
Since for any $x\in \Ls_{p\longrightarrow \infty}Co^p(A)$ such a collection $\mathcal C$   exists  and since the set of all the possible pairs $(I,J)$ is finite, the result follows immediately. $\Box$\\

\begin{prop}\label{fundinclus}Let $A=\{x ^{(1)},...,x^{(m)}\}$ be a finite subset of $\Real^n$. Then there exists $m$ subsets of $\mathbb N$, $I_1, I_2...I_m$ such that
$$Co^\infty (A)\subset \Big\{\bigboxplus_{\substack{ k_j\in I_j\\j\in [m]}}t_{k_j}x^{(j)}: \bigboxplus_{k_j\in I_j}t_{k_j}\geq 0,\max_{\substack{  k_j\in I_j\\j\in [m]}} t_{ k_j}=1,  j\in [m] \Big\}. $$

\end{prop}
{\bf Proof:} Suppose that $x\in Co^\infty(A)$, there exists an  $n$-dimensional orthant $K$ such that $x\in \Ls_{p\longrightarrow \infty}[Co^p(A)\cap K]$. Therefore there exists a collection $\mathcal C$ of pairs $(I,J)$ such that

$$x\in \Big\{\bigboxplus_{ (I,J) \in  \mathcal C}t_{ _I,_J } \zeta_{_I,_J}:  \max_{_{(I,J)\in \mathcal K}}t_{ _I,_J } =1, t_{ _I,_J }\geq 0\big\} .$$

From Lemma \ref{algform}, for any pair $(I,J)\in \mathcal C$ an $(I ,J )$-intermediate point in limit $\zeta_{I ,J }$ is defined as:
\begin{align*} \zeta_{_I ,_J } =\lim_{p\longrightarrow \infty}\stackrel{\varphi_p}
{\sum_{j\in J }}\frac{\big | \Lambda_{_I  ,_J }^{(\nu_{J }(j))}\big
|_p}{\big |\Lambda_{_I,_J }\big |_p}x^{(j)}
 =\bigboxplus_{\substack{j\in
J\\\sigma\in \mathfrak S_{n_J}} }\alpha_{ \sigma,
j}x^{(j)},\end{align*} where  $$\alpha_{ \sigma,j}=\frac{1}{\big |\Lambda_{_I  ,_J }\big |_\infty}\sgn(\sigma)
 \prod_{i\in [n_J]}\lambda^{(\nu_{J }(j))}_{i,  \sigma(i) },
 \quad {\bigboxplus_{\substack{j\in
J \\\sigma\in \mathfrak S_{n_J}} }}\alpha_{ \sigma, j}=1$$ and
 $\bigboxplus_{\sigma\in \mathfrak S_{n_J} }\alpha_{\sigma, j}\geq 0$, for all $j$.  For all  $(I,J)\in \mathcal C$ it is easy to check $|\alpha_{\sigma,j}|\leq 1$. Therefore, the condition  $\bigboxplus\limits_{\substack{\sigma\in \mathfrak S_{n_J}\\ j\in J }}\alpha_{\sigma,j}=1$ implies that
 $\max\limits_{\substack{\sigma \in \mathfrak S_{n_J}\\j\in  J }}\alpha_{\sigma,j}=1$.
 
Let us  denote $\ell_I=(n_{I }+1) !$ and consider a bijective index map $i:[\ell_I] \longrightarrow \mathfrak S_{n_J} $  where $n_{I }=\Card\, I $, which associates to any number in $[\ell_I]$ a permutation $\sigma \in  \mathfrak S_{n_J}$. For all $j\in J$, and all 
$q\in [\ell_I]$, let us denote  $\alpha_{ \sigma,j}=s_{   i(\sigma ), j}=s_{  q,j}$. By construction, it follows that 
$$\zeta_{I ,J }=\bigboxplus_{\substack{j\in
J \\q\in [\ell_I] }}s_{  q,j}x^{(j)}$$
with 
$$\max\limits_{\substack{q\in [\ell_I]\\j\in J }}s_{  q,j}=1\quad \text{ and }\quad \bigboxplus\limits_{q\in [\ell_I]}s_{  q,j}
\geq 0\quad (A)$$ for all $j$. However, by hypothesis $x$ can be written as:
$$x=\bigboxplus_{{ (    I ,J )}\in \mathcal C}r_{_I,_J}\zeta_{_I,_J},$$
with $r_{_I,_J}\geq 0$ for all $(I,J)$ and $\max_{(_I,_J)\in \mathcal C}r_{_I,_J}=1.$ Therefore $x$ can be rewritten
$$x=\bigboxplus_{ (I,J) \in \mathcal C}r_{_I,_J}\Big(\bigboxplus_{\substack{j\in
J  \\q\in [\ell_{I}] }}s_{q,j}x^{(j)}\Big).$$
Since all the intermediate points in $K$ are copositive, we deduce that:
$$x=\bigboxplus_{\substack{j\in
J \\q\in [\ell_{I}]\\ (I,J) \in \mathcal C }}r_{_{_I,_J}}s_{  q,j}x^{(j)} .$$
Let us denote $\mathcal C_j=\{(I,J)\in \mathcal C: j\in J\}$. We equivalently have:
$$x=\bigboxplus_{\substack{ q\in [\ell_{I}]\\(I,J) \in \mathcal C_j\\ j\in [m] }}r_{_{_I,_J}}s_{  q,j}x^{(j)} .$$

For any $j$, the set  $$\Delta_j=\Big\{(I,J, q ): q\in [\ell_I], (I,J)\in \mathcal C_j\Big\}$$ is isomorphic to a finite subset $I_j$ of $\mathbb N$. Therefore, for all $(I,J, q )\in \Delta_j$, there is some $k_j\in I_j$ and a real number $t_{k_j}$ such that 
$$t_{k_j}=r_{_{_I,_J}}s_{ q,j }. $$  

Moreover, we have $\max\limits_{(I,J)\in \mathcal C}r_{_{_I,_J}}=1$ with $r_{_{_I,_J}}\geq 0$, hence:
 $$\max\limits_{(I,J)\in \mathcal C}r_{_{_I,_J}} =\max\limits_{\substack{(I,J)\in \mathcal C_j\\ j\in [m]}}r_{_{_I,_J}}=1.$$  Thus,    there exists $(I_0, J_0)\in \mathcal C$ such that $r_{_{I_0,J_0}}=1$. From $(A)$  there exists $j_0\in J_0$ and $q_0\in [\ell_{I_0}]$ such that 
 $$\max\limits_{\substack{q \in [\ell_{I_0}]\\j\in  J_0 }}s_{q,j}=1=s_{q_0,j_0}.$$
Hence there is some $I_0,J_0,q_0, j_0$ such that $r_{_{I_0,J_0}}=1$ and $s_{k_0,j_0}=1$. Thus:
$$r_{I_0,J_0}s_{k_0,j_0}=1=\max_{\substack{(I,J,q)\in \Delta_{j}\\j\in [m]}}r_{_{_I,_J}}s_{ q,j }=\max_{\substack{k_j\in I_j\\  j\in [m]}}t_{k_j}=1.$$ 

Moreover since $r_{_{_I,_J}}\geq 0$   and $\bigboxplus\limits_{q\in [\ell_I]}s_{ q,j}\geq 0$ we have $\bigboxplus\limits_{k_j\in I_j}t_{k_j}=\bigboxplus\limits_{\substack{q\in [\ell_I]\\(I,J)\in \mathcal C_j}}r_{_I,_J}s_{ q,j}\geq 0$. Since this is true for all $x\in Co^\infty(A)$, and any $n$-dimensional orthant $K$ we deduce the result. $\Box$\\

In subsection \ref{exampleimply} we give an example of decomposition of a $\mathbb B$-polytope in a three-dimensional space. This example is depicted in Figure \ref{Polalg}.2.

\subsection{Converse Inclusion and Algebraic Formulation of a Polytope}
In this section, we complete Proposition~\ref{fundinclus} by providing an algebraic formulation of the limit polytope $Co^\infty(A)$. The central intuition is as follows.

Let $n_1, n_2, \dots, n_m$ be positive integers, and let $A = \{x^{(j)}\}_{j \in [m]}$ be a finite set in $\mathbb{R}^n$. The standard convex hull of $A$ can be equivalently expressed as:
\begin{equation}\label{initform}
Co(A) = \Big\{ \sum_{\substack{k \in [n_j] \\ j \in [m]}} t_{j,k} x^{(j)} \;:\; \sum_{k \in [n_j]} t_{j,k} \geq 0,\; \sum_{\substack{k \in [n_j] \\ j \in [m]}} t_{j,k} = 1,\; j \in [m] \Big\}.
\end{equation}

Setting \( s_j = \sum_{k \in [n_j]} t_{j,k} \), it is straightforward to see that any element \( z \) of the set defined above can be written as
\[
z = \sum_{j \in [m]} s_j x^{(j)},
\]
where \( s_j \geq 0 \) and \( \sum_{j \in [m]} s_j = 1 \). Note that the number of times each \( x^{(j)} \) appears in the expression does not affect the resulting convex combination. However, this invariance no longer holds when replacing the standard addition with the \(\mathbb{B}\)ox-plus operation $\boxplus$, due to its non-associativity. Since the expression above does not depend on the specific values of the $n_j$, we can generalize it further by allowing all possible choices of positive integers $n_j$:
\begin{equation}\label{initform2}
Co(A) = \Big\{ \sum_{\substack{k \in [n_j] \\ j \in [m]}} t_{j,k} x^{(j)} \;:\; \sum_{k \in [n_j]} t_{j,k} \geq 0,\; \sum_{\substack{k \in [n_j] \\ j \in [m]}} t_{j,k} = 1,\; n_j \in \mathbb{N},\; j \in [m] \Big\}.
\end{equation}

We now formally substitute $\sum \mapsto \bigboxplus$ in this expression. For any finite subset $A = \{x^{(j)}\}_{j \in [m]}$ of $\mathbb{R}^n$, we define:
\begin{equation}
Co^\boxplus[A] = \Big\{ \bigboxplus_{\substack{k \in [n_j] \\ j \in [m]}} t_{j,k} x^{(j)} \;:\; \bigboxplus_{k \in [n_j]} t_{j,k} \geq 0,\; \max_{\substack{k \in [n_j] \\ j \in [m]}} t_{j,k} = 1,\; n_j \in \mathbb{N},\; j \in [m] \Big\}.
\end{equation}

Let $\langle \mathbb{N} \rangle$ denote the set of all finite subsets of $\mathbb{N}$. By replacing each $[n_j]$ with a finite subset $I_j \in \langle \mathbb{N} \rangle$, we obtain the equivalent formulation:
\begin{equation}
Co^\boxplus[A] = \Big\{ \bigboxplus_{\substack{k_j \in I_j \\ j \in [m]}} t_{k_j} x^{(j)} \;:\; \bigboxplus_{k_j \in I_j} t_{k_j} \geq 0,\; \max_{\substack{k_j \in I_j \\ j \in [m]}} t_{k_j} = 1,\; I_j \in \langle \mathbb{N} \rangle,\; j \in [m] \Big\}.
\end{equation}

This representation encompasses, as a special case, the formulation given in Proposition~\ref{fundinclus}, with the key difference that the index sets $I_j$ are no longer fixed a priori.

 We first establish the following property. 
\begin{prop}\label{InclusLim}For all finite subsets $A =\{x ^{(1)},...,x^{(m)}\}$  of $\Real^n$, we have:
$$Co^\boxplus[A]\subset \Linf_{p\longrightarrow \infty}Co^{(p)}(A). $$
\end{prop}
{\bf Proof:}  Suppose that $z=\bigboxplus\limits_{\substack{ k_j\in I_j\\j\in [m]}}t_{k_j}x^{(j)}\in Co^\boxplus [A]$ which implies that $\bigboxplus\limits_{k_j\in I_j}t_{ k_j}\geq 0$, $\max\limits_{\substack{ k_j\in I_j\\j\in [m]}} t_{ k_j}=1$ for all $I_j\in \langle \mathbb N \rangle, j\in [m]$. First, note that we have for all natural numbers $p$
\begin{equation} \stackrel{\varphi_p}{\sum_{\substack{ k_j\in I_j\\j\in [m]}}}t_{k_j}x^{(j)}=\stackrel{\varphi_p}{\sum_{\substack{ j\in [m]}}}\stackrel{\varphi_p}{\sum_{\substack{ k_j\in I_j }}}t_{ k_j}x^{(j)}. \end{equation}
Let us denote  $t^{(p)}_j=\stackrel{\varphi_p}{\sum\limits_{\substack{ k_j\in I_j }}}t_{k_j}$. It follows that 
$$\stackrel{\varphi_p}{\sum_{\substack{ k_j\in I_j\\j\in [m]}}}t_{k_j}x^{(j)}=\stackrel{\varphi_p}{\sum_{\substack{ j\in [m]}}}t_{j}^{(p)}x^{(j)}.$$
Suppose that for any $j$, $\bigboxplus\limits_{k_j\in I_j}t_{k_j}\geq 0$. If $\bigboxplus\limits_{k_j\in I_j}t_{k_j}= 0$, then for all $p$, we have $\stackrel{\varphi_p}{\sum\limits_{\substack{ k_j\in I_j }}}t_{k_j}=0$ which implies that $t_j^{(p)}=0. $ If $\bigboxplus\limits_{k_j\in I_j}t_{k_j}> 0$, then there exists some $p^{(j)}$ such that for all $p\geq p^{(j)}$, we have $t_j^{(p)}=\stackrel{\varphi_p}{\sum\limits_{\substack{ k_j\in I_j }}}t_{k_j}>0$. Therefore, for all $p\geq \max_{j\in [m]}\{p^{(j)}\}$ we have $t_j^{(p)}\geq 0$. Let us consider the sequence $\{z^{(p)}\}_{p\in \mathbb N}$ defined as
$$z^{(p)}= \stackrel{\varphi_p}{\sum_{\substack{ j\in [m]}}}s_{j}^{(p)}x^{(j)},$$
where for any $j$ and all $p\geq \max_{j\in [m]}\{p^{(j)}\}$
$$s_{j}^{(p)}=\frac{t_{j}^{(p)}}{ \stackrel{\varphi_p}{\sum_{\substack{ i\in [m]}}}t_{i}^{(p)} }.$$
We obviously have: $   \stackrel{\varphi_p}{\sum_{\substack{ j\in [m]}}}s_{j}^{(p)}=1$. Moreover  $t_i^{(p)}\geq 0$ for all $i$ implies that  $s_j^{(p)}\geq 0$.
Consequently, for all $p$, $z^{(p)}\in Co^{(p)}(A).$ However, since $\max_{\substack{ k_j\in I_j\\j\in [m]}} t_{k_j}=1$, we deduce that $\max\limits_{\substack{  j\in [m]}} \big(\max\limits_{\substack{ k_j\in I_j}} t_{k_j}\big)=1$. Now, note that   $\bigboxplus_{k_j\in I_j}t_{k_j}\geq 0$ implies that $\max\limits_{\substack{ k_j\in I_j}} t_{k_j}\geq 0$ and we deduce that: 
$$\lim_{p\longrightarrow \infty}\stackrel{\varphi_p}{\sum_{\substack{ j\in [m]}}}t_{j}^{(p)}=\max\limits_{\substack{  j\in [m]}}\big(\max\limits_{\substack{ k_j\in I_j}} t_{k_j}\big)=1.$$
Therefore $$ \lim_{p\longrightarrow \infty}s_j^{(p)}= \lim_{p\longrightarrow \infty}t_j^{(p)}=\lim_{p\longrightarrow \infty}\stackrel{\varphi_p}{\sum\limits_{\substack{ k_j\in I_j }}}t_{k_j}=\bigboxplus_{k_j\in I_j}t_{k_j}.$$
Moreover:
$$\lim_{p\longrightarrow \infty} \stackrel{\varphi_p}{\sum_{\substack{ k_j\in I_j\\j\in [m]}}}t_{k_j}x^{(j)}=\bigboxplus_{\substack{ k_j\in I_j\\j\in [m]}}t_{k_j}x^{(j)}.$$
Hence $z =\lim_{p\longrightarrow \infty}z^{(p)}\in \Linf\limits_{p\longrightarrow \infty } Co^{(p)}(A)$. Thus, we deduce that $Co^{\boxplus}[A]\subset \Linf\limits_{p\longrightarrow \infty } Co^{(p)}(A)\subset Co^\infty(A)$,
which ends the proof. $\Box$\\

\begin{prop} \label{EquivAlg}For all finite subsets $A =\{x ^{(1)},...,x^{(m)}\}$ of $\Real^n$,
$$Co^\infty(A)=Co^\boxplus[A]= \Big\{\bigboxplus_{\substack{ k_j\in I_j\\j\in [m]}}t_{k_j}x^{(j)}: \bigboxplus_{k_j\in I_j}t_{k_j}\geq 0,\max_{\substack{ k_j\in I_j\\j\in [m]}} t_{k_j}=1, I_j\in \langle N\rangle ,  j\in [m] \Big\}. $$
\end{prop}
{\bf Proof:} From Proposition \ref{fundinclus}, we have shown that   there exists $m$ subsets of $\mathbb N$, $I_1, I_2...I_m$ such that
$$Co^\infty (A)\subset\Big\{\bigboxplus_{\substack{ k_j\in I_j\\j\in [m]}}t_{k_j}x^{(j)}: \bigboxplus_{k_j\in I_j}t_{k_j}\geq 0,\max_{\substack{ k_j\in I_j\\j\in [m]}} t_{k_j}=1,  j\in [m] \Big\}. $$ 
Therefore, we deduce that $Co^\infty (A) \subset Co^\boxplus [A]. $ Since from Proposition \ref{InclusLim} we have $Co^\boxplus [A]\subset \Linf\limits_{p\longrightarrow \infty } Co^{(p)}(A)\subset Co^\infty(A)$ the converse inclusion is immediate. $\Box$\\

It was shown in \cite{b19} that for all finite non-empty subsets $A$ of $\mathbb{R}^n$, the sequence $\{ Co^p(A) \}_{p \in \mathbb{N}}$ converges to $Co^\infty(A)$  with respect to the Hausdorff{-Pompeiu} metric. Proposition \ref{EquivAlg} provides an equivalent formulation of this result, since $Co^\infty(A)$ is compact.

In what follows, we show that, in the case of two points, this formula is equivalent to the one obtained in \cite{b15}.

\begin{lem}\label{AltformTwo} For all $x,y\in\Real^n$,
$$Co^\infty(x,y)=\Big\{\bigboxplus_{\substack{ i_z\in I_z\\z\in \{ x,y\}}}t_{i_z}z: \bigboxplus_{ i_z\in I_z }t_{i_z}\geq 0,\max_{\substack{ i_z\in I_z\\z\in \{x,y\}}} t_{i_z}=1,  I_x\in \langle \mathbb N \rangle, I_y\in \langle \mathbb N \rangle \Big\}.  $$
 
\end{lem}
{\bf Proof:} From Proposition \ref{InclusLim}:
$$\Big\{\bigboxplus_{\substack{ i_z\in I_z\\z\in \{ x,y\}}}t_{i_z}z: \bigboxplus_{i_z\in I_z}t_{i_z}\geq 0,\max_{\substack{ i_z\in I_k\\z\in \{x,y\}}} t_{i_z}=1,  I_x\in \langle \mathbb N \rangle, I_y\in \langle \mathbb N \rangle \Big\}\subset Co^\infty (x,y) . $$ 
However, we also have from \cite{b15} $$Co^\infty(x,y)=\Big\{t_xx\boxplus s_xx\boxplus t_yy\boxplus s_y y:\max\{t_x,s_x,t_y,s_y\}=1, t_x,s_x,t_y,s_y\geq 0\Big\}.$$ 
Thus, by definition, $Co^\infty(x,y)\subset Co^\boxplus[\{x,y\}]$, which ends the proof. $\Box$\\

\begin{rem}Lemma \ref{AltformTwo} is in fact equivalent to Proposition \ref{fond} (see \cite{b15}). This equivalence stems from the fact that, under certain coordinate symmetries, a negative weight on one vector can be substituted by a positive weight on its symmetric counterpart. Consider a subset \( I \subset [n] \) such that for each \( i \in I \), \( t_x x_i = -t_y y_i \), with \( t_x, t_y > 0 \). Assume that for some \( i \in I \),
$$
(t_x x \boxplus s_x x \boxplus t_y y \boxplus s_y y)_i = s_x x_i
\quad \text{with} \quad s_x < 0 \ \text{and} \ s_x > -1.
$$
Then symmetry implies \( x_i = -\frac{t_y}{t_x} y_i \), so
$s_x x_i = -s_x \frac{t_y}{t_x} y_i.$ If \( t_y = 1 \), we obtain \( s_x x_i = -\frac{s_x}{t_x} y_i \), where \( 0 < -\frac{s_x}{t_x} \leq 1 \); if \( 0 < t_y < 1 \), the same conclusion holds since \( 0 < -s_x \frac{t_y}{t_x} \leq 1 \).

In both cases, the negative weight on \( x \) can be replaced by a positive weight on \( y \), without altering components outside \( I \) or breaking symmetry. This argument also applies if 
$(t_x x \boxplus s_x x \boxplus t_y y \boxplus s_y y)_i = t_y y_i,$
or when more than three weights are involved on either \( x \) or \( y \).

However, the situation becomes more intricate when more than two points are considered, as the possible symmetries may not correspond to the same pair of vectors.\end{rem}

By definition a subset $C$ of $\Real^n$ is $\mathbb B$-convex if for any finite subset $A$ of $C$, $Co^\infty(A)\subset C$. The next result shows that $\mathbb{B}$-convexity (over the entire Euclidean vector space) can be redefined in algebraic terms as follows. In \cite{b15}, it was mentioned that the concept of convexity defined in \cite{bh} can be relaxed by using a definition that involves only two points, instead of any finite number. The advantage of the former is that it is formulated in algebraic terms, whereas the latter is expressed only in topological terms. The next proposition overcomes this difficulty by providing an algebraic definition of $\mathbb{B}$-convexity over $\mathbb{R}^n$ that involves an arbitrary number of points.

\begin{prop}\label{ReformPoly} A subset $C$ of $\Real^n$ is $\mathbb B$-convex if for all finite subsets $A=\{x^{(1)}, ...,x^{(m)}\}$ of $C$
\begin{equation}  Co^\boxplus[A]= \Big\{\bigboxplus_{\substack{ k\in [n_j]\\j\in [m]}}t_{j,k}x^{(j)}: \bigboxplus_{k\in [n_j]}t_{j,k}\geq 0, \max_{\substack{ k\in [n_j]\\j\in [m]}} t_{j,k}=1, n_j\in \mathbb N, j\in [m]   \Big\}\subset C.\end{equation}
\end{prop}
We have demonstrated that this definition encompasses, as a special case, the two-point definition. In the following remark, we show that it also recovers, as a special case, the description of a $\mathbb{B}$-polytope when it is a subset of $\mathbb{R}_+^n$ (see \cite{bh}).

\begin{rem} Suppose that $A$ is a finite subset of $\Real_+^n$. For all $j$ the condition $\bigboxplus\limits\limits_{k\in [n_j]}{t_{j,k}}\geq 0 $ ensures that one can find some $t_j=\bigboxplus\limits_{k\in [n_j]}{t_{j,k}}\geq 0$ such that 
\begin{equation}\bigboxplus_{ k\in [n_j] }t_{j,k}x^{(j)}=\bigvee_{j\in [m]}t_jx^{(j)}\end{equation}
Moreover $\max\limits_{\substack{ k\in [n_j]\\j\in [m]}} t_{j,k}=1$ implies that $\max\limits_{j\in [m]}t_j=1$. Therefore:
\begin{align*}  Co^\infty(A)=&\Big\{\bigboxplus_ {\substack{ k\in [n_j]\\j\in [m]}}t_{j,k}x^{(j)}: \bigboxplus_{k\in [n_j]}t_{j,k}\geq 0, \max_{\substack{ k\in [n_j]\\j\in [m]}} t_{j,k}=1, n_j\in \mathbb N, j\in [m]   \Big\}\\&=\Big\{\bigvee_{j\in [m]}t_jx^{(j)}: \max_{j\in [m]}t_j=1 , t_j\geq 0,j\in [m]\Big\}.\end{align*}

\end{rem}

Proposition \ref{ReformPoly}, yields to the following equivalent definition of $\mathbb B$-convexity over the whole Euclidean vector space. However, this definition is now purely algebraic  and expressed with respect to the non-associative $\mathbb B$oxplus operation $\boxplus$ and the real scalar mutiplication. It extends the notion of an idempotent symmetric convex set from a multiary perspective.

\begin{defn}A subset \( C \) of \( \mathbb{R}^n \) is  $  \mathbb{B} $-convex if  for any finite subset \( A = \{x^{(1)}, \dots, x^{(m)}\} \) of \( C \) and for all  
\( t \in \prod_{j \in [m]} \mathbb{R}^{n_j} \) with \( n_j \in \mathbb{N} \) for each \( j \), the conditions  
$$
\bigboxplus\limits_{k\in [n_j]} t_{j,k} \geq 0, \quad \max\limits_{\substack{k\in [n_j]\\ j\in [m]}} t_{j,k} = 1 \quad \text{ imply that }\quad 
 \bigboxplus_{\substack{k\in [n_j]\\ j\in [m]}} t_{j,k} x^{(j)} \in C.
$$
\end{defn}
 
The advantage of this formulation is that it is sufficiently general to extend the class of functions defined in \cite{adil}. Likewise, it would be interesting to explore a harmonic version of this convexity structure in order to generalize certain results established in \cite{adilYe}.  Note that the convexity formalism introduced in this paper is consistent with the concept of idempotent symmetric pseudo-fields and idempotent symmetric spaces defined in \cite{b24}. It is easy to see that an idempotent symmetric space is $\mathbb{B}$-convex. Notice also that a non-associative symmetrization of Max-Plus convex sets was proposed in \cite{b20, b24}, and similar investigations could be carried out. This is briefly  done in the following example. 

\begin{expl}We consider the case of the non-associative symmetrisation of the Max-Plus semi-module analysed in \cite{b20}. Let $\mathbb M=\Real\cup \{-\infty\}$ and let us denote $(\mathbb M, \oplus, \otimes)$ the Maslov’s semi-module where we replace the operations $\vee$ with $\oplus$ and $+$ with $\otimes$. More precisely one can define on $\mathbb M$ the operations $\oplus$ and $\otimes$ respectively as 
$\lambda\oplus \mu=\max\{\lambda ,\mu\}$ and $\lambda \otimes \mu=\lambda +\mu$, where  $-\infty$ is the neutral element of the operation $\oplus$.  Suppose now
that $\mu\in \Real_-$ and let us extend the logarithm function to
the whole set of real numbers. This we do by introducing the set
\begin{equation}\widetilde {\mathbb M}=\mathbb M \cup (\Real +\mathrm{i}\pi)\end{equation}
where $\mathrm{i}$ is the complex number such that $\mathrm{i}^2=-1$ and $\Real
+\mathrm{i}\pi=\{\lambda +\mathrm{i}\pi:\lambda\in \Real\}$. {Note that this approach is not an extension to the complex numbers. In fact, it designates a copy of the real numbers.  This  formalism     is however convenient for introducing the following extended logarithmic function to $\widetilde {\mathbb M}$ and transferring the algebraic structure of the real set. }Let $\psi_{\ln}
:{\mathbb R} \longrightarrow \widetilde {\mathbb M}$ be the map defined by:
\begin{equation}
\psi_{ \ln}(\lambda)=\left\{\begin{matrix}\ln(\lambda)&\text{ if
}\lambda >0\\-\infty&\text{ if }\lambda=0\\\ln(-\lambda)+\mathrm{i}\pi&\text{ if
}\lambda<0.
\end{matrix}\right.
\end{equation}
The map $\lambda\mapsto \psi_{ \ln}(\lambda)$  is an isomorphism from $\mathbb
\Real$ to $\widetilde {\mathbb M}$. Let $\psi_{ \exp}: \widetilde
{\mathbb M}\longrightarrow \mathbb \Real$ be its inverse. Notice that
$\psi_{ \ln}(-1)=\mathrm{i}\pi$. By
definition we have $\lambda\widetilde{\boxplus}
\mu=\psi_{\ln}\big(\psi_{\exp}(\lambda)\boxplus \psi_{\exp}(\mu)\big )$. Moreover, we have $\lambda \widetilde{\otimes}
\mu=\psi_{\ln}\big(\psi_{\exp}(\lambda)\otimes \psi_{\exp}(\mu)\big )$. It follows that $(\widetilde{\mathbb M}, \tilde \oplus, \tilde \otimes)$ is isomorphic to $(\Real, \boxplus, \cdot)$ and is also an idempotent pseudo-field with the neutral elements $0_{\tilde {\mathbb M}}=-\infty$ and $1_{\tilde {\mathbb M}}=0$. 

{In Chapter 3 of \cite{koloma} it was shown that Max-Plus algebra can also be viewed as a limit algebraic structure via 
 the dequantization principle. }
For all $\lambda =(\lambda_1,...,\lambda_n)\in \widetilde{\mathbb M}^n$, let  us  denote:
 \begin{equation}\widetilde{\bigboxplus _{i\in [n]}}\lambda_i=
\psi_{\ln}\big(\bigboxplus_{i\in [n]}\psi_{\exp}(\lambda_i)\big ).\end{equation}
From \cite{b20}, it follows that for all $(\lambda_1,...,\lambda_n)\in \widetilde{\mathbb M}^n$, we have
$\widetilde{\bigboxplus}_{i\in [n]} \lambda_i=\lim_{p\longrightarrow \infty}\frac{1}{2p+1}\psi_{\ln}\Big(\sum_{i\in [n]}\psi_{\exp}\big((2p+1)\lambda_i\big)\Big). $
Suppose that $\tilde A=\{z^{(1)}, ...,z^{ (m)} \}$ is a finite subset of $\widetilde{\mathbb M}^n$. Setting $s_{j,k}=\psi_{\mathrm{ln}}(t_{j,k})$, for each pair $(j,k)$ we obtain the limit polytope of $\widetilde{\mathbb M}$ defined as:
\begin{equation}
Co^{\widetilde{\boxplus}}[A] = \Big\{ \widetilde{\bigboxplus_{\substack{k \in [n_j] \\ j \in [m]}}} s_{j,k} \otimes z^{(j)} \;:\; \widetilde{\bigboxplus_{k \in [n_j]}} s_{j,k}\in  \Real \cup\{-\infty\},\; \max_{\substack{k \in [n_j] \\ j \in [m]}} s_{j,k} = 0,\; n_j \in \mathbb{N},\; j \in [m] \Big\},
\end{equation}
which involves a convex structure on  $\widetilde{\mathbb M}^n$.

\end{expl}

\subsection{An Example and some Implications}\label{exampleimply}

In this section, we present a detailed example illustrating the construction of a 3-dimensional simplex.  Specifically, we provide a counterexample demonstrating that for any finite part $A$,  $\mathrm{IS}[A] \nsubseteq Co^\infty(A)$. In other words, in general, a 
$\mathbb B$-polytope is not idempotent symmetric convex, i.e. it does not satisfy the condition of convexity at the limit for two points, as defined in  \cite{b15} (which represents a special case of 
$\mathbb B$-convexity).

It was established in \cite{b19} that if $A=\{x ^{(1)},...,x^{(m)}\}$ is a finite subset of $\Real^n$, then for all natural numbers $p$, and all $K\in \mathcal K(A)$, since $\mathcal Z^{(p)}(A)\cap K$ contains all the $\varphi_p$-extreme points of $Co^p(A)\cap K$, that we have $Co^p\big(\mathcal Z^{(p)}(A)\cap K\big)=Co^p(A)\cap K.$ Moreover
\begin{equation}Co^p(A)=\bigcup_{K\in \mathcal K(A)}Co^p\big (\mathcal Z^{(p)}(A)\cap K\,\big ),\end{equation}
and
\begin{equation}Co^\infty(A)=\Lim_{p\longrightarrow \infty}Co^p(A)
= \bigcup_{K\in \mathcal K(A)}Co^\infty\big(\mathcal Z^\infty(A)\cap K\,\big).\end{equation}
In the following the proposed decomposition uses the fact that when there is an $n$-dimensional orthant $K$ such that $A\subset K$ then $Co^\infty(A)=\Big\{\bigboxplus\limits_{j\in [m]}t_jx^{(j)}: \max_{j\in [m]}t_j=1, t_j\geq 0, j\in [m]\Big\}$.

\begin{expl}\label{fundex}  {\small Suppose that $n=3$ and lets us consider the vectors $x^{(1)}=(1,3,3)$, $x^{(2)}=(3,2,3)$ and $x^{(3)}=(-1,-1,-1)$. Let us calculate the intermediate points related to the $\mathbb B$-polytope $Co^\infty(A)$. From the results above, we have:
\begin{align*}Co^\infty(A) =\Bigg\{t_{1,1}\begin{pmatrix}
1\\3\\3
\end{pmatrix}&\boxplus \cdots \boxplus t_{1,n_1}\begin{pmatrix}
1\\3\\3
\end{pmatrix} \boxplus t_{2,1}\begin{pmatrix}
3\\2\\3
\end{pmatrix}\boxplus \cdots \boxplus t_{2,n_2}\begin{pmatrix}
3\\2\\3
\end{pmatrix}    \boxplus t_{3,1}\begin{pmatrix}
-1\\-1\\-1
\end{pmatrix}\boxplus \cdots \boxplus t_{2,n_3}\begin{pmatrix}
-1\\-1\\-1
\end{pmatrix}:\\& \bigboxplus_{k\in [n_1]}t_{1,k}\geq 0, \bigboxplus_{k\in [n_2]}t_{2,k}\geq 0, \bigboxplus_{k\in [n_3]}t_{3,k}\geq 0
\bigboxplus_{\substack{j\in \{1,2,3\}\\k\in [n_j]}}t_{k,j}=1, n_j\in \mathbb N, j=1,2,3. \Bigg\}\end{align*}

Let us compute the intermediate point between two points. Since $x^{(1)}$ and $x^{(2)}$ are copositive the have only to compute the intermediate points between  $x^{(1)}$ and $x^{(3)}$ and  $x^{(2)}$ and $x^{(3)}$. We have:

$(i)$ $I=\{1\}$, $I=\{2\}$, $I=\{3\}$ and $J=\{1,3\}$.   
$$t_{\{1\},\{1,3\}}=\big(1,1\big),\; t_{\{2\},\{1,3\}}=\big(\frac{1}{3},1\big), \; t_{\{3\},\{1,3\}}=\big(\frac{1}{3},1\big).$$
It follows that:
$$\zeta_{\{1\},\{1,3\}}= x^{(1)}\boxplus x^{(3)}=\begin{pmatrix}
0\\3\\3
\end{pmatrix}\quad \zeta_{\{2\},\{1,3\}}=\frac{1}{3}x^{(1)}\boxplus x^{(3)}=\begin{pmatrix}
-1\\0\\0
\end{pmatrix} \quad \zeta_{\{3\},\{1,3\}}=\frac{1}{3}x^{(1)}\boxplus x^{(3)}=\begin{pmatrix}
-1\\0\\0
\end{pmatrix}.$$
 
$(ii)$ $I=\{1\}$, $I=\{2\}$, $I=\{3\}$ and $J=\{2,3\}$.
$$t_{\{2\},\{2,3\}}=\big(\frac{1}{3},1\big),\, t_{\{2\},\{2,3\}}=\big(\frac{1}{2},1\big), \; t_{\{3\},\{2,3\}}=\big(\frac{1}{3},1\big).$$
$$\zeta_{\{1\},\{2,3\}}=\frac{1}{3}x^{(2)}\boxplus x^{(3)}=\begin{pmatrix}
0\\-1\\0
\end{pmatrix}\quad \zeta_{\{2\},\{2,3\}}=\frac{1}{2}x^{(2)}\boxplus x^{(3)}=\begin{pmatrix}
\frac{3}{2}\\0\\\frac{3}{2}
\end{pmatrix} \quad \Gamma_{\{3\},\{2,3\}}=\frac{1}{3}x^{(2)}\boxplus x^{(3)}=\begin{pmatrix}
0\\-1\\0
\end{pmatrix}.$$

$(iii)$    $I=\{1,2\}$ and $J=\{1,2,3\}$.

We have 
{\scriptsize $$ t_{\{1,2\},1} =\frac{\begin{vmatrix}0&3&-1\\0&2&-1\\1&1&1

\end{vmatrix}_\infty}{\begin{vmatrix}1&3&-1\\1&2&-1\\1&1&1

\end{vmatrix}_\infty}=\frac{(-3)\boxplus 2}{(-3)\boxplus 2\boxplus (-3)\boxplus 2\boxplus (-9)\boxplus 1}=\frac{1}{3};$$

$$ t_{\{1,2\},2} =\frac{\begin{vmatrix}1&0&-1\\3&0&-1\\1&1&1

\end{vmatrix}_\infty}{\begin{vmatrix}1&3&-1\\1&2&-1\\1&1&1

\end{vmatrix}_\infty}=\frac{(-3)\boxplus 1}{(-3)\boxplus 2\boxplus (-3)\boxplus 2\boxplus (-9)\boxplus 1}=\frac{-3}{-9}=\frac{1}{3};$$

$$ t_{\{1,2\},3} =\frac{\begin{vmatrix}1&3&0\\3&2&0\\1&1&1

\end{vmatrix}_\infty}{\begin{vmatrix}1&3&-1\\1&2&-1\\1&1&1

\end{vmatrix}_\infty}=\frac{(-9)\boxplus 2}{(-3)\boxplus 2\boxplus (-3)\boxplus 2\boxplus (-9)\boxplus 1}=\frac{-9}{-9}= {1}.$$}

It follows that:

 $$\zeta_{\{1,2\},\{1,2,3\}}= -\frac{ 1}{9}\Big[
 (-3)x^{(1)}\boxplus  2  x^{(1)}\boxplus   (-3)x^{(2)}\boxplus   x^{(2)}\boxplus (-9)  x^{(3)}\boxplus 2 x^{(3)}\Big]$$
 
 Hence:
 {\scriptsize  $$\zeta_{\{1,2\},\{1,2,3\}}=  -\frac{ 1}{9}\Bigg[ (-3)\begin{pmatrix}
1\\3\\3
\end{pmatrix} \boxplus  2\begin{pmatrix}
1\\3\\3
\end{pmatrix}  \boxplus   (-3) \begin{pmatrix}
3\\2\\3
\end{pmatrix}\boxplus  \begin{pmatrix}
3\\2\\3
\end{pmatrix}\boxplus (-9) \begin{pmatrix}
-1\\-1\\-1
\end{pmatrix}\boxplus  2 \begin{pmatrix}
-1\\-1\\-1
\end{pmatrix}\Bigg]= -\frac{ 1}{9} \begin{pmatrix}
(-3)\boxplus 2 \boxplus (-9)\boxplus 3 \boxplus 9 \boxplus  (-2) \\(-9)\boxplus 6\boxplus (-6)\boxplus 2 \boxplus 9 \boxplus (-2)\\(-9)\boxplus 6 \boxplus (-9)\boxplus 3 \boxplus 9 \boxplus (-2)
\end{pmatrix}=\begin{pmatrix}
0\\0\\1
\end{pmatrix}.$$}

$(iv)$ $I=\{2,3\}$ and $J=\{1,2,3\}$.
{\scriptsize $$ t_{\{2,3\},1} =\frac{\begin{vmatrix}0&2&-1\\0&3&-1\\1&1&1

\end{vmatrix}_\infty}{\begin{vmatrix}3&2&-1\\3&3&-1\\1&1&1

\end{vmatrix}_\infty}=\frac{(-2)\boxplus 3}{(-6)\boxplus 3\boxplus 3\boxplus (-2)\boxplus  9 \boxplus (-3)}=\frac{ 3}{ 9}=\frac{1}{3}$$

$$ t_{\{2,3\},2} =\frac{\begin{vmatrix}3&0&-1\\3&0&-1\\1&1&1

\end{vmatrix}_\infty}{\begin{vmatrix}3&2&-1\\3&3&-1\\1&1&1

\end{vmatrix}_\infty}=\frac{(-3)\boxplus 3}{(-6)\boxplus 3\boxplus 3\boxplus (-2)\boxplus  9 \boxplus (-3)}=\frac{ 0}{ 9}=0;$$

$$ t_{\{2,3\},3} =\frac{\begin{vmatrix}3&2&0\\3&3&0\\1&1&1

\end{vmatrix}_\infty}{\begin{vmatrix}3&2&-1\\3&3&-1\\1&1&1

\end{vmatrix}_\infty}=\frac{(-6)\boxplus 9}{(-6)\boxplus 3\boxplus 3\boxplus (-2)\boxplus  9 \boxplus (-3)}=\frac{ 9}{ 9}=1.$$}

$$\zeta_{\{2,3\},\{1,2,3\}}= \frac{1}{9}\Big[3 x^{(1)}\boxplus   ( -2 ) x^{(1)} \boxplus    (- 6) x^{(3)}\boxplus 9 x^{(3)}\Big].$$

Hence:
 {\scriptsize  $$\zeta_{\{2,3\},\{1,2,3\}}= \frac{1}{9}\Bigg[ 3\begin{pmatrix}
1\\3\\3
\end{pmatrix} \boxplus   (-2) \begin{pmatrix}
1\\3\\3
\end{pmatrix}  \boxplus  (-6) \begin{pmatrix}
-1\\-1\\-1
\end{pmatrix} \boxplus  9 \begin{pmatrix}
-1\\-1\\-1
\end{pmatrix}\Bigg]=\frac{1}{9} \begin{pmatrix}
3\boxplus (-2)\boxplus 6 \boxplus (-9)\\9\boxplus (-6)\boxplus 6 \boxplus (-9)\\9\boxplus (-6)\boxplus 6 \boxplus (-9)
\end{pmatrix}=\begin{pmatrix}
-1\\0\\0
\end{pmatrix}.$$}

$(v)$ $I=\{1,3\}$, $J=\{1,2,3\}$.\\

{\scriptsize $$ t_{\{1,3\},1} =\frac{\begin{vmatrix}0&3&-1\\0&3&-1\\1&1&1

\end{vmatrix}_\infty}{\begin{vmatrix}1&3&-1\\3&3&-1\\1&1&1

\end{vmatrix}_\infty}=\frac{(-3)\boxplus 3}{(-9)\boxplus 1\boxplus 3\boxplus 1\boxplus  3 \boxplus (-3)}=\frac{0}{- 9}=0;$$

$$ t_{\{1,3\},2} =\frac{\begin{vmatrix}1&0&-1\\3&0&-1\\1&1&1

\end{vmatrix}_\infty}{\begin{vmatrix}1&3&-1\\3&3&-1\\1&1&1

\end{vmatrix}_\infty}=\frac{(-3)\boxplus 1}{(-9)\boxplus 1\boxplus 3\boxplus 1\boxplus  3 \boxplus (-3)}=\frac{ -3}{- 9}=\frac{1}{3};$$

$$ t_{\{1,3\},3} =\frac{\begin{vmatrix}1&3&0\\3&3&0\\1&1&1

\end{vmatrix}_\infty}{\begin{vmatrix}1&3&-1\\3&3&-1\\1&1&1

\end{vmatrix}_\infty}=\frac{3\boxplus (-9)}{(-9)\boxplus 1\boxplus 3\boxplus 1\boxplus  3 \boxplus (-3)}=\frac{ -9}{- 9}=1.$$}

 $$\zeta_{\{1,3\},\{1,2,3\}}=    -\frac{1}{9}\Big[(-3)x^{(2)}\boxplus   x^{(2)}\boxplus  3x^{(3)}\boxplus  (-9) x^{(3)}\Big]$$
 
Hence:
 {\scriptsize  $$\zeta_{\{1,3\},\{1,2,3\}}=   -\frac{1}{9} \Bigg[(-3) \begin{pmatrix}
3\\2\\3
\end{pmatrix}\boxplus    \begin{pmatrix}
3\\2\\3
\end{pmatrix}\boxplus 3 \begin{pmatrix}
-1\\-1\\-1
\end{pmatrix}\boxplus (-9) \begin{pmatrix}
-1\\-1\\-1
\end{pmatrix}\Bigg]= -\frac{1}{9} \begin{pmatrix}
(- 9)\boxplus 3\boxplus (-3)\boxplus 9\\(-6)\boxplus 2\boxplus (-3)\boxplus 9\\(-9)\boxplus 3\boxplus (-3)\boxplus 9\end{pmatrix}  
=\begin{pmatrix}
0\\-1\\0
\end{pmatrix}.$$}
For all $\epsilon=(\epsilon_1, \epsilon_2, \epsilon_3)\in \{-1,1\}^3$ let us define  
$$K_\epsilon=\{x\in \Real^3: \epsilon_ix_i\geq 0, i=1,2,3\}.$$  We have $x^{(1)}=(1,3,3),x^{(2)}=(3,2,3)\in K_{1,1,1} $. Moreover $\zeta_{\{1\},\{1,3\}}=(0,3,3),  \zeta_{\{2\},\{2,3\}}=(\frac{3}{2}, 0,\frac{3}{2} ), \zeta_{\{1,2\},\{1,2,3\}}=(0,0,1)\in K_{1,1,1}$.   It follows that 
$$Co^\infty(A)\cap K_{1,1,1}=Co^\infty\Big(\big\{x^{(1)},x^{(2)} ,\zeta_{\{1\},\{1,3\}} ,  \zeta_{\{2\},\{2,3\}} , \zeta_{\{1,2\},\{1,2,3\}}\big\}\Big).$$
Moreover: $\zeta_{\{2\},\{1,3\}}=(-1,0,0),  \zeta_{\{1\},\{2,3\}}=(0, -1,0 ), \zeta_{\{1,2\},\{1,2,3\}}=(0,0,1)\in K_{-1,-1,1}$, hence:
$$Co^\infty(A)\cap K_{-1,-1,1}=Co^\infty\Big(\big\{ \zeta_{\{2\},\{1,3\}} ,  \zeta_{\{1\},\{2,3\}} , \zeta_{\{1,2\},\{1,2,3\}}\big\}\Big).$$
Now note that $\zeta_{\{2\},\{1,3\}}=(-1,0,0),  \zeta_{\{1\},\{1,3\}}=(0,3,3),  \zeta_{\{1,2\},\{1,2,3\}}=(0,0,1)\in K_{-1, 1,1}$. Thus $$Co^\infty(A)\cap K_{-1, 1,1}=Co^\infty\Big(\big\{\zeta_{\{2\},\{1,3\}} , \zeta_{\{1\},\{1,3\}}, \zeta_{\{1,2\},\{1,2,3\}}\big\}\Big).$$
 Moreover $\zeta_{\{1\},\{2,3\}}=(0,-1,0),  \zeta_{\{2\},\{2,3\}}=(\frac{3}{2},0,\frac{3}{2}),   \zeta_{\{1,2\},\{1,2,3\}}=(0,0,1)\in K_{ 1,- 1,1}$. Therefore:
$$  Co^\infty(A)\cap K_{ 1,- 1,1}=Co^\infty\Big(\big\{\zeta_{\{1\},\{2,3\}} , \zeta_{\{2\},\{2,3\}},\zeta_{\{1,2\},\{1,2,3\}}\big\}\Big).$$
Finally, $x^{(3)}=(-1,-1,-1), \zeta_{\{2\},\{1,3\}}=(-1,0,0), \zeta_{\{1\},\{2,3\}}=(0, -1,0 )\in K_{-1,-1,-1} $. Therefore:
$$Co^\infty(A)\cap K_{-1,- 1,-1}=Co^\infty \Big (\big\{x^{(3)}, \zeta_{\{2\},\{1,3\}}, \zeta_{\{1\},\{2,3\}}\big\}\Big ).$$
It follows that:
{\scriptsize \begin{align*}Co^\infty(A)&=Co^\infty\Bigg(\Bigg\{\begin{pmatrix}
1\\3\\3
\end{pmatrix} ,\begin{pmatrix}
3\\2\\3
\end{pmatrix} , \begin{pmatrix}
0\\3\\3
\end{pmatrix}  , \begin{pmatrix}
\frac{3}{2}\\0\\\frac{3}{2}
\end{pmatrix}  ,\begin{pmatrix}
0\\0\\1
\end{pmatrix} \Bigg\}\Bigg) \cup  Co^\infty\Bigg( \Bigg\{\begin{pmatrix}
-1\\0\\0
\end{pmatrix} , \begin{pmatrix}
0\\-1\\0
\end{pmatrix} ,\begin{pmatrix}
0\\0\\1
\end{pmatrix}\Bigg\}\Bigg)\\&\cup Co^\infty\Bigg(\Bigg\{\begin{pmatrix}
-1\\0\\0
\end{pmatrix} , \begin{pmatrix}
0\\3\\3
\end{pmatrix}, \begin{pmatrix}
0\\0\\1
\end{pmatrix}\Bigg\}\Bigg) \cup Co^\infty\Bigg(\Bigg\{\begin{pmatrix}
0\\-1\\0
\end{pmatrix} , \begin{pmatrix}
\frac{3}{2}\\0\\\frac{3}{2}
\end{pmatrix}, \begin{pmatrix}
0\\0\\1
\end{pmatrix}\Bigg\}\Bigg)\cup Co^\infty \Bigg (\Bigg\{\begin{pmatrix}
-1\\-1\\-1
\end{pmatrix}, \begin{pmatrix}
-1\\0\\0
\end{pmatrix}, \begin{pmatrix}
0\\-1\\0
\end{pmatrix}\Bigg\}\Bigg).\end{align*}}
It is easy to check that $(0,0,0)\notin Co^\infty(A)$ (see also Example \ref{hypcontains} below). However $x^{(1)}\boxplus x^{(2)}=(1,3,3)\vee (3,2,3)=(3,3,3)\in Co^\infty(A)$.   Since 
$\frac{1}{3}\big(x^{(1)}\boxplus x^{(2)})\boxplus x^{(3)}=(1,1,1)\boxplus(-1,-1,-1)=(0,0,0)$ we deduce that $Co^\infty(A)$ is not idempotent symmetric convex since it does not contain $0$. Hence $ \mathrm{IS}[A]\nsubseteq Co^\infty(A). $ }

\end{expl}

\begin{rem}\label{NonConv}The proposed example illustrates a typical situation where the convergence of a sequence of pairs $\{(x^{(p)}, y^{(p)})\}_{p \in \mathbb{N}}$ to $(x, y)$ does not imply that $Co^p(x^{(p)}, y^{(p)})$ converges to $Co^\infty(x, y)$ when the sequence $\{(x^{(p)}, y^{(p)})\}_{p \in \mathbb{N}}$ is not contained within the same $n$-dimensional orthant.
For example, the sequence $\{z^{(p)}\}_{p \in \mathbb{N}}$ with  
$$z^{(p)} = \frac{1}{2^{\frac{1}{2p+1}}}x^{(1)} \stackrel{p}{+} \frac{1}{2^{\frac{1}{2p+1}}}x^{(2)}= \frac{1}{2^{\frac{1}{2p+1}}}\begin{pmatrix}
1\\3\\3
\end{pmatrix} \stackrel{p}{+} \frac{1}{2^{\frac{1}{2p+1}}}\begin{pmatrix}
3\\2\\3
\end{pmatrix}$$  
converges to $x^{(1)} \boxplus x^{(2)} = (3, 3, 3)$. However, there is no sequence, each term of which belongs to $Co^p(\{x^{(1)}, x^{(2)}, x^{(3)}\})$, that converges to $0$.
 Since  
$$0 = \frac{1}{3}(x^{(1)} \boxplus x^{(2)}) \boxplus x^{(3)},$$  
this shows that  
$$\Lim_{p \to \infty} Co^p(z^{(p)}, x^{(3)}) \neq Co^\infty(x^{(1)} \boxplus x^{(2)}, x^{(3)}).$$

\end{rem}

\begin{rem}We observed in the example above that, in general, a $\mathbb{B}$-polytope is not idempotent symmetric convex, although this holds true when the set $A$ contains only two points. 
We have not been able to to find in dimension 3 a counter-example of the fact that $\mathrm{IS}[A]\supset Co^\infty(A)$.  However, in the general case, there is no guarantee that this inclusion holds, particularly due to the pathological cases highlighted in Remark \ref{NonConv}. On the contrary, if this inclusion were valid, it would imply that idempotent symmetric convexity and $\mathbb{B}$-convexity are equivalent.

However, it is important to emphasize that the proposed formulation, which involves any arbitrary finite subset of the considered set, provides numerous advantages. Firstly, it is also defined in algebraic terms and is very useful for solving the separation problem, as shown in the next section. Moreover, it explicitly provides a polytope structure that, while not coinciding with the convex hull, is nonetheless the limit of a sequence of $\varphi_p$-convex hulls. 

 \end{rem}

In the following, we show that when the dimension of the Euclidean vector space is 2, the $ \mathbb{B} $-polytopes are idempotent  symmetric  convex. It follows that when $ n = 2 $, $ \mathbb{B} $-convexity and idempotent symmetric convexity are equivalent notions. This property also holds on each $n$-dimensional orthant, as shown in \cite{bh}.

\begin{prop}Suppose that $n=2$. Then, any subset $C$ of $\Real^2$ is $\mathbb B$-convex if and only if it idempotent symmetric convex.  It follows that that for any finite subset $A$ of $\Real^2$, $Co^\infty(A)=\mathrm{IS}[A]=\mathbb B[A]$.  

\end{prop}
{\bf Proof:} To prove this result, we need to show that  for all finite subsets     $A=\{x^{(1)},...,x^{(m)}\}\subset \mathbb R^2$ we have: 
$$ Co^\infty(A) \subset \mathrm{IS}[A]. $$

We first establish that all the intermediate points of \( Co^\infty(A) \) belong to \( \mathrm{IS}[A] \). By definition, all the intermediate points of \( Co^\infty(A) \) that have only one null component have the form  
$$ \zeta_{\{i\},\{j,k\}} = t_{\{i\},\{j,k\},j} x^{(j)} \boxplus t_{\{i\},\{j,k\},k} x^{(k)} $$   
with $ [\zeta_{\{i\},\{j,k\}}]_i = 0 $,  
$  \max \{t_{\{i\},\{j,k\},j}, t_{\{i\},\{j,k\},k}\} = 1, $   
and  $ t_{\{i\},\{j,k\},j}, t_{\{i\},\{j,k\},k} \geq 0, $  
for some $ i \in \{1,2\} \) and \( j,k \in [m] $ with $ j \neq k $. It follows that all these points belong to $ \mathrm{IS}[A] $. The only intermediate point having two components equal to zero is $ 0 $.

If $0 \notin Co^\infty(A) $, the problem is solved. Suppose that \( 0 \in Co^\infty(A) \) and that $ A \neq \{0\}$. Pick some $ v \in Co^\infty(A) $ and define  
$  \lambda_+ = \max \{\lambda : \lambda v \in Co^\infty(A)\} $   
and  $  \lambda_- = \min \{\lambda : \lambda v \in Co^\infty(A)\}. $ Since \( Co^\infty(A) \) is compact, the points  
$  v_- = \lambda_- v \quad \text{and} \quad v_+ = \lambda_+ v $   
lie in $ \partial Co^\infty(A) $. This implies that there exist some intermediary points 
$  z_-, w_-, z_+, w_+ \in \partial Co^\infty(A) $   
such that  
 $ v_- \in \mathrm{IS}[z_-, w_-] \quad \text{and} \quad v_+ \in \mathrm{IS}[z_+, w_+]. $ However, $ 0 \in \mathrm{IS}[v_-, v_+] $. Consequently, since $v_-, v_+\in\mathrm{IS}[A] $ we deduce that $ 0 \in \mathrm{IS}[A] $. It follows that $Co^\infty(A)\subset \mathrm{IS}[A]$. Therefore, idempotent symmetric convexity implies $\mathbb{B}$-convexity. Since $\mathbb{B}$-convexity implies  idempotent symmetric convexity, these convexities are   equivalent, which concludes the proof. We have shown that all the intermediate points of $Co^\infty(A)$ are in $\mathrm{IS}[A]$. It follows that for any 2-dimensional orthant $K$, $\mathrm{IS}[A] \cap K \subset Co^\infty(A) \cap K$. Therefore, $\mathrm{IS}[A] = \mathbb{B}[A] \subset Co^\infty(A)$, which proves that $\mathbb{B}[A] = Co^\infty(A)$.
 $\Box$.\\

A consequence of this property is that constructing a limit polytope from the convex hull of any pair of points yields a $\mathbb  B$-polytope in a two-dimensional space, as illustrated in Figure \ref{Polalg}.1.

 \section{From the  {Separation} of Polytopes to the Separation of Convex Sets }\label{DualStruc}

In this section, we establish a general separation theorem for $\mathbb{B}$-convex sets in $\mathbb{R}^n$. A separation result was previously derived in \cite{b17} for fixed $\mathbb{B}$-convex sets, defined as those subsets $S \subset \mathbb{R}^n$ that satisfy $Co^\infty(S) = S$. Equivalently, a set $S$ is fixed $\mathbb{B}$-convex if for every (not necessarily finite) subset $E \subset S$, we have $Co^\infty(E) \subset S$. This definition is notably restrictive and, in particular, implies that all fixed $\mathbb{B}$-convex sets are closed. The results presented here extend the scope of separation theorems to the broader, topologically-defined notion of $\mathbb{B}$-convexity introduced in \cite{bh}, which is now characterized algebraically.
 
\subsection{Limit of Idempotent Lower (Upper) Symmetric Halfspaces}\label{Approx}

Let $\langle \cdot, \cdot \rangle_p$ denote the $\varphi_p$-inner product defined for any $p\in \mathbb N$ and all $x,y\in \Real^n$ as $\langle x,y\rangle _p=\varphi_p^{-1}\big(\langle \phi_p(x), \phi_p(y)\rangle\big)$ (see \cite{b15, b17}). In the following we introduce the operation $ \langle \cdot,\cdot\rangle_\infty :\Real^n\times \Real^n\longrightarrow
 \Real$ defined for all $x,y\in \Real^n$ by $\langle x,y\rangle_\infty =\lim_{p\longrightarrow \infty }\langle x, y\rangle_p=\bigboxplus_{i\in
 [n]}x_iy_i$. Let $\|\cdot\|_\infty$ be the Tchebychev
 norm defined by $\|x\|_\infty=\max_{i\in [n]}|x_i|$. It {was} established in \cite{b15} that
for all $x,y\in \Real^n$, we have:
 $(i)$  $\sqrt{\langle x,x\rangle_\infty} =\|x\|_\infty$;
  $(ii)$  $|\langle x,y\rangle_\infty| \leq \|x\|_\infty \|y\|_\infty$;
   $(iii)$ For all $\alpha\in \Real$, $\alpha \langle x,y\rangle_\infty= \langle \alpha x,y\rangle_\infty=\langle  x,\alpha y\rangle_\infty$.

In the following, we say that a map $f:\Real^n\longrightarrow
\Real$ is an {\bf  idempotent symmetric form } if there exists some
$a\in\Real^n$ such that:
\begin{equation}f(x)=\bigboxplus_{i\in [n]}a_ix_i=\langle a,x\rangle_\infty.
\end{equation}

We intend to analyze the dual relation between $\mathbb B$-convex
sets and   the maps of the form $x\mapsto \langle
a,x\rangle_\infty$, where $a\in \Real^n$.

 The function above is
depicted in Figure \ref{Rec}.2. All the points such that $\langle
a,x\rangle_{\infty}= 0$ are represented by the doted line. In the
following, for all subsets $A$ of $\Real^n$ $\cl  (A)$ and $\inter
(A)$ respectively stand for the closure and the interior of $A$.

\begin{center}{\scriptsize % This is a LaTeX picture output by TeXCAD.

\unitlength 1mm % = 2.845pt
\linethickness{0.4pt}
\ifx\plotpoint\undefined\newsavebox{\plotpoint}\fi % GNUPLOT compatibility
\begin{picture}(152.625,60.375)(0,0)
{\linethickness{0.15mm} 
}
%\emline(69.571,52.005)(52.863,41.295)
\multiput(69.571,52.005)(-.0525396226,-.0336792453){318}{\line(-1,0){.0525396226}}
%\end
\put(52.863,41.295){\line(-1,0){33.415}}
\put(19.448,41.295){\line(0,-1){21.42}}
%\emline(19.448,19.875)(8.31,12.735)
\multiput(19.448,19.875)(-.052539623,-.033679245){212}{\line(-1,0){.052539623}}
%\end
\put(55.285,34.907){$\cl{\big[ \langle a,.\rangle_{\infty}\leq c\big]}$}
\put(52.863,41.295){\line(-1,0){33.415}}
\put(19.448,41.295){\line(0,-1){21.42}}
%\emline(19.448,19.875)(8.31,12.735)
\multiput(19.448,19.875)(-.052539623,-.033679245){212}{\line(-1,0){.052539623}}
%\end
\put(32.977,60.287){$x_2$}
\put(70.25,27.875){$x_1$}
\put(5,0){ {\bf Figure  \ref{DualStruc}.1:}The domain $\cl{\big[\langle a,.\rangle_{\infty}\leq c\big]}$ as a limit of $\varphi_p$-halfspaces.}
{
}
\put(150.625,58.375){\line(-1,-1){15}}
\put(140.625,48.675){\line(-1,0){40}}
\put(100.625,48.675){\line(0,-1){5}}
\put(100.625,43.675){\line(4,-3){10}}
\put(110.625,36.375){\line(0,-1){15.5}}
\put(135.625,43.675){\line(0,-1){22.5}}
\put(135.625,21.075){\line(-1,0){25.2}}
{ 
\put(106.525,18.875){\line(1,0){28.2}}
\put(134.725,18.875){\line(3,-2){16}}
}
\put(118.125,25.375){$0$}
\put(120.625,28.375){\circle*{1}}
\put(150.685,58.585){\circle*{1}}
\put(100.75,43.5){\circle*{1}}
\put(110.75,20.875){\circle*{1}}
\put(75.75,49.5){\circle*{.951}}
\put(81.625,10){\circle*{1}}
\put(135.75,31.125){\circle*{1}}
\put(100.125,32.463){\circle*{.951}}
\put(120.125,60.375){$x_2$}
\put(152.625,27.875){$x_1$}
\put(96.5,0){{\bf Figure  \ref{DualStruc}.2:} Separation in Limit of two Polytopes}
\put(133.688,5.619){$[f ^-\leq c ]$}
\put(150.625,10.375){$[f^+\geq c]$}
\put(81.875,28.25){\vector(1,0){68.375}}
{
}
{
}
{ 
}
\put(31.477,25.287){$0$}
\put(33.977,29.524){\circle*{1}}
\put(0,29.399){\vector(1,0){68.375}}
\qbezier(36.853,41.474)(53.42,41.875)(68.875,51.738)
\qbezier(106.65,27.649)(105.975,37.645)(89.4,46.969)
\qbezier(36.713,41.474)(19.727,41.607)(19.448,30.496)
\qbezier(36.296,41.563)(47.294,42.723)(67.761,52.095)
\qbezier(36.017,41.474)(23.973,40.849)(22.79,40.581)
\qbezier(19.448,30.675)(20.423,38.394)(20.841,39.153)
\qbezier(22.79,40.492)(21.188,40.135)(20.98,39.064)
\qbezier(60.687,52.725)(51.986,45.807)(33.538,41.568)
\qbezier(120.5,18.75)(134.5,18.75)(150,7.75)
\qbezier(100.625,43.875)(101.438,50.063)(120.5,50)
\qbezier(120.5,50)(143.688,50.188)(150.625,58.625)
\qbezier(110.25,21)(135.75,22.188)(135.75,30.125)
\qbezier(135.75,30.125)(136,43.625)(150.75,58.125)
\qbezier(100.75,43.5)(110.125,36.313)(110,28.375)
\qbezier(110.125,28.25)(110.188,19.625)(110.5,21.5)
\qbezier(75.675,49.678)(81.525,48.765)(81.525,9.523)
%\emline(75.35,49.815)(99.725,38.993)
\multiput(75.35,49.815)(.0759345794,-.0337133956){321}{\line(1,0){.0759345794}}
%\end
%\emline(81.525,9.45)(81.688,47.036)
\multiput(81.525,9.45)(.0326,7.5172){5}{\line(0,1){7.5172}}
%\end
\put(99.875,38.7){\line(0,-1){15.638}}
%\emline(99.875,23.063)(81.25,9.563)
\multiput(99.875,23.063)(-.046446384,-.0336658354){401}{\line(-1,0){.046446384}}
%\end
\qbezier(75.75,49.613)(99.313,38.644)(99.625,33.3)
\qbezier(100,33.188)(99.688,22.163)(81.125,9.788)
\put(120.625,6.375){\vector(0,1){52.5}}
\put(34.102,6.037){\vector(0,1){52.5}}
\qbezier(106.75,27.875)(107.188,23.125)(107.875,21.375)
\qbezier(120.5,18.75)(110.813,19)(109.375,19.75)
\qbezier(109.375,19.75)(108.375,20.313)(107.875,21.125)
%\emline(106.307,18.764)(106.396,36.353)
\put(106.307,18.764){\line(0,1){17.589}}
%\end
%\emline(106.307,36.441)(90.397,47.49)
\multiput(106.307,36.441)(-.0485060976,.0336859756){328}{\line(-1,0){.0485060976}}
%\end
\put(63.863,56.448){\makebox(0,0)[cc]{$\big[ \langle a^{(p)},.\rangle_{p}\leq c_p\big]$}}
\qbezier(33.569,41.463)(21.954,39.395)(19.462,29.054)
\qbezier(3.234,13.038)(16.439,19.349)(19.25,29.266)
\qbezier(4.082,11.871)(18.295,19.295)(19.356,29.266)
\qbezier(4.825,10.598)(18.985,19.508)(19.356,28.841)
%\dashline{1}(19.533,19.897)(52.855,41.198)
\multiput(19.462,19.826)(.0507964,.0324719){16}{\line(1,0){.0507964}}
\multiput(21.088,20.865)(.0507964,.0324719){16}{\line(1,0){.0507964}}
\multiput(22.713,21.904)(.0507964,.0324719){16}{\line(1,0){.0507964}}
\multiput(24.339,22.944)(.0507964,.0324719){16}{\line(1,0){.0507964}}
\multiput(25.964,23.983)(.0507964,.0324719){16}{\line(1,0){.0507964}}
\multiput(27.59,25.022)(.0507964,.0324719){16}{\line(1,0){.0507964}}
\multiput(29.215,26.061)(.0507964,.0324719){16}{\line(1,0){.0507964}}
\multiput(30.841,27.1)(.0507964,.0324719){16}{\line(1,0){.0507964}}
\multiput(32.466,28.139)(.0507964,.0324719){16}{\line(1,0){.0507964}}
\multiput(34.092,29.178)(.0507964,.0324719){16}{\line(1,0){.0507964}}
\multiput(35.717,30.217)(.0507964,.0324719){16}{\line(1,0){.0507964}}
\multiput(37.343,31.256)(.0507964,.0324719){16}{\line(1,0){.0507964}}
\multiput(38.968,32.295)(.0507964,.0324719){16}{\line(1,0){.0507964}}
\multiput(40.594,33.335)(.0507964,.0324719){16}{\line(1,0){.0507964}}
\multiput(42.219,34.374)(.0507964,.0324719){16}{\line(1,0){.0507964}}
\multiput(43.845,35.413)(.0507964,.0324719){16}{\line(1,0){.0507964}}
\multiput(45.47,36.452)(.0507964,.0324719){16}{\line(1,0){.0507964}}
\multiput(47.095,37.491)(.0507964,.0324719){16}{\line(1,0){.0507964}}
\multiput(48.721,38.53)(.0507964,.0324719){16}{\line(1,0){.0507964}}
\multiput(50.346,39.569)(.0507964,.0324719){16}{\line(1,0){.0507964}}
\multiput(51.972,40.608)(.0507964,.0324719){16}{\line(1,0){.0507964}}
%\end
\end{picture}
}
\end{center}

\bigskip

In the following, for all maps $f:\Real^n\longrightarrow \Real$ and
all real numbers $c$, the notation $[\,f\leq c\,]$ stands for the
set $f^{-1}(\,]-\infty, c]\,)$.  Similarly, $[\,f< c\,]$ stands
for  $f^{-1}(\,]-\infty, c[\,)$ and $[\,f\geq c\,]=[\,-f\leq
-c\,]$.

  For all $u,v\in \Real$, let us define the binary operation
\begin{equation*}\label{Bform}u\stackrel{-}{\smile} v=
\left\{\begin{matrix}
u &\hbox{ if } &|u|& > &|v|\\
\min \{u,v\}&\hbox{ if }&|u|&=&|v|\\
v& \hbox{ if }& |u|&<&|v|.\end{matrix}\right.\end{equation*} An
elementary calculus shows that $u\boxplus v=\frac{1}{2}\Big[u
\stackrel{-}{\smile} v-\big[(-u) \stackrel{-}{\smile}(-v)\big]
\Big]$. Note that in \cite{ir10}, a polynomial theory for supertropical algebra was proposed. To achieve this, a similar operation was constructed by extending the tropical semiring with a suitable valuation function. However, in our approach, this binary operation is derived from the idempotent semiring $(\Real_+, \max, \cdot)$ and the absolute value function.

 Similarly one can introduce a symmetric binary operation defined for all $u,v\in \Real$ as
$u\stackrel{+}{\smile}v=-\big[(-u)\stackrel{-}{\smile}(-v)\big]$. This implies that we have $u\boxplus v=\frac{1}{2}\Big[(u\stackrel{-}{\smile} v)+ (u\stackrel{+}{\smile} v)\Big]$.

 Given $m$ elements $u_1, ..., u_m$ of $\Real$, not all of which are $0$, let $I_+$, respectively $I_-$,
 be the set of indices for which $0 < u_i$, respectively $u_i < 0$. We can then write
 $u_1\stackrel{-}{\smile}\cdots\stackrel{-}{\smile} u_m = (\stackrel{-}{\smile}_{i\in I_+}u_i) \stackrel{-}{\smile} (\stackrel{-}{\smile}_{i\in I_-}u_i)
 = (\max_{i\in I_+}u_i)\stackrel{-}{\smile}(\min_{i\in I_-}u_i)$.
 
We define a {\bf lower   idempotent symmetric form }on ${\mathbb R}^n_+$
as a map $g: {\mathbb R}^n\to\Real$ such that for all $(x_1, ...,
x_n)\in{\mathbb R}^n_+$,
 \begin{equation}\label{eqcarBform}
 g(x_1, ..., x_n) = {\langle a,x\rangle}_\infty^{-} =a_1x_1 \stackrel{-}{\smile}\cdots \stackrel{-}{\smile} a_nx_n.
 \end{equation}
Note that a similar closely related It was established in \cite{b15} that for all $c \in \Real$,
$g^{-1}\left(\,]-\infty,c]\right)=\left\{x\in \Real^n: g(x)\leq
c\right\}$ is closed. It follows that  a lower idempotent symmetric form is
lower semi-continuous. It was established in \cite{bh3} that
$g^{-1}\left(\,]-\infty,c]\right)\cap \Real_+^n$ is a $\mathbb
B$-halfspace, that is a $\mathbb B$-convex subset of $\Real_+^n$
whose the complement in $\Real_+^n$ is also $\mathbb B$-convex.

Similarly, one can define an {\bf upper idempotent symmetric form} as  a
map $h: {\mathbb R}^n\to\Real$  such that, for all $(x_1, ...,
x_n)\in{\mathbb R}^n$,
 \begin{equation}\label{eqcarBform}
 h(x_1, ..., x_n) = {\langle a,x\rangle}_\infty^{+} =a_1x_1 \stackrel{+}{\smile}\cdots \stackrel{+}{\smile} a_n x_n.
 \end{equation}
 For all $x\in \Real^n$, we clearly, have the following identities

 \begin{equation}
 \langle a,x\rangle_\infty^{+}=-\langle a,-x\rangle_\infty^{-}
 \text{ and }
 \langle a,x\rangle_\infty^{-}=-\langle a,-x\rangle_\infty^{+}.
 \end{equation}

The largest (smallest) lower (upper) semi-continuous minorant
(majorant) of a map $f$ is said to be the lower (upper)
semi-continuous regularization of $f$. In the next statements it is shown that the lower (upper) idempotent symmetric forms are the lower (upper) semi-continuous regularized of the idempotent symmetric forms.

It was shown in \cite{b15} that if  $g$ is a lower idempotent symmetric form defined by $ g(x_1,
..., x_n) = a_1x_1\stackrel{-}{\smile}\cdots\stackrel{-}{\smile}
a_nx_n, $ for some $a\in \Real^n$, then $g$ is the lower
semi-continuous regularization of the map $x\mapsto \langle
a,x\rangle_{\infty}=\bigboxplus_{i\in [n]}a_ix_i$. Similarly, if $h$ is an upper idempotent symmetric form defined by $ h(x_1,
..., x_n) = a_1x_1\stackrel{+}{\smile}\cdots\stackrel{+}{\smile}
a_nx_n, $ for some $a\in \Real^n$, then $h$ is the upper
semi-continuous regularization of the map $x\mapsto \langle
a,x\rangle_{\infty}=\bigboxplus_{i\in [n]}a_ix_i$.

 In \cite{b17}, the geometric deformation of a sequence of $\varphi_p$-halfspaces, expressed as $H_p = [f_p \leq c_p]$, was analyzed. For each natural number $p$, $f_p: \mathbb{R}^n \longrightarrow \mathbb{R}$ is a $\varphi_p$-linear form defined by $f(x) = \langle a^{(p)}, x \rangle_p$, with $(a^{(p)}, c_p) \in \mathbb{R}^n \times \mathbb{R}$. It was shown that the {Painlevé-Kuratowski} limit of this sequence can be expressed with respect to the idempotent symmetric forms.

\begin{prop}\label{FundRes2} Let $f$ be an
idempotent symmetric form defined by $f(x)=\langle a,x \rangle_{\infty}$
for some $a\in \mathbb R^n\backslash\{0\}$. For all natural numbers
$p$ let $f_{p}:\Real^n\longrightarrow \Real$ be a map defined by
$f_{p}(x)=\langle a^{(p)},x \rangle_{p}$ where $\{a^{(p)}\}_{p\in
\mathbb{N}}$ is a sequence  of $\Real^n\backslash \{0\}$. If there
exists a sequence $\{c_{p}\}_{p\in \mathbb{N}}\subset \Real$  such
that $\lim_{p\longrightarrow \infty}(a^{(p) },c_{p })=(a,c)$, then:
$$\Lim_{p\longrightarrow \infty}\big [f_p\leq c_p\big ]
=\cl{\big [f\leq c\big ]}=\big [f^-\leq c\big ]\quad \text{ and }\quad \Lim_{p\longrightarrow \infty}\big [f_p\geq c_p\big ]
=\cl{\big [f\geq c\big ]}=\big [f^+\geq c\big ].$$

\end{prop}

The subsequent result is closely related to Proposition \ref{FundRes2}. 
Although the latter could, in principle, be employed to establish it via Hausdorff{-Pompeiu} type arguments, 
we present here a fully self-contained proof. 
Such an approach will prove instrumental in extending our analysis 
from the separation of polytopes to the more general framework of compact $\mathbb{B}$-convex sets.

 \begin{prop}Let $f$ be an
idempotent symmetric form defined by $f(x)=\langle a,x \rangle_{\infty}$
for some $a\in \mathbb R^n\backslash\{0\}$. Let $f^-$ and $f^+$ be respectively its lower and upper regularized forms.  For all natural numbers
$k$ let $f_{k}:\Real^n\longrightarrow \Real$ be a map defined by
$f_{k}(x)=\langle a^{(k)},x \rangle_{\infty}$ where $\{a^{(k)}\}_{k\in
\mathbb{N}}$ is a sequence  of $\Real^n\backslash \{0\}$. If there
exists a sequence $\{c_{k}\}_{k\in \mathbb{N}}\subset \Real$  such
that $\lim_{k\longrightarrow \infty}(a^{(k) },c_{k })=(a,c)$, then:
$$\Lim_{k\longrightarrow \infty}\big [f_k^-\leq c_k\big ]
=\big [f^-\leq c\big ]\quad \text{ and }\quad 
 \Lim_{k\longrightarrow \infty}\big [f_k^+\geq c_k\big ]
=\big [f^+\geq c\big ].$$

 \end{prop}
 \noindent{\bf Proof:} Let us prove the first part of the statement for lower forms. By hypothesis we have:
 $$f^-(x)=\langle a,x \rangle_{\infty}^-=\bigsmileminus_{i\in [n]}a_{i}x_i\quad \text{and}\quad f_k^-(x)=\langle a^{(k)},x \rangle_{\infty}^-=\bigsmileminus_{i\in [n]}a_{i}^{(k)}x_i.$$
 We first establish that  
 $\Ls_{k\longrightarrow \infty}\big[\langle a^{(k )}, \cdot \rangle_{\infty }^-\leq c_k\big ]\subset \big [\langle a, \cdot \rangle _\infty^{-}\leq c\big ].$ Suppose that $z\in  \Ls_{k\longrightarrow \infty}\big[\langle a^{(k )}, \cdot \rangle_{\infty }^-\leq c_k\big ]$. This implies that there is an increasing sequence of natural numbers $\{ k_m \}_{m\in \mathbb N}$ a sequence $\{z^{(m)}\}_{m\in \mathbb N}$ with $z^{(m)}\in \big[\langle a^{(k_m )}, \cdot \rangle_{\infty }^-\leq c_{k_m}\big ]$ for any $m$, and $\lim_{m\longrightarrow \infty} z^{(m)}=z.$ Let $I_+=\{i\in [n]: a_iz_i> 0\}$, $I_0=\{i\in [n]: a_iz_i=0\}$ and $I_-=\{i\in [n]: a_iz_i<0\}$. First suppose that $c\geq 0$. Since $z^{(m)}\in \big[\langle a^{(k_m )}, \cdot \rangle_{\infty }^-\leq c_{k_m}\big ]$ for any $m$, we have
 $$\bigsmileminus_{i\in I_0\cup I_+}a_{i}^{k_m}z_i^{(m)}\leq \bigsmileminus_{i\in I_-}\big(-a_{i}^{(k_m)}z_i^{(m)}\big )\stackrel{-}{\smile} c.$$
 However,  since $\lim_{m\longrightarrow \infty }a_{i}^{(k_m)}z_i^{(m)} =a_iz_i$ for any $i$, we deduce that $\lim_{m\longrightarrow \infty}a_{i}^{(k_m)}z_i^{(m)}=0$ for all $i\in I_0$. Therefore
 $$\lim_{m\longrightarrow \infty}\Big[\bigsmileminus_{i\in I_0\cup  I_+}a_{i}^{(k_m)}z_i^{(m)}\Big]=\max_{i\in I_0\cup I_+}a_{i}z_i=\bigsmileminus_{i\in I_0\cup I_+}a_{i}z_i=\bigsmileminus_{i\in   I_+}a_{i}z_i\geq 0.$$ Using the fact that $-a_iz_i>0$ for all $i\in I_-$
 $$ \lim_{m\longrightarrow \infty}\Big[\bigsmileminus_{i\in I_-}(-a_{i}^{(k_m)}z_i^{(m)})\smile c_{k_m}\Big]=\max_{i\in I_-}\big\{-a_i z_i, c\big \}=\bigsmileminus_{i\in I_-}(-a_{i}z_i)\stackrel{-}{\smile} c\geq 0.$$
 Taking the limit on both sides, we obtain 
 $$\bigsmileminus_{i\in I_0\cup I_+}a_{i}z_i\leq \bigsmileminus_{i\in I_-}(-a_{i}z_i)\stackrel{-}{\smile} c\iff \langle a, \cdot \rangle _\infty^{-}\leq c.$$
 Therefore $z\in  \big [\langle a, \cdot \rangle _\infty^{-}\leq c\big ]$. The proof is similar if $c<0$ because we obtain:
  $$\bigsmileminus_{i\in I_0\cup I_+}a_{i}z_i\smile (- c)=\max_{i\in I_0\cup I_+}\{a_{i}z_i,  - c\}\leq \max_{ i\in I_-}\{ -a_{i} z_i\}= \bigsmileminus_{i\in I_-}a_{i}z_i. $$
  Let us now prove that $\big [\langle a, \cdot \rangle _\infty^{-}\leq c\big ]\subset \Linf_{k\longrightarrow \infty}\big[\langle a^{(k )}, \cdot \rangle_{\infty }^-\leq c_k\big ]$. We consider two cases:
  
  $(i)$ $c\not=0$.  Suppose that $z\in \big [\langle a, \cdot \rangle _\infty^{-}\leq c\big ]$. For all $i\in I_-\cup I_+$, $a_iz_i\not=0$, therefore there is some $k_0$ such that for all $k>k_0$, $a_i^{(k)}\not=0$.  Let us construct the sequence $\{z^{(k)}\}_{\in \mathbb N}$ defined as
  $$z^{(k)}_i=\left\{\begin{matrix}\frac{a_i}{ a_i^{(k)}}z_i-\frac{1}{ka_i^{(k)}}& \text{ if }&i\notin I_{0}\\
 z_i & \text{ if }&i\in I_{0} .\end{matrix}\right. $$
 Clearly, we have $\lim_{k\longrightarrow \infty}z^{(k)}=z.$ For $k$ sufficiently large $a_i^{(k)}z_i$ is closed to $0$ for all $i\in I_0$. By construction,  we have:
 $$\langle a^{(k)}, z^{(k)}\rangle_{\infty}^-= \bigsmileminus_{i\in [n]}a_{i}^{(k) }z_i^{(k)}= \bigsmileminus_{i\in I_-\cup I_+}\big(a_{i}z_i-\frac{1}{k}\big) \leq c. $$
Let us consider now the sequence $\{w^{(k)}\}_{k\in\mathbb N}$ defined as $w^{(k)}=\frac{c_k}{c}z^{(k)}$. Since $\lim_{k\longrightarrow \infty} c_k=c\not=0$, it follows that $\lim_{k\longrightarrow \infty} w^{(k)}=\lim_{k\longrightarrow \infty} z^{(k)}=z$. Moreover, since for  $k$ sufficiently large we have $\frac{c_k}{c}>0$, the condition $\langle a^{(k)}, z^{(k)}\rangle_{\infty}^-\leq c $ implies that $\langle a^{(k)}, w^{(k)}\rangle_{\infty}^-\leq c_k$. It follows that $z\in  \Linf_{k\longrightarrow \infty}\big[\langle a^{(k )}, \cdot \rangle_{\infty }^-\leq c_k\big ]$. 

$(ii)$ $c=0$. In such a case, if  $\langle a , z \rangle_{\infty}^-<0$ then the result is immediate,  since $\lim_{k\longrightarrow \infty }c_k= 0$ implies that there is some $k_0$ such  that for each $k>k_0$,  $z\in \big[\langle a^{(k )}, \cdot \rangle_{\infty }^-\leq c_k\big ]$.  Suppose now that $\langle a , z \rangle_{\infty}^-=0$. This implies that for all $i\in [n]$ we have $a_iz_i=0$.  Since $a\not=0$,   there is some $i_0$ such that $a_{i_0}\not=0$ and $z_{i_0}=0$, therefore there is some $k_{i_0}$ such that for all $k>k_{i_0}$ $a_{i_0}^{(k)}\not=0. $ 
Moreover, note that since $\bigsmileminus\limits_{i\not=i_0}a_iz_i=0$,  $\lim_{k\longrightarrow \infty}\bigsmileminus\limits_{i\not=i_0}a_i^{(k)}z_i=0$. For each $k>k_{i_0}$, let $z^{(k)}$ be the vector defined as
$$z^{(k)}_i=\left\{\begin{matrix} 
-\big(a_{i}^{(k)}\big)^{-1}  \max\big\{|\bigsmileminus\limits_{i\not=i_0}a_i^{(k)}z_i|,|c_k|\big\} & \text{if }& i=i_0\\z_i& \text{if }& i\not=i_0.
\end{matrix}\right. $$ 
Then, since $-\max\big \{|\bigsmileminus\limits_{i\not=i_0}a_i^{(k)}z_i|,|c_k|\big\}\leq c_k$, we deduce that:
$$\bigsmileminus\limits_{i\in [n]}a_i^{(k)}z_i^{(k)}=\big(-\max\big \{|\bigsmileminus\limits_{i\not=i_0}a_i^{(k)}z_i|,|c_k|\big\}\big)\stackrel{-}{\smile}\big(\bigsmileminus\limits_{i\not=i_0}a_i^{(k)}z_i\big)\leq c_k. $$
Since $\lim_{k\longrightarrow \infty}\bigsmileminus\limits_{i\not=i_0}a_i^{(k)}z_i= \lim_{k\longrightarrow \infty}c_k=0$, it follows that $\lim_{k\longrightarrow \infty}z^{(k)}_{i_0}=0=z_{i_0}$. Since $z^{(k)}_i=z_i$ for all $i\not=i_0$, it follows tha $\lim_{k\longrightarrow \infty}z^{(k)}=z.$ Thus $z\in \Linf_{k\longrightarrow \infty}[\langle  a^{(k)},\cdot \rangle_{\infty}^-\leq c_k]$. Consequently, $[\langle  a ,\cdot \rangle_{\infty}^-\leq c ]\subset \Linf_{k\longrightarrow \infty}[\langle  a^{(k)},\cdot \rangle_{\infty}^-\leq c_k].$ Therefore:
$$\Ls_{k\longrightarrow \infty}[\langle  a^{(k)},\cdot \rangle_{\infty}^-\leq c_k]\subset [\langle  a ,\cdot \rangle_{\infty}^-\leq c ]\subset \Linf_{k\longrightarrow \infty}[\langle  a^{(k)},\cdot \rangle_{\infty}^-\leq c_k],$$
which ends the proof of the first statement. The proof of the second statement for upper forms is symmetric.   $\Box$\\

\subsection{Separation of $\mathbb B$-Convex Sets over $\Real^n$}

The next result states a separation theorem that was mentionned in \cite{b19}, however the method   of proof was not properly mentioned and just outlined from the separation property in limit of fixed $\mathbb B$-convex sets.  We give here a detailled proof. It is straightforwardly obtained from Proposition \ref{FundRes2} and from the standard separation theorem. 

\begin{prop}\label{PresepPoly}Let $A$ and $E$ two finite subsets of   $\Real^n$ such that $Co^\infty(A)\cap Co^\infty(E)=\emptyset$. Then:\\

\noindent $(a)$ There exists a natural number $p_0$ such that for all $p\geq p_0$ 
$$Co^{p }(A )\cap Co^{p}(E)=\emptyset.$$
$(b)$ For all $p\geq p_0$, there exists $a^{(p)}\in \Real^n$ and $c_p\in \Real$ such that:
$$Co^{p }(A )\subset \big [\langle a^{(p)}, \cdot \rangle _p\leq c_p\big ]\quad \text{and}\quad  Co^{p }(E)\subset \big[\langle a^{(p)}, \cdot \rangle _p\geq c_p \big]. $$
$(c)$ There exists $a\in  \Real^n$ and $c\in \Real$ such that:
$$Co^{\infty }(A)\subset \big [\langle a, \cdot \rangle _\infty^{-}\leq c\big ]\quad \text{and}\quad  Co^{\infty }(E)\subset \big[\langle a, \cdot \rangle_\infty^+\geq c \big]. $$

\end{prop} 
{\bf Proof:} $(a)$ Suppose that this is not the case and let us show a contradiction. This implies that 
 there is an increasing sequence of natural numbers $\{p_{q}\}_{q\in {\mathbb N}}$  such that for all $q$ 
$$Co^{(p_{q}) }(A )\cap Co^{(p_{q})}(E ) \not=\emptyset.$$
Therefore, there is a sequence $\{z_{q}\}_{q\in \mathbb N}$ such that $z_{q}\in Co^{(p_{q}) }(A )\cap Co^{(p_{q})}(E )$
for all $q$. Now, since $Co^\infty(A)$ and $Co^\infty(E)$ are compacts,  there exists a compact  $K\subset \Real^n$ that contains the sequence   $\{ z_{q}\}_{q\in \mathbb N}$. Therefore,  there exists a subsequence $\{z_{q_k}\}_{k\in \mathbb N}$ which converges to some   $\bar z \in \Real^n$. It follows that  
$ \bar z \in \Ls_{p\longrightarrow \infty}Co^p(A)=Co^\infty(A)$ and $ \bar z \in \Ls_{p\longrightarrow \infty}Co^p(E)=Co^\infty(E)$.   
However, this implies that $Co^\infty(A)\cap Co^\infty(E)\not=\emptyset$, that is a contradiction, which proves $(a)$. $(b)$ is  is a consequence of the convex separation theorem. $(c)$ Let $C_\infty(0,1)$ denote the circle centered at $0$ of radius $1$. From $(b)$, there exists a sequence $\{a^{(p )}\}_{p\in \mathbb N}\subset C_\infty(0,1)$ and a sequence $\{c_{p}\}_{p\in \mathbb N}$ such that for all $p\geq p_0$:
 $$Co^{p }(A)\subset \big[\langle a^{({p })}, \cdot\rangle_{p }\leq c_{p }\big]\quad \text{ and }\quad Co^{p }(E )\subset \big[\langle a^{({p })}, \cdot\rangle_{p }\geq c_{p }\big].$$
 Now, pick $x_0\in B$ and $y_0\in C$.  We have:
 $$\langle a^{({p })}, x_0\rangle_{p }\leq c_{p_m}\leq \langle a^{({p })}, y_0\rangle_{p }.$$
 Hence there exists a compact interval $I$ of $\Real$ such that $c_{p }\in I$ for all $p$. Therefore, since $C(0,1)\times I$ is compact we can extract an increasing sequence $\{p_{q}\}_{q\in \mathbb N}$ such that there is some $(a,c)\in C(0,1)\times K$ such that
 $\lim_{q\longrightarrow \infty}(a^{(p_{q})},c_{p_{q}})=(a,c).$
Suppose now that $z\in Co^\infty(A)=\Lim_{p\longrightarrow \infty}Co^p(A)$.  This implies that there exists a sequence $\{z^{(q)}\}_{q\in \mathbb N}$ with $z^{(q)}\in \big[\langle a^{({p_q })}, \cdot\rangle_{p_q }\leq c_{p_q }\big]$ for any $q$, and such that $\lim_{q\longrightarrow \infty}z^{(q)}=z$. However, from Proposition \ref{FundRes2}
$$\Lim_{p\longrightarrow \infty}\big[\langle a^{(p )}, \cdot \rangle_{p }\leq c_p\big ]= \big [\langle a, \cdot \rangle _\infty^{-}\leq c\big ].$$
 We deduce that $z\in \big [\langle a, \cdot \rangle _\infty^{-}\leq c\big ]$. Consequently $Co^\infty(A)\subset \big [\langle a, \cdot \rangle _\infty^{-}\leq c\big ] $. The proof of the second part of the statement is symmetric. $\Box$\\

 It was established in \cite{b20} that one can compute the hyperplane limit passing trough a finite number of points in a finite dimensional space. This property provides a useful example of the separation property described above.   Given $n$ points    $v_{1},..., v_n$  in $\Real^n$, let $V$ be the $n\times n$ matrix whose each column is a vector $v_i$. If $|v_1,..., v_n|_p\not=0$ then let $H_p(V)$ denote the $\varphi_p$-hyperplane passing trough $v_1,..., v_n$.  {Let $V$ be the $n\times n$ matrix with  $v_i$ as $i$-th column for each $i$}. Let $V_{(i)}$ be the matrix obtained from $V$ by replacing line $i$ with the transpose of the unit vector $1\!\!1_n$. Suppose that $|V|_\infty\not=0$. Then
\begin{equation}H_\infty(V):=\Lim_{p\longrightarrow \infty}H_p(V)=\Big\{x\in \Real^n: \bigsmileminus_{i\in [n]} |V_{(i)}|_\infty x_i\leq |V|_\infty\leq \bigsmileplus_{i\in [n]} |V_{(i)}|_\infty x_i\Big\}.\end{equation}
 In the situation where we consider $n$ points with 
$A = \{x^{(1)}, \dots, x^{(n)}\} \subset \mathbb{R}^n$ 
and $v_j = x^{(j)}$ for all $j$, we should have 
$ {Co}^p(A) \subset H_p(V)$ and this implies that $
 {Co}^\infty(A) \subset H_\infty(V).$
 In the following, returning to example \ref{fundex}, we can verify 
that all the intermediate points belong to the limiting hyperplane 
that contains $x^{(1)}, x^{(2)}$, and $x^{(3)}$.   Therefore, since example \ref{fundex}  showed that $(0,0,0)\notin Co^\infty(A)$, it is shown below that $\{(0,0,0)\}$ can be separated from $Co^\infty(A)$.
 
 \begin{expl}\label{hypcontains}{\small In example \ref{fundex} we consider the vectors $x^{(1)}=(1,3,3)$, $x^{(2)}=(3,2,3)$ and $x^{(3)}=(-1,-1,-1)$. We have
$$V=\begin{pmatrix}1&3&-1\\3&2&-1\\3&3&-1\end{pmatrix}, V_{(1)}=\begin{pmatrix}1&1&1\\3&2&-1\\3&3&-1\end{pmatrix}, V_{(2)}=\begin{pmatrix}1&3&-1\\1&1&1\\3&3&-1\end{pmatrix} \text{ and }\; V_{(3)}=\begin{pmatrix}1&3&-1\\3&2&-1\\1&1&1\end{pmatrix}. $$
Hence:

 $|V|_\infty=1\cdot 2\cdot (-1)\boxplus 3\cdot 3\cdot (-1)\boxplus 3\cdot (-1)\cdot 3 \boxplus (-1)\cdot(-1)\cdot 2\cdot (-1)\boxplus (-1)\cdot(-1)3\cdot 3\cdot (-1) \boxplus (-1)\cdot(-1)\cdot 3\cdot 1=-9$; 

$|V_{(1)}|_\infty=1\cdot 2\cdot  (-1)\boxplus 3\cdot 3\cdot 1\boxplus 1\cdot (-1) \cdot 3\boxplus (-1)\cdot 1\cdot 2\cdot 3 )\boxplus (-1)\cdot (-1)\cdot 3\cdot 1 \boxplus (-1)\cdot 3\cdot 1\cdot (-1)=9$; 

$|V_{(2)}|_\infty=1\cdot 1\cdot  (-1) \boxplus  3\cdot 1\cdot  3 \boxplus 1\cdot 3\cdot (-1) \boxplus (-1)\cdot (-1)\cdot 1\cdot 3\boxplus (-1)\cdot 1\cdot 3\cdot 1 \boxplus (-1)\cdot (-1)\cdot  3\cdot (-1)=9$; 

$|V_{(3)}|_\infty=1\cdot 2\cdot    1\boxplus 3\cdot (-1)\cdot 1\boxplus 3\cdot (-1)\cdot 1 \boxplus (-1)\cdot 1\cdot 2\cdot 1\boxplus (-1)\cdot (-1)\cdot 1\cdot 1 \boxplus (-1) 3\cdot 3\cdot 1=-9$.

$$H_\infty(V)=\big\{(x_1,x_2,x_3): 9x_1\stackrel{-}{\smile} 9 x_2 \stackrel{-}{\smile} (-9) x_3\leq -9\leq 9x_1\stackrel{+}{\smile} 9 x_2 \stackrel{+}{\smile} (-9) x_3\big \}.$$

Equivalently:

$$H_\infty(V)=\big\{(x_1,x_2,x_3):  x_1\stackrel{-}{\smile}   x_2 \stackrel{-}{\smile} (-1) x_3\leq -1\leq  x_1\stackrel{+}{\smile}  x_2 \stackrel{+}{\smile} (-1) x_3\big \}.$$

Let us check that all the initial points and all the intermediate points belong to the limit hyperplane.

$(i)$ $x^{(1)}=(1,3,3)$. We have $1\stackrel{-}{\smile}3\stackrel{-}{\smile}(-3)=-3\leq -1$ and $1\stackrel{+}{\smile}3\stackrel{+}{\smile}(-3)=3\geq -1$, which implies that $x^{(1)}\in H_\infty(V)$.
 
$(ii)$  $x^{(2)}=(3,2,3)$. We have $3\stackrel{-}{\smile}3\stackrel{-}{\smile}(-3)=-3\leq -1$ and $3\stackrel{+}{\smile}3\stackrel{+}{\smile}(-3)=3\geq -1$. Hence $x^{(2)}\in H_\infty(V)$.
  
$(iii)$    $x^{(3)}=(-1,-1,-1)$. We have $(-1)\stackrel{-}{\smile}(-1)\stackrel{-}{\smile}(1)=-1\leq -1$ and $(-1)\stackrel{+}{\smile}(-1)\stackrel{+}{\smile}(1)=1\geq -1$. Therefore  $x^{(3)}\in H_\infty(V)$.

$(iv)$ $\zeta_{\{1\},\{1,3\}}=(
0,3,3)$.  We have $0\stackrel{-}{\smile}3\stackrel{-}{\smile}(-3)=-3\leq -1$ and $0\stackrel{+}{\smile}3\stackrel{+}{\smile}(-3)=3\geq -1$. Thus  $\zeta_{\{1\},\{1,3\}}\in H_\infty(V)$.

$(v)$ $\zeta_{\{2\},\{1,3\}}= \zeta_{\{3\},\{1,3\}}=(
-1,0,0).$ We have $(-1)\stackrel{-}{\smile}0\stackrel{-}{\smile}0=-1\leq -1$ and $(-1)\stackrel{+}{\smile}0\stackrel{+}{\smile}0=-1\geq -1$. Hence  $\zeta_{\{2\},\{1,3\}}= \zeta_{\{3\},\{1,3\}}\in H_\infty(V)$.
 
$(vi)$ $ \zeta_{\{2\},\{2,3\}}=(
\frac{3}{2},0,\frac{3}{2})$. We have $\frac{3}{2}\stackrel{-}{\smile}0\stackrel{-}{\smile}(-\frac{3}{2})=- \frac{3}{2}\leq -1$ and $\frac{3}{2}\stackrel{+}{\smile}0\stackrel{+}{\smile}(-\frac{3}{2})= \frac{3}{2}\geq -1$. Hence  $\zeta_{\{2\},\{2,3\}}\in H_\infty(V)$.

$(vii)$ $ \zeta_{\{3\},\{2,3\}}=\zeta_{\{1,3\},\{1,2,3\}}=\zeta_{\{2,3\},\{1,2,3\}}=(0,-1,0).$ We have $0\stackrel{-}{\smile}(-1)\stackrel{-}{\smile}0=- 1\leq -1$ and $0\stackrel{+}{\smile}(-1)\stackrel{+}{\smile}0= -1\geq -1$. Hence  $\zeta_{\{3\},\{2,3\}}=\zeta_{\{1,3\},\{1,2,3\}}=\zeta_{\{2,3\},\{1,2,3\}}\in H_\infty(V)$.

$(viii)$   $\zeta_{\{1,2\},\{1,2,3\}}=(0,0,1).$ We have $0\stackrel{-}{\smile}0\stackrel{-}{\smile}(-1)=- 1\leq -1$ and $0\stackrel{+}{\smile}0\stackrel{+}{\smile}(-1)= -1\geq -1$. Hence  $\zeta_{\{1,2\},\{1,2,3\}}\in H_\infty(V)$.}

We   note that point $(0,0,0)$ does not belong to the hyperplane $H_\infty(V)$ since 
\begin{equation}0\stackrel{+}{\smile}0\stackrel{+}{\smile}(-0)=0\stackrel{-}{\smile}0\stackrel{-}{\smile}(-0)=0>-1.\end{equation}
Let us denote $a=(1,1,-1)$ and $c=-1$. It follows that 
$$Co^\infty(A)\subset  \big [\langle a, \cdot \rangle_{ \infty}^-\leq c\big]\quad \text{ and }\quad \{(0,0,0)\}\subset  \big[\langle a, \cdot \rangle_{ \infty}^+\geq c\big ]. $$This is coherent with example \ref{fundex},  since we showed that $(0,0,0)\notin Co^\infty(A)$ and $H_\infty(V)$ must contain $Co^\infty(A)$. 

\end{expl}

 \begin{lem}\label{approx}Let $C$ be a $\mathbb B$-convex compact subset of $\Real^n$. There exists a sequence of finite subsets of $C$, $\{A_{m}\}_{m\in \mathbb N}$, with $A_{m+1}\subset A_m$ for any positive integer $m$, such that
 $$\Lim_{m\longrightarrow \infty}Co^\infty(A_m)=C. $$
 \end{lem}
 {\bf Proof:} Let $ C \subset \mathbb{R}^n $ be a compact set. From the compactness of \( C \), which implies that it can be covered by finitely many \( \|\cdot\|_\infty \)-balls of arbitrarily small radius, centered at points of \( C \), for every positive integer $ m \in \mathbb{N} $, there exists a finite subset $ A_m \subset C $ such that:
$$
d_{H,\infty}(A_m, C) \leq \frac{1}{m},
$$
where $d_{H,\infty} $ denotes the Hausdorff{-Pompeiu} distance induced by the $ \|\cdot\|_\infty $ norm. One can then construct sequence  $\{A_m\}_{m\in \mathbb N}$ a sequence of finite subsets of points of $C$, such that
 $A_{m}\subset A_{m+1}$ 
and $d_{H,\infty}( A_m, C)\leq \frac{1}{m}$ for any positive integer $m$.
 Since $C $ is $\mathbb B$-convex, for any $m$, $Co^\infty(A_m)$ is a subset of $C $. Therefore, since $C\supset Co^\infty(A_m)\supset A_m$, we deduce that:
 $$d_{H,\infty}(Co^\infty(A_m), C )\leq d_{H,\infty}( A_m, C)\leq \frac{1}{m}.$$
  It follows that $\lim_{m\longrightarrow \infty}d_{H,\infty}(Co^\infty(A_m), C)=0$. This implies that $$\Lim_{n\longrightarrow \infty}Co^\infty(A_m)=C. \quad \Box$$

For any $x\in \Real^n$  , let $B_\infty(x , d\,]$ denote the closed ball centered at $x$
of radius $d$, with respect to the $ \ell_\infty$-norm. 

\begin{prop}Let $C$ and $D$ two   non-empty closed  $\mathbb B$-convex subsets of $\Real^n$  with $C\cap D=\emptyset$. Then there exists an idempotent symmetric form $f:\Real^n\longrightarrow \Real$ and a real number $c$, such that:
$$C\subset[f^-\leq c]\quad \text{ and }\quad D\subset[f^+\geq c]$$
where $f^-$ and $f^+$ are respectively the lower and upper semi-continuous regularized of $f. $ 
\end{prop}
{\bf Proof:} We first suppose that $C$ and $D$ are compacts. From Lemma \ref{approx} there     exists a sequence of finite subsets of $C$, $\{A_{m}\}_{m\in \mathbb N}$, with $A_{m+1}\subset A_m$ for any positive integer $m$, such that $$\Lim_{n\longrightarrow \infty}Co^\infty(A_m)=C. $$
Similarly, there exists an increasing sequence of finite subsets of $D $ denoted $\{E_m\}_{m\in \mathbb N}$ such that
$$\Lim_{n\longrightarrow \infty}Co^\infty(E_m)=D. $$

Let $C_\infty(0,1)$ denote the circle centered at $0$ of radius $1$. Now,   since $C \cap D=\emptyset$, there exists $m_0\in \mathbb N$ such that for all $m\geq m_0$ $Co^\infty(A_m)\cap Co^\infty(E_m)=\emptyset.$ Moreover, from Proposition \ref{PresepPoly}  there exists for each $m$, $a ^{( m)}\subset C_\infty(0,1)$ and $ c_{m} $ such that
$$Co^\infty(A_m) \subset [\langle a^{(m)}, \cdot \rangle_\infty^-\leq c_{ m}]\quad \text{ and }\quad Co^\infty(E_m) \subset [\langle a^{(m)}, \cdot \rangle_\infty^+\geq c_{m}]. $$
Pick   $x_0\in A_0$ and $y_0\in  E_0 $. Now, note that since $x_0\in C$ and $y_0\in D$ we have:
 $$-\|x_0\|_\infty\leq  \langle a^{(m)} , x_1\rangle_{\infty}^-\leq c_{m }\leq \langle a^{(m)}, y_1\rangle_{\infty}^+\leq  \|y_0\|_\infty.$$
 Hence there exists a compact interval $I$ of $\Real$ such that $c_{m}\in I$ for each $m$. Now, let us consider the sequence $\{(a^{(m)},c_m)\}_{m\in \mathbb N}$. Therefore, since $C(0,1)\times I$ is compact we can extract an increasing sequence $\{m_q\}_{q\in \mathbb N}$ such that there is some $(a ,c )\in C(0,1)\times K$ such that
 $\lim_{\longrightarrow \infty}(a^{(m_q)},c_{m_q})=(a ,c ).$
 It follows taking the limit that 
 $$C=\Lim_{m\longrightarrow \infty} Co^\infty(A_m)\Big)\subset [\langle a , \cdot \rangle_\infty^-\leq c ]\quad \text{ and }\quad D=\Lim_{m\longrightarrow \infty}Co^\infty(E_m))\subset [\langle a , \cdot \rangle_\infty^-\geq c ]. $$
 
 Let us prove the case where $C$ and $D$ are not compact. Let $x_0\in C$ and let us consider the ball $B_\infty(x_0, k\,]$. Let $k_0$ be a positive natural number such that for all $k\geq k_0$ we have $B_\infty(x_0, k\,]\cap D\not=\emptyset$. Let us denote 
$$C_k=C\cap B_\infty(x_0, k\,]\quad \text{ and } \quad D_k=D\cap B_\infty(x_0, k].  $$
Suppose now that for all $k \geq k_0$, $y_0\in D_k.$ First note that for all $k$,  $B_\infty(x_0, k]$ is compact and $\mathbb B$-convex set. It follows that $C_k$ and $D_k$ are compact and $\mathbb B$-convex if $k\geq k_0$. Moreover since $C$ and $D$ are non proximate, $C_k\cap D_k=\emptyset $ for each $k$. 
Then fore, for each $k$, there exists   $(a^{(k)}, c_k)\in C_\infty(0,1)\times [-\|x_0\|_\infty,\|y_0\|_\infty ]$ for each $k$, such that:
$$C_k  \subset [\langle a^{(k)} , \cdot \rangle_\infty^-\leq c_k ]\quad \text{ and }\quad D_k \subset [\langle a^{(k)} , \cdot \rangle_\infty^-\geq c_k ]. $$
Now, note that 
$$C=\Lim_{k\longrightarrow \infty}C_k\quad \text{ and }\quad D=\Lim_{k\longrightarrow \infty}D_k. $$
Using similar arguments, there exists a subsequence $\{(a^{(k_l)}, c_{k_l})\}_{l\in \mathbb N}$ which converges to some $(a,c)\in C_\infty(0,1)\times [-\|x_0\|, \|y_0\| ]$. 
Therefore we deduce that:
$$C\subset [\langle a , \cdot \rangle_\infty^-\leq c ]\quad \text{ and }\quad D \subset [\langle a , \cdot \rangle_\infty^-\geq c]  $$
which ends the proof. $\Box$\\

  \begin{prop}  Let $C$ be a closed  $\mathbb B$-convex set. $C$ is the intersection of all the lower halfspaces which contain it. 
  \end{prop}
{\bf Proof:} Suppose that $x\notin C$. There exists some $\epsilon >0$  such that $B_\infty(x, \epsilon]\cap B=\emptyset$. Since $B_\infty(x, \epsilon]$ is closed and $\mathbb B$-convex, there  exists some $a\in \Real^n$ and some  $\alpha\in \Real$ such that:
$$C\subset [\langle a, \cdot \rangle_\infty^-  \leq \alpha ] \quad \text{ and }\quad B_\infty(x, \epsilon]\subset  [\langle a, \cdot \rangle_\infty^+\geq \alpha].$$ 
It follows that $x\in \mathrm{int}[\langle a, \cdot \rangle_\infty^+\geq \alpha].$ This implies that $x\notin [\langle a, \cdot \rangle_\infty^-  \leq \alpha ]$ which ends the proof. $\Box$\\

In \cite{b19}, it was claimed that an external representation of a $\mathbb{B}$-polytope can be given from a finite collection of lower halfspaces. In the following, we slightly correct this result, which only holds in the case where the intersection set is topologically regular. {A closed subset $E \subset \mathbb{R}^n$ is said to be {\bf topologically regular} if
$E = \mathrm{cl}\big({\operatorname{int}(E)}\big).$
In other words, $E$ coincides with the closure of its interior.}

\begin{prop} For any finite subset $A$ of $\Real^n$, there exists $\ell$ lower idempotent symmetric forms $f_1^-,\cdots,f_\ell^-$ and $\ell $ real numbers $c$, $j=1,\cdots, \ell$ such that
$$\mathrm{int}\Big(  \bigcap_{k\in [\ell]}[f_j^-\leq c_j]\Big)\subset Co^\infty(A)\subset \bigcap_{k\in [\ell]}[f_j^-\leq c_j]  .$$
If the intersection set   $\bigcap\limits_{k\in [\ell]}[f_j^-\leq c_j]$ is topologically regular then 
$$Co^\infty(A)= \bigcap_{k\in [\ell]}[f_j^-\leq c_j]  .$$
\end{prop}
{\bf Proof:}
From \cite{b19}, there exists an
increasing subsequence $\{p_q\}_{q\in\mathbb{N}}$ and $\ell $
subsequences of linear forms $\{f_{j,p_q}\}_{q\in \mathbb{N}} $
respectively defined for all natural numbers $q$ by
$f_{j,p_q}(x)=\langle a_{j}^{(p_q)},x\rangle_{p_q}$ with $
{{Co}}^{p_q}(A)= \bigcap_{j\in [\ell ]} {\big [f_{j,p_q}\leq
c_{j,p_q}\big ]}$. Hence for each $j$:
\begin{align*}c_{j,p_q}&
=\sup_{x}\left\{\langle a_{j}^{(p_q)},x\rangle_{p_{q}} : x\in
{Co}^{{p_q}}(A)\right\}=\max_{k=1,...,m}\langle
a_{k}^{(p_q)},x^{(k)}\rangle_{p_q}.\end{align*} For all vectors $u$
of $\Real^n$, let us denote $|u|=(|u_1|,...,|u_n|)$. Fix $\bar
c=n\max_{k=1,..., m}\|x^{(k)}\|_\infty$. For $j=1,..., \ell $, we
have from the H\"{o}lder inequality:
$$ -\bar c
\leq -\max_{k=1,...,m}\langle |a_j|,|x^{(k)}|\rangle_{p_{q}}\leq
c_{j,p_q}\leq \max_{k=1,...,m}\langle
|a_j|,|x^{(k)}|\rangle_{p_{q}}\leq \bar c.$$
 Since $C_{\infty}(0,1)\times [-\bar c,
 \bar c]$ is compact,  $\Big (C_{\infty}(0,1)\times [-\bar c,
 \bar c]\Big )^{\ell }$ is compact and we deduce that it contains a subsequence $$\Big\{\Big((a_{1}^{(p_{q_r})},c_{1,p_{q_r}}),..., (a_{\ell  }^{(p_{q_r})},c_{\ell ,p_{q_r}})
\Big)\Big\}_ {r\in \mathbb{N}}$$ which converges to some
$\big((a_{1},c_{1}),..., (a_{\ell },c_{\ell })\big)\in
\Big(C_{\infty}(0,1)\times [-\bar c,
 \bar c]\Big)^{\ell }$.  Let us denote $f_j: x\mapsto \langle a^{(j)}, x\rangle_\infty.$ Moreover, let $f_j^-$  and $f_j^+$ be respectively the lower and upper semi-continuous regularized of $f_j$. Hence
\begin{align*}\Linf_{r\longrightarrow \infty} {Co}^{{p_{q_r}}}(A) =
\Linf_{r\longrightarrow \infty}\bigcap_{j\in [ \ell ]}[f_{j,p_{q_r}}
\leq c_{j,p_{q_r}}] \subset \bigcap_{j\in [ \ell ]}\Linf_{r\longrightarrow \infty}
[\ f_{j, p_{q_r}}\leq c_{j,p_{q_r}}] \subset \bigcap_{j\in [ \ell ]} \cl{[f_j\leq c_{j}]}=\bigcap_{j\in [ \ell ]} {[f_j^-\leq c_{j}]}.
\end{align*}
  However, by definition
  $\Linf_{p\longrightarrow \infty}Co^p(A)
  \subset \Linf_{r\longrightarrow \infty} {Co}^{{p_{q_r}}}(A)$. Since
  we have established that $Co^\infty (A)=\Lim_{p\longrightarrow \infty}Co^p(A)=
  \Linf_{p\longrightarrow \infty}Co^p(A)$, this yields the  inclusion.

Let us show the converse inclusion.
Assume that $y\notin 
{Co}^{\infty}(A)$. In such a case there is an increasing sequence of natural numbers $\{p_k\}_{k\in \mathbb N}$ and some $j$
such that for all $k$, $y\in [f_{j, p_k}>c_{j,p_k}]$. It follows from Proposition \ref{FundRes2} that $$y\in \Lim_{p\longrightarrow \infty}[f_{j, p}>c_{j,p}]\subset \Lim_{p\longrightarrow \infty}[f_{j, p}\geq c_{j,p}]= [f_j^+\geq c_j].$$  Therefore $y\in \Real^n\backslash \mathrm{int}[f_j^-\leq c_j]$. Hence $y\notin \bigcap\limits_{j\in [\ell]}\mathrm{int}[f_j^-\leq c_j]=\mathrm{int}\bigcap_{j\in [\ell]}[f_j^-\leq c_j]$, which proves the second part of the statement.  To complete the proof,   if $\bigcap\limits_{k\in [\ell]}[f_j^-\leq c_j]$ is topologically regular, then  we deduce that:
$$\bigcap\limits_{k\in [\ell]}[f_j^-\leq c_j]=\cl\Big(\mathrm{int}\Big(  \bigcap_{k\in [\ell]}[f_j^-\leq c_j]\Big)\Big)\subset \cl \Big(\mathrm{int}\Big(Co^\infty(A)\Big)\Big)\subset {Co}^{\infty}(A), $$
which yields the result.    $\Box$\\

In the following, we consider the example of a tetrahedron in limit by adding the point $x^{(0)}=(0,0,0)$ to the triple $A=\{x^{(1)},x^{(2)},x^{(3)}\}$ considered in Example \ref{hypcontains}. 

\begin{expl}{\small Let us consider the set $B=\{x^{(0)},x^{(1)},x^{(2)},x^{(3)}\}$. We want to the hyperplanes support of $T^\infty=Co^\infty(B)$ and provide an external representation of $T^\infty$. To do that we need to find the hyperplan in limit for the triples $\{x^{(0)}x^{(1)},x^{(2)}\} $, $ \{x^{(0)}x^{(2)},x^{(3)}\} $ and $ \{x^{(0)}x^{(1)},x^{(3)}\} $.

$(i)$ $\{x^{(0)}x^{(1)},x^{(2)}\} $. We have the matrices:

$$V'=\begin{pmatrix}0&1&3 \\0&3&2 \\0&3&3 \end{pmatrix}, V'_{(1)}=\begin{pmatrix}1&1&1\\0&3&2 \\0&3&3 \end{pmatrix}, V'_{(2)}=\begin{pmatrix}0&1&3 \\1&1&1\\0&3&3 \end{pmatrix} \text{ and }\; V'_{(3)}=\begin{pmatrix}0&1&3 \\0&3&2 \\1&1&1\end{pmatrix}. $$
Hence:

 $|V'|_\infty=0$; $|V'_{(1)}|_\infty= 9$; $|V'_{(2)}|_\infty= 9$;  $|V'_{(3)}|_\infty=-9$.

$$H_\infty(V')=\big\{(x_1,x_2,x_3): 9x_1\stackrel{-}{\smile} 9 x_2 \stackrel{-}{\smile} (-9) x_3\leq 0\leq 9x_1\stackrel{+}{\smile} 9 x_2 \stackrel{+}{\smile} (-9) x_3\big \}.$$
Equivalently:

$$H_\infty(V)=\big\{(x_1,x_2,x_3):  x_1\stackrel{-}{\smile}   x_2 \stackrel{-}{\smile} (-1) x_3\leq 0\leq  x_1\stackrel{+}{\smile}  x_2 \stackrel{+}{\smile} (-1) x_3\big \}.$$

$(ii)$ $\{x^{(0)}x^{(2)},x^{(3)}\} $. We have the matrices:

We have
$$V''=\begin{pmatrix}0&3&-1\\0&2&-1\\0&3&-1\end{pmatrix}, V''_{(1)}=\begin{pmatrix}1&1&1\\0&2&-1\\0&3&-1\end{pmatrix}, V''_{(2)}=\begin{pmatrix}0&3&-1\\1&1&1\\0&3&-1\end{pmatrix} \text{ and }\; V''_{(3)}=\begin{pmatrix}0&3&-1\\0&2&-1\\1&1&1\end{pmatrix}. $$
Hence: $|V''|_\infty= 0$; $|V''_{(1)}|_\infty=3$; $|V''_{(2)}|_\infty=0$; $|V''_{(3)}|_\infty=-3$.
Thus

$$H_\infty(V'')=\big\{(x_1,x_2,x_3): 3x_1\stackrel{-}{\smile}  (-3) x_3\leq 0\leq 3x_1\stackrel{+}{\smile}  (-3)x_3;, (- x_1)\stackrel{-}{\smile} (-  x_2) \stackrel{-}{\smile}   x_3\leq 0, \big \}.$$

Equivalently:

$$H_\infty(V'')=\big\{(x_1,x_2,x_3):  x_1\stackrel{-}{\smile}  (- x_3)\leq 0\leq  x_1\stackrel{+}{\smile}  (- x_3)\big \}.$$

$(iii)$ $\{x^{(0)}x^{(1)},x^{(3)}\} $. We have the matrices:

We have
$$V'''=\begin{pmatrix}0&1& -1\\0&3& -1\\0&3& -1\end{pmatrix}, V'''_{(1)}=\begin{pmatrix}1&1&1\\0&3& -1\\0&3& -1\end{pmatrix}, V'''_{(2)}=\begin{pmatrix}0&1& -1\\1&1&1\\0&3& -1\end{pmatrix} \text{ and }\; V'''_{(3)}=\begin{pmatrix}0&1& -1\\0&3& -1\\1&1&1\end{pmatrix}. $$
Hence: $|V'''|_\infty=0$; $|V'''_{(1)}|_\infty=0$; $|V'''_{(2)}|_\infty= -3$; $|V'''_{(3)}|_\infty=  3$.

$$H_\infty(V''')=\big\{(x_1,x_2,x_3): (-3) x_2 \stackrel{-}{\smile} 3 x_3\leq 0\leq  (-3) x_2 \stackrel{+}{\smile} 3 x_3\big \}.$$

Equivalently:

$$H_\infty(V''')=\big\{(x_1,x_2,x_3):   (-   x_2 )\stackrel{-}{\smile}   x_3\leq 0\leq  (- x_2) \stackrel{+}{\smile}  x_3\big \}.$$

Now  note that  $-1\leq  x_1\stackrel{+}{\smile}  x_2 \stackrel{+}{\smile} (-1) x_3$ if and only if $(- x_1)\stackrel{-}{\smile}  (-x_2) \stackrel{-}{\smile}   x_3\leq 1$ which ensures that $x^{(0)}=(0,0,0)$ satisfies this last inequality. Moreover $0\leq  (- x_2) \stackrel{+}{\smile}  x_3$ if and only if $ x_2  \stackrel{-}{\smile}  (-x_3)\leq 0$ ensuring that all the points in $Co^\infty(A)$ satify this inequality. We deduce that:
$$T^\infty=Co^\infty(B)\subset\{(x_1,x_2,x_3):(- x_1)\stackrel{-}{\smile}  (-x_2) \stackrel{-}{\smile}   x_3\leq 1,  x_1 \stackrel{-}{\smile}   x_2  \stackrel{-}{\smile}   (-x_3)\leq 0,x_1\stackrel{-}{\smile}  (- x_3)\leq 0, x_2  \stackrel{-}{\smile}  (-x_3)\leq 0 \}.$$}

\end{expl}

% ------------------------------------------------------------------------

% ------------------------------------------------------------------------


\begin{thebibliography}{99}

\bibitem{adil}
{\sc G. Adilov and A. M. Rubinov},
$\mathbb{B}$-convex sets and functions,
{\it Numer. Funct. Anal. Optim.} {\bf 27} (2006), 237--257.

\bibitem{adilYe}
{\sc G. Adilov and I. Ye\c{s}ilce},
$\mathbb{B}^{-1}$-convex functions,
{\it J. Convex Anal.} {\bf 24} (2017), 505--517.

\bibitem{avr1}
{\sc M. Avriel},
$R$-convex functions,
{\it Math. Programming} {\bf 2} (1972), 309--323.

\bibitem{avr2}
{\sc M. Avriel},
{\it Nonlinear Programming: Analysis and Methods},
Prentice Hall, New Jersey, 1976.

\bibitem{b62}
{\sc A. Beck},
A convexity condition in Banach spaces and the strong law of large numbers,
{\it Proc. Amer. Math. Soc.} {\bf 13} (1962), 329--334.

\bibitem{ben}
{\sc A. Ben-Tal},
On generalized means and generalized convex functions,
{\it J. Optim. Theory Appl.} {\bf 21} (1977), 1--13.

\bibitem{b10}
{\sc P. Butkovi\v{c}},
{\it Max-linear Systems: Theory and Algorithms},
Springer, 2010.

\bibitem{butheged84}
{\sc P. Butkovi\v{c} and G. Heged\"us},
An elimination method for finding all solutions of the system of linear equations over an extremal algebra,
{\it Ekonom.-Mat. Obzor} {\bf 20} (1984).

\bibitem{b15}
{\sc W. Briec},
Some remarks on an idempotent and non-associative convex structure,
{\it J. Convex Anal.} {\bf 22} (2015), 259--289.

\bibitem{b17}
{\sc W. Briec},
Separation properties in some idempotent and symmetrical convex structures,
{\it J. Convex Anal.} {\bf 24} (2017), 1143--1168.

\bibitem{b19}
{\sc W. Briec},
On some class of polytopes in an idempotent, symmetrical and non-associative convex structure,
{\it J. Convex Anal.} {\bf 26} (2019), 823--853.

\bibitem{b20}
{\sc W. Briec},
Determinants and limit systems in some idempotent and non-associative algebraic structures,
{\it Linear Algebra Appl.} {\bf 651} (2022), 162--208.

\bibitem{b24}
{\sc W. Briec},
Remarks on some limit geometric properties related to an idempotent and non-associative algebraic structure,
{\it J. Math. Sci.} (2025). https://doi.org/10.1007/s10958-025-07795-0.

\bibitem{bh}
{\sc W. Briec and C. D. Horvath},
$\mathbb{B}$-convexity,
{\it Optimization} {\bf 53} (2004), 103--127.

\bibitem{bh2}
{\sc W. Briec and C. D. Horvath},
Nash points, Ky Fan inequality and equilibria of abstract economies in max-plus and $\mathbb{B}$-convexity,
{\it J. Math. Anal. Appl.} {\bf 341} (2008), 188--199.

\bibitem{bhr05}
{\sc W. Briec, C. D. Horvath and A. Rubinov},
Separation in $\mathbb{B}$-convexity,
{\it Pacific J. Optim.} {\bf 1} (2005), 13--30.

\bibitem{bh3}
{\sc W. Briec and C. D. Horvath},
On the separation of convex sets in some idempotent semimodules,
{\it Linear Algebra Appl.} {\bf 435} (2011), 1542--1548.

\bibitem{by22}
{\sc W. Briec and I. Yesilce},
Ky Fan inequality, Nash equilibria in some idempotent and harmonic convex structures,
{\it J. Math. Anal. Appl.} {\bf 508} (2022).

\bibitem{hlp}
{\sc G. H. Hardy, J. E. Littlewood and G. P\'olya},
{\it Inequalities},
Cambridge Univ. Press, 1952.

\bibitem{ir10}
{\sc Z. Izhakian and L. Rowen},
Supertropical algebra,
{\it Adv. Math.} {\bf 225} (2010), 2222--2286.

\bibitem{kya15}
{\sc S. Kemali, I. Ye\c{s}ilce and G. Adilov},
$\mathbb{B}$-convexity, $\mathbb{B}^{-1}$-convexity, and their comparison,
{\it Numer. Funct. Anal. Optim.} {\bf 36} (2015), 133--146.

\bibitem{F86}
{\sc S. Friedland},
Limit eigenvalues of nonnegative matrices,
{\it Linear Algebra Appl.} {\bf 74} (1986), 173--178.

\bibitem{koloma}
{\sc V. N. Kolokoltsov and V. P. Maslov},
{\it Idempotent Analysis and its Applications},
Kluwer Academic Publishers, 1997.

\bibitem{ms}
{\sc V. P. Maslov and S. N. Samborski} (eds.),
{\it Idempotent Analysis},
Amer. Math. Soc., Providence, 1992.

\bibitem{my24}
{\sc C. Mo and Y. Yang},
The unified description of abstract convexity structures,
{\it Axioms} {\bf 13} (2024), 506.

\bibitem{p94}
{\sc J. Pin},
Tropical semirings,
in {\it Idempotency} (Bristol, 1994), Cambridge Univ. Press, 1998, 50--69.

\bibitem{MP}
{\sc M. Plus},
Linear systems in (max,+)-algebra,
Proc. 29th Conf. Decision and Control, Honolulu, 1990.

\bibitem{rw98}
{\sc R. T. Rockafellar and R. J.-B. Wets},
{\it Variational Analysis},
Springer, 1998.

\bibitem{rubi}
{\sc A. Rubinov},
{\it Abstract Convexity and Global Optimization},
Kluwer, 2000.

\bibitem{s88}
{\sc I. Simon},
Recognizable sets with multiplicities in the tropical semiring,
in {\it MFCS 1988}, Lecture Notes in Comput. Sci. {\bf 324}, Springer, 1988, 107--120.

\bibitem{t23}
{\sc L. Thibault},
{\it Unilateral Variational Analysis in Banach Spaces}, Vol.~I--II,
World Scientific, 2023.

\bibitem{adilYeTi}
{\sc G. Tinaztepe, I. Ye\c{s}ilce and G. Adilov},
Separation of $\mathbb{B}^{-1}$-convex sets by $\mathbb{B}^{-1}$-measurable maps,
{\it J. Convex Anal.} {\bf 21} (2014), 571--580.

\bibitem{ut24}
{\sc M. Uzun and T. Tun\c{c}},
On some positive linear operators preserving the $\mathbb{B}^{-1}$-convexity of functions,
{\it Fundam. Contemp. Math. Sci.} {\bf 5} (2024), 134--142.

\bibitem{vdv}
{\sc M. van de Vel},
{\it Theory of Convex Structures},
North-Holland, 1993.

\bibitem{iy23}
{\sc I. Ye\c{s}ilce, G. Tinaztepe, S. Kemali and G. Adilov},
Inequalities involving general fractional integrals of $p$-convex functions,
{\it Turkish J. Math.} {\bf 47} (2023), 2028--2042.

\end{thebibliography}
\end{document}